# THÈSE

présentée à l'Université de Provence
pour obtenir le titre de
Docteur de l'université
Spécialité : Mathématiques pures
par

François Ziegler

# MÉTHODE DES ORBITES
# ET
# REPRÉSENTATIONS QUANTIQUES

Soutenue le 20 décembre 1996
devant le jury composé de Messieurs :

| | |
|---|---|
| Jean-Philippe Anker | Jimmy Elhadad |
| Patrick Delorme | Étienne Ghys |
| Michel Duflo | Bertram Kostant |
| Christian Duval | Jean-Marie Souriau, directeur |

Rapporteurs :

| | |
|---|---|
| Michel Duflo | Victor W. Guillemin |

*Au lecteur inconnu*

# Table

> *Puis il a été reçu docteur. Là, le comique du sérieux a commencé, pour faire suite au sérieux du comique qui avait précédé. Il est devenu grave, s'est caché pour faire de minces fredaines, s'est acheté définitivement une montre et a renoncé à l'imagination (textuel); comme la séparation a dû être pénible!*
> —G. Flaubert

# Remerciements

# QUANTUM REPRESENTATIONS
# AND THE ORBIT METHOD

François Ziegler
CNRS-CPT Luminy, Case 907
13288 Marseille Cedex 09
ziegler@cpt.univ-mrs.fr

**Abstract.** J. M. Souriau has suggested an *a priori* characterization of the unitary representations 'attached' to a given coadjoint orbit of a Lie group. When the group is compact, we show that his condition effectively selects the expected representation within sections of the line bundle over the orbit. When the group is noncompact we find many unexpected representations, but we show how Souriau's condition can be refined so as to either recover the traditional theory for exponential groups, or characterize some new discontinuous representations in which states can be 'localized' on lagrangian submanifolds of the orbit.

Contents

# 0. Introduction

THE IDEA of a correspondence between the unitary representations of a Lie group, $G$, and the symplectic manifolds on which $G$ acts, can be traced to the early days of quantum mechanics ('correspondence principle'). As formalized by Kirillov, Kostant, and Souriau during the 1960s, it has brought geometrical understanding into much of harmonic analysis. Thus the following two-step scheme for producing a representation from a symplectic $G$-manifold $X$ has become familiar [K66, S66]:[1]

(0.1) *Prequantization*: Construct over $X$ a hermitian line bundle with connection, $L$, whose curvature is the given symplectic form. When an equivariant momentum map exists, $G$ lifts to act on $L$ and one obtains a unitary $G$-module of sections, $H(L)$, which however is 'too large'.

(0.2) *Polarization*: Somehow extract a smaller module $H^{\bullet}(L)$ from $H(L)$. This usually involves looking at sections constant along a lagrangian subbundle of $TX^{\mathbf{C}}$, or at the cohomology of such sections; $H^{\bullet}(L)$ should be largely independant of the subbundle, unitarizable, and irreducible when $G$ is transitive on $X$—i.e. when $X$ is, up to covering, a coadjoint orbit of $G$.

Step (0.1) is quite well understood. When $X$ is $\mathbf{R}^{2n}$ it goes back to Sophus Lie [L90, p. 270]; in general the existence of $L$ requires certain integrality conditions, and when nonunique $L$ seems to be just the extra structure needed to specify a representation.

---

[1]The outline (0.1–2) is perhaps abrupt; definitions and examples are forthcoming.

In contrast, (0.2) only gives *principles* which all suffer exceptions and require debated exegesis. Suitably implemented, they do indeed attach almost all known irreducible unitary representations to coadjoint orbits [B54, K62, A71, D84, V87], sometimes in many ways. But a few representations especially prominent in physics, like the metaplectic, still escape the picture, and altogether the established correspondence remains somewhat 'magic' [V87]; we have an eclectic book of answers but ... to what question exactly? Is there, in general, something like the Stone-von Neumann characterization of the Schrödinger representation?

With an eye on these and other problems, J. M. Souriau [S88, S92] started a quite different attack on the whole matter. Rather than pursue the spell of (0.1–2), he proposed to write down minimal intrinsic properties that a representation should have to be 'attached' to a given coadjoint orbit, $X$. He arrived at the following definition: a *quantum representation* (for $X$) is a unitary $G$-module $\mathcal{H}$ such that, for each unit vector $\varphi$ in $\mathcal{H}$, the matrix element $m(g) = (\varphi, g\varphi)$ satisfies

$$\left|\sum_{j=1}^n c_j m(\exp(Z_j))\right| \leq \sup_{x \in X} \left|\sum_{j=1}^n c_j e^{i\langle x, Z_j\rangle}\right| \qquad (0.3)$$

for any given complex numbers $c_j$ and commuting $Z_j$ in the Lie algebra of $G$. Such functions $m$ are called *quantum states*, and are the basic object of study. As we will show in (2.9), the inequality can be interpreted as requiring that

> *the quantum spectrum of 'commuting observables'*
> *is concentrated on their classical range, suitably compactified.* (0.4)

Our purpose here is to investigate the consequences of this principle, and determine its relation with known orbit constructions in a reasonably wide range of cases. We emphasize from the outset that one should *not* expect to directly recover the traditional theories from (0.3), be it only because this definition, originally made with infinite-dimensional groups in mind, does not impose any continuity property. Therefore we shall concurrently also study two successively stronger variants of it, to be called *0-quantum representations* (wherein one simply forgets about the compactification in (0.4), which among other things ensures continuity) and *∞-quantum representations* (which builds in a certain recursivity with respect to 'reduction by commuting integrals').

$$*\quad *\quad *\quad *$$

**Notation.** All manifolds and Lie groups in this paper are understood to be Hausdorff and countable at $\infty$; they may be disconnected. If $G$ is a Lie group, we reserve the corresponding german letter for its Lie algebra, $\mathfrak{g}$. If $G$ acts on a manifold, $X$, so that we have a morphism $g \mapsto g_X$ of $G$ into the diffeomorphisms of $X$ (with $g_X(x)$ a $C^\infty$ function of the pair $(g, x)$), we define the corresponding infinitesimal action

$$Z \mapsto Z_X, \qquad \mathfrak{g} \to \text{vector fields on } X \qquad (0.5)$$

by $Z_X(x) = \frac{d}{dt} \exp(tZ)_X(x)\big|_{t=0}$. This is a morphism of Lie algebras, provided we define the bracket of vector fields with minus its usual sign. Whenever possible, we drop the subscripts to write $g(x)$ and $Z(x)$, instead of $g_X(x)$ and $Z_X(x)$. We use

$$G(x), \qquad \mathfrak{g}(x), \qquad G_x, \qquad \mathfrak{g}_x, \qquad (0.6)$$

to denote respectively the $G$-orbit of $x$, its tangent space at $x$, the stabilizer of $x$ in $G$, and the stabilizer of $x$ in $\mathfrak{g}$.

It will be convenient to have a concise notation concerning the translation of tangent and cotangent vectors to the group; thus for fixed $g, q \in G$ we will let

$$\begin{array}{ccc} T_q G \to T_{gq} G & & T_q^* G \to T_{gq}^* G \\ v \mapsto gv, & \text{resp.} & p \mapsto gp \end{array} \qquad (0.7)$$

denote the derivative of $q \mapsto gq$, respectively the transposed map so that $\langle gp, v\rangle = \langle p, g^{-1}v\rangle$. Likewise we define $vg$ and $pg$ with $\langle pg, v\rangle = \langle p, vg^{-1}\rangle$. Finally we recall the coadjoint action, $g(x) = gxg^{-1}$, on $\mathfrak{g}^* = T_e^* G$; infinitesimally it gives $Z(x) = \langle x, [\,\cdot\,, Z]\rangle$.

*It became imperative for me to brush up my knowledge of Poisson brackets.*
—P. A. M. Dirac

# 1. Preliminaries

## 1.1 Preliminaries on symplectic manifolds

**A. Hamiltonian vector fields.** Let $X$ be a symplectic manifold—a manifold with a nondegenerate, closed 2-form $\sigma$. A vector field $\eta$ on $X$ is called *symplectic* if its flow preserves the 2-form: $\pounds(\eta)\sigma = 0$. By Cartan's formula for the Lie derivative, $\pounds(\eta)\sigma = i(\eta)d\sigma + di(\eta)\sigma$, this holds iff the 1-form $i(\eta)\sigma = \sigma(\eta, \cdot)$ is closed. When this 1-form is actually exact, so that there is a function $H$ on $X$ with

$$i(\eta)\sigma = -dH, \tag{1.1}$$

we say that $\eta$ is *hamiltonian*, $\eta \in \mathrm{Ham}(X)$, and we write $\eta = \mathrm{drag}\, H$ ('symplectic gradient'). The field $\eta$ determines $H$ up to a local additive constant, so if $X$ has $c$ connected components we have an exact sequence

$$0 \longrightarrow \mathbf{R}^c \longrightarrow \mathrm{C}^\infty(X) \xrightarrow{\mathrm{drag}} \mathrm{Ham}(X) \longrightarrow 0. \tag{1.2}$$

This is a central extension of Lie algebras when $\mathrm{C}^\infty(X)$ is endowed with *Poisson bracket*: $\{H, H'\} = \sigma(\mathrm{drag}\, H', \mathrm{drag}\, H)$.

**B. Hamiltonian $G$-spaces.** The $\sigma$-preserving action of a Lie group $G$ on $X$ is called hamiltonian if (0.5) takes its values in $\mathrm{Ham}(X)$; or in other words, if there is a *momentum map*, $\Phi : X \to \mathfrak{g}^*$, such that

$$i(Z_X)\sigma = -dH_Z \qquad \text{where} \qquad H_Z = \langle \Phi(\cdot), Z \rangle. \tag{1.3}$$

The momentum map is called *equivariant*, and the pair $(X, \Phi)$ a *hamiltonian G-space*, if $\Phi$ intertwines the action of $G$ on $X$ with the coadjoint action on $\mathfrak{g}^*$. Infinitesimally this means that $Z \mapsto H_Z$ is a morphism:

$$\{H_Z, H_{Z'}\} = H_{[Z,Z']}. \tag{1.4}$$

(As one knows, this requirement is no true restriction; it can always be met by passing to a central extension of $\mathfrak{g}$ by $\mathbf{R}^c$ (1.2).) This morphism pulls the first two terms of (1.2) back to ideals of $\mathfrak{g}$ which we shall often encounter and denote

$$\mathfrak{o} = \mathrm{Ker}(Z \mapsto H_Z), \qquad \mathfrak{c} = \mathrm{Ker}(Z \mapsto \mathrm{drag}\, H_Z). \tag{1.5}$$

**C. Basic examples.** (a) If $G$ acts on a manifold $Q$, we get an action on $X = T^*Q$ which preserves the canonical 1-form $\varpi = \text{``}\langle p, dq\rangle\text{''}$ and $\sigma = d\varpi$. In this case, $i(Z_X)d\varpi + di(Z_X)\varpi = 0$ shows that a momentum map, which one can see is equivariant, is given by E. Nœther's formula: $H_Z = i(Z_X)\varpi$, i.e.,

$$\langle \Phi(p), Z\rangle = \langle p, Z(q)\rangle, \qquad p \in T_q^*Q. \tag{1.6}$$

(b) If $X$ is an orbit in $\mathfrak{g}^*$ for the coadjoint action of $G$, then the 2-form defined on it by $\sigma(Z(x), Z'(x)) = \langle Z(x), Z'\rangle$ makes $(X, X \hookrightarrow \mathfrak{g}^*)$ into a homogeneous hamiltonian $G$-space. Conversely (Kirillov-Kostant-Souriau) every such space covers a coadjoint orbit:

**1.7 Theorem.** *Let $(X, \Phi)$ be a hamiltonian $G$-space, and suppose that $G$ acts transitively on $X$. Then $\Phi$ is a symplectic covering of its image, which is a coadjoint orbit of $G$.*

By *symplectic* covering we mean of course that $\Phi$ pulls the symplectic structure of the orbit back to the given one on $X$; this is just a restatement of (1.4). That $\Phi$ is a covering (has discrete fibers) follows from the first of two informative consequences of (1.3), valid for any momentum map:

$$\mathrm{Ker}(D\Phi(x)) = \mathfrak{g}(x)^\sigma, \qquad \mathrm{Im}(D\Phi(x)) = \mathrm{orth}(\mathfrak{g}_x). \tag{1.8}$$

Here the superscript means orthogonal subspace relative to $\sigma$, and if $(\cdot)$ is a subset of either $\mathfrak{g}$ or $\mathfrak{g}^*$, $\mathrm{orth}(\cdot)$ denotes its annihilator in the other.

**D. Symplectic induction.** A symplectic version of Mackey's unitary induction process was introduced by Kazhdan, Kostant and Sternberg [K78]. Given a closed subgroup $H$ of $G$ and a hamiltonian $H$-space $(Y, \Psi)$, it produces a hamiltonian $G$-space $(\operatorname{Ind}_H^G Y, \Phi_{\text{ind}})$ as follows. (We use the notation (0.7).)

First endow $N = T^*G \times Y$ with the symplectic form $\omega = d\varpi + \tau$, where $\varpi$ is the canonical 1-form on $T^*G$ and $\tau$ the given 2-form on $Y$; and let $H$ act 'diagonally' on $N$ by $h(p, y) = (ph^{-1}, h(y))$. This is hamiltonian, with momentum map $\psi$:

$$\psi(p, y) = \Psi(y) - q^{-1}p_{|\mathfrak{h}} \qquad (1.9)$$

for $p \in T_q^*G$. Here the second term denotes the restriction to $\mathfrak{h}$ of $q^{-1}p \in \mathfrak{g}^*$; it comes from (1.6) with $Z(q) = -qZ$. The induced manifold is now defined as the Marsden-Weinstein reduction of $N$ at zero, i.e.,

$$\operatorname{Ind}_H^G Y := \psi^{-1}(0)/H. \qquad (1.10)$$

In more detail: the action of $H$ is free and proper (because it is free and proper on the factor $T^*G$, where it is the right action of $H$ regarded as a subgroup of $T^*G$ [B60, ch. III, p. 106]); so $\psi$ is a submersion (1.8b), $\psi^{-1}(0)$ is a submanifold, and (1.10) is a manifold; moreover $\omega_{|\psi^{-1}(0)}$ degenerates exactly along the $H$-orbits (1.8a), so it comes from a uniquely defined symplectic form, $\sigma_{\text{ind}}$, on the quotient.

To make (1.10) into a $G$-space, we let $G$ act on $N$ by $g(p, y) = (gp, y)$. This action commutes with the diagonal $H$-action and preserves $\psi^{-1}(0)$. Moreover it is hamiltonian and its momentum map $\phi : N \to \mathfrak{g}^*$, given by (1.6) with now $Z(q) = Zq$:

$$\phi(p, y) = pq^{-1}, \qquad p \in T_q^*G, \qquad (1.11)$$

is constant on the $H$-orbits. Passing to the quotient, we obtain the required $G$-action on $\operatorname{Ind}_H^G Y$ and momentum map $\Phi_{\text{ind}} : \operatorname{Ind}_H^G Y \to \mathfrak{g}^*$.

We record a few elementary properties of the construction. First of all, (1.10) has dimension $\dim(N) - \dim(\mathfrak{h}^*) - \dim(H)$, i.e.,

$$\dim(\operatorname{Ind}_H^G Y) = 2\dim(G/H) + \dim(Y). \qquad (1.12)$$

Secondly we have $\mathrm{Im}(\Phi_{\mathrm{ind}}) = \phi(\psi^{-1}(0))$, which can be read as follows: if $W$ is a coadjoint orbit of $G$, then

$$W \text{ intersects } \mathrm{Im}(\Phi_{\mathrm{ind}}) \quad \Leftrightarrow \quad W_{|\mathfrak{h}} \text{ intersects } \mathrm{Im}(\Psi) \qquad (1.13)$$

('Frobenius reciprocity'). Third we have the 'stages theorem': if $K$ is an intermediate closed subgroup, then

$$\mathrm{Ind}_K^G \mathrm{Ind}_H^K Y = \mathrm{Ind}_H^G Y. \qquad (1.14)$$

Indeed, the left-hand side is by construction a space of $K \times H$-orbits within $T^*G \times T^*K \times Y$, and it is not hard to verify that an isomorphism from left to right is obtained by sending the $K \times H$-orbit of $(p, p', y)$, $p' \in T^*_{q'}K$, to the $H$-orbit of $(pq', y)$.

## 1.2 Preliminaries on positive-definite functions

**A. Positive-definite functions; states.** Let $G$ be a group, with identity element $e$. A complex-valued function $m$ on $G$ is called *positive-definite* if the sesquilinear form

$$(c, d)_m := \sum_{g,h \in G} \bar{c}_g d_h m(g^{-1}h), \qquad (1.15)$$

defined on $\mathbf{C}[G] = \{$complex-valued functions with finite support on $G\}$, is positive: $(c, c)_m \geq 0$. If further $m(e) = 1$, then $m$ is called a *state* of $G$.

We can identify each function $m$ on $G$ with the linear functional on $\mathbf{C}[G]$ defined by $m(\delta^g) = m(g)$, where $\delta^g$ denotes the basis function which is one at $g$ and zero elsewhere; then (1.15) writes

$$(c, d)_m = m(c^* \cdot d), \qquad (1.16)$$

where we are using the $^*$-algebra structure of $\mathbf{C}[G]$: $\delta^g \cdot \delta^h = \delta^{gh}$, $\delta^{g*} = \delta^{g^{-1}}$; so states are the same as normalized positive linear functionals on $\mathbf{C}[G]$.

They are covariant objects: the pull-back, $m \circ \alpha$, of a state $m$ by a group morphism $\alpha$, is a state. Conversely, if $\alpha$ is onto and $m \circ \alpha$ is a state, then $m$ is a state.

**B. Remarkable inequalities.** Any state gives rise to a Cauchy-Schwarz inequality, $|(c,d)_m|^2 \leq (c,c)_m (d,d)_m$. Three inequalities of Herglotz, Kreĭn, and Weil [H63]:

$$|m(g)| \leq 1, \qquad (1.17)$$

$$|m(g) - m(h)| \leq \sqrt{2\mathrm{Re}(1 - m(g^{-1}h))}, \qquad (1.18)$$

$$|m(gh) - m(g)m(h)| \leq \sqrt{1 - |m(g)|^2}\sqrt{1 - |m(h)|^2}, \qquad (1.19)$$

follow for all $g, h$ in $G$ by taking the pair $c^*, d$ to be $\delta^e, \delta^g$; resp. $\delta^e, \delta^g - \delta^h$; resp. $\delta^g - m(g)\delta^e$, $\delta^h - m(h)\delta^e$. Because of (1.19), the equation $m(h) = 1$ defines a subgroup $H$ of $G$, on whose cosets $m$ is constant; so $m$ comes from a function $\dot{m}$ on $G/H$ (which is a state when $H$ is normal).

**C. States and representations.** States arise from unitary representations and conversely (Gel'fand-Naĭmark-Segal):

**1.20 Theorem.** *A function $m$ on $G$ is a state iff there are a unitary $G$-module $\mathcal{H}$, and a unit vector $\varphi$ in $\mathcal{H}$, such that $m(g) = (\varphi, g\varphi)$.*

Sufficiency is immediate. Conversely if $m$ is a state, one observes that the form (1.16) on $\mathbf{C}[G]$ is invariant under the regular action, $gc = \delta^g \cdot c$; dividing out null vectors and completing, one gets a unitary $G$-module $\mathcal{H}_m$ in which the stated identity holds with $\varphi$ the class of $\delta^e$.

A practical way to complete here is to take the *antidual* [S62]: we realize $\mathcal{H}_m$ as the space of all antilinear functionals $f$ on $\mathbf{C}[G]$, such that the quantity

$$\|f\|^2 := \sup_{c \in \mathbf{C}[G]} \frac{|f(c)|^2}{(c,c)_m} \qquad \text{is finite.} \qquad (1.21)$$

(We understand that the numerator must vanish when the denominator does, so that $f$ factors through the null vectors.) Each $c \in \mathbf{C}[G]$ defines an element $e_c = (\,\cdot\,, c)_m \in \mathcal{H}_m$; these are dense in $\mathcal{H}_m$, and are also characterized by

$$f(c) = (e_c, f) \qquad \forall f \in \mathcal{H}_m. \qquad (1.22)$$

We will abbreviate $f(\delta^g)$ to $f(g)$ and $e_{\delta^g}$ to $e_g$. In this way $\mathcal{H}_m$ is a Hilbert space of functions on $G$, with cyclic vector $e_e = \overline{m}$, in which the norm is given by (1.21) and $G$ acts by $(gf)(g') = f(g^{-1}g')$.

This realization is especially well suited to discuss intertwining operators: if $\tilde{m}$ is another state, each bounded intertwining operator $J : \mathcal{H}_m \to \mathcal{H}_{\tilde{m}}$ will be characterized by the single function $Je_e$. More precisely, writing $\mathcal{H}_m^\vee$ for the image of $\mathcal{H}_m$ under $f \mapsto f^\vee := \overline{f(\,\cdot\,^{-1})}$, we have

**1.23 Proposition.** $J \mapsto Je_e$ is an injection $\mathrm{Hom}_G(\mathcal{H}_m, \mathcal{H}_{\tilde{m}}) \to \mathcal{H}_m^\vee \cap \mathcal{H}_{\tilde{m}}$.

**Proof.** By hypothesis the function $j = Je_e$ is in $\mathcal{H}_{\tilde{m}}$ and satisfies $gj = Je_g$. Thus, by (1.22), the adjoint of $J$ is given by $(J^*f)(g) = (e_g, J^*f) = (gj, f)$. In particular it maps the cyclic vector $\tilde{e}_e$ of $\mathcal{H}_{\tilde{m}}$ onto $j^\vee$. Therefore $j^\vee$ is in $\mathcal{H}_m$, and determines $J$ just as above: $(Jf)(g) = (gj^\vee, f)$. **q.e.d.**

**D. Application: discrete induction.** As an example to be used later, we compute the space $\mathcal{H}_m$ when $\chi$ is a character of a subgroup $H \subset G$ and

$$m(g) = \chi^\bullet(g) := \begin{cases} \chi(g) & \text{if } g \in H, \\ 0 & \text{otherwise.} \end{cases} \qquad (1.24)$$

This is a state, as one sees by splitting the sum (1.15) over the cosets of $H$. We claim that $\mathcal{H}_m$ is the induced representation $\mathrm{ind}_H^G \chi$, where $G$ and $H$ are regarded as discrete groups.[2] In other words:

**1.25 Proposition.** $\mathcal{H}_{\chi^\bullet}$ consists exactly of all $f : G \to \mathbf{C}$ such that
 (a) $f(gh) = \overline{\chi}(h)f(g)$ for all $h \in H$;
 (b) $\|f\|_\star^2 := \sum_{gH \in G/H} |f(gH)|^2$ is finite.

---

[2] 'ind' will always denote the induction functor for discrete groups, as opposed to the usual 'Ind' when $G$ already has another topology.

**Proof.** Suppose that $f$ satisfies (1.21). If we take $c = \delta^{gh} - \chi(h)\delta^g$ there, we obtain $(c,c)_{\chi\bullet} = 0$; therefore $f(c) = 0$, whence (a). On the other hand if we take $c = \sum_{g \in \Gamma} f(g)\delta^g$ where $\Gamma$ is a finite subset of $G$ with at most one point in each $H$-coset, then the quotient in (1.21) equals $\sum_{g \in \Gamma} |f(g)|^2$. Therefore we have $\|f\|_\star^2 \leq \|f\|^2$, which proves (b).

Conversely assume that $f$ satisfies (a, b). Then a short computation gives $f(c) = (e_c, f)_\star$, whence $\|f\|^2 \leq \|f\|_\star^2$ by Cauchy-Schwarz.    **q.e.d.**

**1.26 Remarks.**  (a) This Proposition, implicit in [G49], is much generalized in [B63]. For instance we can replace $\chi$ in (1.24) by any state $n$ of $H$, and again $\mathcal{H}_{n\bullet}$ is isomorphic to $\mathrm{ind}_H^G \mathcal{H}_n$ (under $f \mapsto F$, $F(g)(h) = f(gh)$).

(b) (1.25a), which means that we are really dealing with sections of an associated bundle, remains true with the same proof even when $\chi$ is extended by something nonzero outside $H$. For an example see (6.13).

(c) Given a character $\eta$ of another subgroup $K \subset G$, (1.23) and (1.25) combine to give the Mackey-Shoda bounds on intertwining numbers [M51]:

**1.27 Proposition.** *The dimension of* $\mathrm{Hom}_G(\mathcal{H}_{\chi\bullet}, \mathcal{H}_{\eta\bullet})$ *is bounded by the number of double cosets* $D = HaK$ *such that*

   (a) $\chi(h) = \eta(a^{-1}ha)$ *for all* $h \in H \cap aKa^{-1}$;
   (b) $D$ *projects onto finite sets in both* $G/K$ *and* $H\backslash G$.

**Proof.** By (1.23), this dimension does not exceed that of $\mathcal{H}_{\chi\bullet}^\vee \cap \mathcal{H}_{\eta\bullet}$, whose members $j$ satisfy $j(h^{-1}ak) = \chi(h)\overline{\eta}(k)j(a)$ by virtue of (1.25a).

Such a function is determined by one value per double coset $D = HaK$. This value must vanish when (a) fails: try $k = a^{-1}ha$; also when (b) fails: indeed $|j|^2$ is constant in $D$, and this constant occurs $\sharp(D/K)$ times in the series for $\|j\|^2$, resp. $\sharp(H\backslash D)$ times in the series for $\|j^\vee\|^2$ (1.25b).    **q.e.d.**

**E. The abelian case.** Suppose that $G$ is a locally compact abelian group, and let $\hat{G}$ denote its dual, i.e. the group of all continuous characters $\chi$ of $G$ with the topology of uniform convergence on compact sets. *Bochner's Theorem* [H63] asserts that the Fourier transform $\mu \mapsto m$,

$$m(g) = \int_{\hat{G}} \chi(g) \, d\mu(\chi), \tag{1.28}$$

gives a bijection between all continuous positive-definite functions $m$ on $G$, and all positive bounded Radon measures $\mu$ on $\hat{G}$. Likewise, *Stone's Theorem* [F88] asserts that the continuous unitary representations $U$ of $G$ correspond to all projection-valued measures (p.v.m.) $E$ on $\hat{G}$ under

$$U(g) = \int_{\hat{G}} \chi(g) \, dE(\chi). \tag{1.29}$$

Here we recall that a *p.v.m.* on a locally compact space $X$ is a *-representation $E$ of the *-algebra $C_0 = \{$continuous functions on $X$ vanishing at $\infty\}$ (pointwise operations) in a Hilbert space $\mathcal{H}$, such that $E(C_0)\mathcal{H}$ is dense in $\mathcal{H}$. Then $(\varphi, E(\cdot)\varphi)$ is a positive bounded Radon measure for each $\varphi \in \mathcal{H}$, which allows us to extend $E$ to all Borel functions. Finally $\int_X f(x) \, dE(x)$ is another notation for $E(f)$. (See [F88], pp. 125, 387, 116.)

We shall refer to $\mu$ and $E$ as the *spectral measures* of $m$ and $U$.

**1.30 Remark.** This correspondence is covariant: if $\alpha$ is a continuous morphism $H \to G$ with dual $\beta : \hat{G} \to \hat{H}$, the spectral measure of $m \circ \alpha$ is $\beta(\mu)$. Likewise $U \circ \alpha$ has spectral measure $\beta(E)$, where $\beta(E)(f) = E(f \circ \beta)$.

**F. Concentration of spectral measures.** We keep the above notation. Following [B65], we say that $\mu$ is *concentrated* on a set $M$ (written $\mu \Subset M$) if its complement $M^c$ is locally $\mu$-null. For a positive bounded $\mu$ this means that $\inf\{\mu(\Omega) : \Omega \text{ open}, M^c \subset \Omega\} = 0$. (See [B65], ch. IV, cor. 1 p. 172.)

Likewise we say that $E$ is concentrated on $M$ if $(\varphi, E(\cdot)\varphi)$ is for all $\varphi$, and we call $m$ or $U$ concentrated on $M$ *(in $\hat{G}$)* when its spectral measure is. We collect here what we shall need concerning this notion. First we have some facts from [B65], ch. V, pp. 70, 125, 109, 59:[†]

**1.31 Proposition.** *Let $\beta : X \to Y$ be a continuous map of locally compact spaces, $\nu$ a positive measure on $Y$, and $N$ a subset of $Y$.*
*(a) If $\nu = \beta(\mu)$ for some positive $\mu$ on $X$, then $\nu \Subset N \Leftrightarrow \mu \Subset \beta^{-1}(N)$.*
*(b) If $X$ is countable at $\infty$, then $\nu \Subset \beta(X) \Rightarrow \nu = \beta(\mu)$ for some positive $\mu$.*

---

[†]Voir aussi l'Annexe A.

(c) If $\nu$ is the vague limit of measures $\nu_n \in N$, then $\nu \in \text{closure}(N)$.
(d) If $N$ is closed, then $\nu \in N \Leftrightarrow \text{supp}(\nu) \subset N$.

Of course we have corresponding statements for projection-valued measures. Secondly we will need:

**1.32 Proposition.** *Let $E$ be a projection-valued measure on $X$ and $M$ a subset of $X$. Then $\mathcal{H}_M := \{\,\varphi \in \mathcal{H} : (\varphi, E(\cdot)\varphi) \text{ is concentrated on } M\,\}$ is a closed linear subspace of $\mathcal{H}$.*

**Proof.** By definition $\varphi$ is in $\mathcal{H}_M$ iff given $\varepsilon > 0$ we can find an open set $\Omega$, containing $M^c$, such that $(\varphi, E(\Omega)\varphi) = \|E(\Omega)\varphi\|^2 < \varepsilon$. Clearly $\mathcal{H}_M$ is stable under scalar multiplication. It is stable under addition, because if $\varphi = \varphi_1 + \varphi_2$ and $\Omega = \Omega_1 \cap \Omega_2$ then $\|E(\Omega)\varphi\| \leq \|E(\Omega_1)\varphi_1\| + \|E(\Omega_2)\varphi_2\|$. Finally it is closed, for if $\varphi$ lies in its closure then we can find $\psi$ in $\mathcal{H}_M$ with $\|\varphi - \psi\| < \varepsilon$ and $\Omega$ with $\|E(\Omega)\psi\| < \varepsilon$, whence $\|E(\Omega)\varphi\| < 2\varepsilon$.     **q.e.d.**

# 2. Quantum states

From now on, $X$ denotes a coadjoint orbit of a Lie group $G$ with Lie algebra $\mathfrak{g}$. Thus $X$ is a symplectic manifold, and we have a morphism $Z \mapsto H_Z$ of $\mathfrak{g}$ into the smooth functions on $X$ under Poisson bracket:

$$H_Z(x) = \langle x, Z \rangle, \qquad \{H_Z, H_{Z'}\} = H_{[Z,Z']}. \tag{2.1}$$

## 2.1 What is a quantum state?

We start with the following heuristic consideration. If we think of $X$ as describing some classical dynamical system, then a *statistical state* is simply a probability measure, $\mu$, on $X$. We can regard $\mu$ as a measure on $\mathfrak{g}^*$, and characterize it by a state (as in §1.2) of the additive group $\mathfrak{g}$, namely its Fourier transform

$$\hat{\mu}(Z) = \int_{\mathfrak{g}^*} e^{i\langle x, Z\rangle} d\mu(x). \tag{2.2}$$

So statistical mechanics on $X$ can be formulated in terms of states of $\mathfrak{g}$. Together with (1.20), this may suggest that 'quantum mechanics on $X$' should be the study of appropriate states of $G$ itself. Here is the class of states proposed in [S88, S92]:

**2.3 Definition.** A *quantum state for $X$* is a state $m$ of $G$ such that, for all choices of $n \in \mathbf{N}$, $c_1, \ldots, c_n \in \mathbf{C}$ and $Z_1, \ldots, Z_n \in \mathfrak{g}$ with $\{H_{Z_j}, H_{Z_k}\} = 0$, we have

$$\left| \sum_{j=1}^n c_j m(\exp(Z_j)) \right| \leq \sup_{x \in X} \left| \sum_{j=1}^n c_j e^{i\langle x, Z_j \rangle} \right|. \tag{2.4}$$

**2.5 Remark.** This definition is actually slightly stronger than Souriau's (0.3): we require the inequality whenever the *hamiltonians* $H_{Z_j}$ commute, rather than the $Z_j$ themselves; this is more natural, as (2.25) will show. To help us think modulo the "extraneous" ideal $\mathfrak{o} = \text{orth}(X) = \text{Ker}(Z \mapsto H_Z)$ (1.5, 2.1), which is what this change amounts to, we will call a subalgebra $\mathfrak{a}$ of $\mathfrak{g}$

$$X\text{-abelian if } [\mathfrak{a}, \mathfrak{a}] \subset \mathfrak{o}, \qquad X\text{-central if } [\mathfrak{g}, \mathfrak{a}] \subset \mathfrak{o}. \qquad (2.6)$$

We note also that a choice of units and of $\sqrt{-1}$ is implicit in the definition; many authors would replace $i$ by $-i$ (or $-i/\hbar$) in all our formulas.

## 2.2 A geometrical criterion

The statistical states (2.2) were *concentrated on* $X$, in the sense of §1.2F. In this section, we want to substantiate the idea that Definition 2.3 means something similar. More precisely, if $\mathfrak{a}$ is $X$-abelian and $m$ is a quantum state for $X$, we are going to see that the function

$$m \circ \exp_{|\mathfrak{a}} \qquad (2.7)$$

is a state of the additive group $\mathfrak{a}$, and is concentrated somewhere.

To make sense of this, we need a topology on $\mathfrak{a}$. Implicitly we already have the usual one, for which the dual group may—and will—be identified with $\mathfrak{a}^*$ by letting $p \in \mathfrak{a}^*$ stand for the character $e^{i\langle p, \cdot \rangle}$ of $\mathfrak{a}$. But these continuous characters are not enough to deal with our possibly discontinuous states; rather than $\mathfrak{a}^*$ we must use its Bohr compactification $\hat{\mathfrak{a}}$:

$$\hat{\mathfrak{a}} = \{ \textit{all } \text{characters of } \mathfrak{a} \} = \text{dual of the } \textit{discrete } \text{group } \mathfrak{a}. \qquad (2.8)$$

This is a compact group, in which $\mathfrak{a}^*$ is dense [H63]. Likewise we define $\hat{\mathfrak{g}}$ and regard $\mathfrak{g}^*$ as a dense subgroup; in doing so we must be careful to distinguish between usual closure in $\mathfrak{g}^*$ and closure in $\hat{\mathfrak{g}}$, which we shall denote by $X \to bX$, for Bohr closure or compactification.

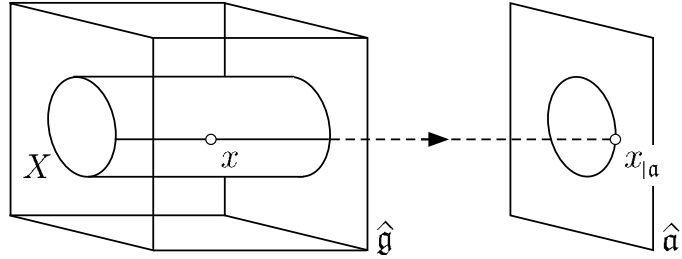

**Fig. 1** Projection of $X$ in the dual of an $X$-abelian subalgebra.

Now restricting characters defines a continuous projection $\hat{\mathfrak{g}} \to \hat{\mathfrak{a}}$ (Fig. 1), and we can assert:

**2.9 Proposition.** *A state $m$ of $G$ is quantum for $X$ iff for each $X$-abelian $\mathfrak{a} \subset \mathfrak{g}$, $m \circ \exp_{|\mathfrak{a}}$ is a state concentrated on $bX_{|\mathfrak{a}}$ (in $\hat{\mathfrak{a}}$).*

**Proof.** First we reassure ourselves that the notation $bX_{|\mathfrak{a}}$ is unambiguous, i.e. $(bX)_{|\mathfrak{a}} = b(X_{|\mathfrak{a}})$: the projection of $bX$ lies in the closure of $X_{|\mathfrak{a}}$ by continuity; moreover it is compact and so contains this closure.

Assume that $m$ is quantum for $X$, and let $\mathfrak{a}$ be $X$-abelian. We must first verify that $m \circ \exp_{|\mathfrak{a}}$ is a state. This is clear for abelian $\mathfrak{a}$, but in general we must argue as follows: Let $\tilde{m}$ denote the pull-back of $m$ to the universal covering $\tilde{G}$ of the identity component of $G$, and $\tilde{O}$ the closed [B60] subgroup of $\tilde{G}$ generated by $\mathfrak{o}$ (2.6). Taking $(c_1, Z_1, c_2, Z_2) = (1, Z, -1, 0)$ in (2.4) gives

$$\left|\tilde{m}(\exp(Z)) - 1\right| \leq \sup_{x \in X} \left|e^{i\langle x, Z\rangle} - 1\right| = 0 \qquad \forall Z \in \mathfrak{o}. \qquad (2.10)$$

So $\tilde{m}$ is trivial on $\tilde{O}$, and consequently comes from a state $\dot{m}$ of $\dot{G} = \tilde{G}/\tilde{O}$ (§ 1.2B). Writing $\dot{\mathfrak{a}}$ for the (abelian) projection of $\mathfrak{a}$ in $\mathfrak{g}/\mathfrak{o}$, we obtain a diagram

$$\begin{array}{ccccc} \mathfrak{a} & \xrightarrow{\exp} & \tilde{G} & \longrightarrow & G \\ \downarrow & & \downarrow & & \downarrow m \\ \dot{\mathfrak{a}} & \xrightarrow{\exp} & \dot{G} & \xrightarrow{\dot{m}} & \mathbf{C}. \end{array} \qquad (2.11)$$

This exhibits $m \circ \exp_{|\mathfrak{a}}$ as the pull-back of $\dot{m}$ by a morphism $(\mathfrak{a} \to \dot{\mathfrak{a}} \to \dot{G})$ and hence as a state, as claimed. Now by Bochner's theorem (1.28) applied with the discrete topology on $\mathfrak{a}$, this state has a spectral measure $\nu$:

$$(m \circ \exp_{|\mathfrak{a}})(Z) = \int_{\hat{\mathfrak{a}}} \chi(Z) \, d\nu(\chi). \tag{2.12}$$

So (2.4) says that

$$\left| \int_{\hat{\mathfrak{a}}} f(\chi) \, d\nu(\chi) \right| \leq \sup_{\chi \in X_{|\mathfrak{a}}} |f(\chi)| \tag{2.13}$$

for all trigonometric polynomials $f(\chi) = \sum_j c_j \chi(Z_j)$ with $c_j \in \mathbf{C}$, $Z_j \in \mathfrak{a}$. By Stone-Weierstrass, these are uniformly dense in the continuous functions on $\hat{\mathfrak{a}}$, so therefore (2.13) still holds for all continuous $f$. In particular if $f$ vanishes on $bX_{|\mathfrak{a}}$ then $\nu(f) = 0$, which is to say that

$$\mathrm{supp}(\nu) \subset bX_{|\mathfrak{a}}, \tag{2.14}$$

or in other words, that $\nu$ is concentrated on $bX_{|\mathfrak{a}}$ (1.31d).

Conversely let $c_j$ and $Z_j$ be given as in Definition 2.3. Then $\mathfrak{o}$ and the $Z_j$ span an $X$-abelian subalgebra $\mathfrak{a}$ of $\mathfrak{g}$, and $f(\chi) = \sum_j c_j \chi(Z_j)$ defines a continuous function on $\hat{\mathfrak{a}}$. Assuming (2.14) for $\mathfrak{a}$, the mean value inequality gives us (2.13) and hence (2.4).     **q.e.d.**

**2.15 Remark.** As the proof shows, it would be equivalent to consider only *maximal* $X$-abelian subalgebras in (2.9).

**2.16 Interpretation.** (2.14) means that the spectral measure $\nu$ looks as if it came from a measure on $X$, or perhaps $bX$ (cf. 1.31b); thus $X$ allows us to 'statistically interpret' the state as saying, for any given $\mathfrak{a}$, in which fiber of the projection $\pi : X \to X_{|\mathfrak{a}}$ the point $x$ lies (Fig. 1). Aptly, the need for such an interpretation is the main reason why each quantum system must have an underlying classical system. Note that:

  (a) Our use of the Bohr closure $bX_{|\mathfrak{a}}$, embedded in $\hat{\mathfrak{a}}$, suppresses here the need for the generalized probabilities and statistical spaces of [S88, S92]; it essentially adds 'fibers at $\infty$'.

(b) The most that $\nu$ might tell us, according to this interpretation, is that $x$ lies in a given fiber of $\pi$. The 'uncertainty principle' finds its expression then in the fact that the fibers above regular values are coisotropic: $\mathrm{Ker}(D\pi(x))^\sigma \subset \mathrm{Ker}(D\pi(x))$ (cf. 1.8); in particular they have at least half the dimension of $X$.

(c) We can tentatively rephrase the interpretation of $\nu$ as follows: If a function $H$ on $X$ factors through a continuous function $h$ on $\hat{\mathfrak{a}}$ (Fig. 1), we regard the image $h(\nu)$ as the probability distribution of the variable $H(x)$ in the state $m$. By (2.14) this measure only depends on $h$ via $H$, as it should; it *might* depend on the particular $\mathfrak{a}$ used in factoring $H$, although we hope not.

## 2.3 0-quantum states

The main lesson of (2.9) is that our states are the 'quantum analogues' of probability measures on $bX$, rather than $X$. This will result in a much looser correspondence than usual (see e.g. (4.10c)), and suggests that we also study a more restrictive definition, omitting the '$b$' in (2.9):

**2.17 Definition.** A *0-quantum state for $X$* is a state $m$ of $G$ such that, for each $X$-abelian $\mathfrak{a} \subset \mathfrak{g}$, $m \circ \exp_{|\mathfrak{a}}$ is a state concentrated on $X_{|\mathfrak{a}}$ (in $\hat{\mathfrak{a}}$).

Clearly 0-quantum $\Rightarrow$ quantum; but in general the converse will fail unless $X$ is compact (so that $bX = X$). Indeed we shall meet discontinuous quantum states, whereas we have:

**2.18 Proposition.** *0-quantum states are automatically continuous. (In particular, if $X$ is compact, all quantum states for $X$ are continuous.)*

**Proof.** For $m$ to be 0-quantum means that the measure $\nu$ of (2.12) is concentrated on $X_{|\mathfrak{a}}$ (in $\hat{\mathfrak{a}}$), for each $X$-abelian $\mathfrak{a}$. By (1.31), it is equivalent to say that $\nu$ is the image of a measure $\mu$ concentrated on $X_{|\mathfrak{a}}$ (in $\mathfrak{a}^*$). So we have

$$(m \circ \exp_{|\mathfrak{a}})(Z) = \int_{\mathfrak{a}^*} e^{i\langle p, Z\rangle}\, d\mu(p), \qquad (2.19)$$

which shows that $m \circ \exp_{|\mathfrak{a}}$ is continuous (1.28). The continuity of $m$ at $g \in G$ follows by writing $\mathfrak{g}$ as a direct sum of lines $\mathfrak{a}_1, \ldots, \mathfrak{a}_n$ and using the chart $(Z_1, \ldots, Z_n) \mapsto g \exp(Z_1) \cdots \exp(Z_n)$, together with the inequality

$$\big|m(gg_1 \cdots g_n) - m(g)\big| \leq \sqrt{2\mathrm{Re}(1 - m(g_1))} + \cdots + \sqrt{2\mathrm{Re}(1 - m(g_n))} \quad (2.20)$$

which is obtained from (1.18) by induction on $n$. Finally if $X$ is compact then $bX = X$; so quantum $\Leftrightarrow$ 0-quantum, and the above applies. **q.e.d.**

In accordance with (2.5) and (2.9), we have used $X$-abelian subalgebras rather than just abelian ones in Definition 2.17. Let us, however, take note of a situation where this makes no difference:

**2.21 Proposition.** *When $\mathfrak{g}$ is reductive, it is enough to verify the conditions of (2.9) [resp. (2.17)] for abelian subalgebras only.*

**Proof.** Assume that the conditions hold for every abelian subalgebra, and let $\mathfrak{a}$ be $X$-abelian. Since $\mathfrak{g}$ is reductive, we can write $\mathfrak{g} = \mathfrak{g}_0 \oplus \mathfrak{g}_1$ where $\mathfrak{g}_0 = \mathfrak{o}$ and $\mathfrak{g}_1 \simeq \mathfrak{g}/\mathfrak{o}$ is a supplementary ideal. If $Z = Z_0 + Z_1$ accordingly, then (2.10) and (1.19) give $m(\exp(Z)) = m(\exp(Z_1))$. This shows that

$$m \circ \exp_{|\mathfrak{a}} = m \circ \exp_{|\mathfrak{a}_1} \circ \pi_1 \qquad (2.22)$$

where $\pi_1$ denotes the projection of $\mathfrak{a}$ onto its (abelian) image $\mathfrak{a}_1$ in $\mathfrak{g}_1$. Writing $\pi_1^*$ for the dual injection $\hat{\mathfrak{a}}_1 \hookrightarrow \hat{\mathfrak{a}}$, we conclude from (1.30), (1.31a) and our assumption that (2.22) is a state concentrated on $\pi_1^*(X_{|\mathfrak{a}_1}) = X_{|\mathfrak{a}}$ [resp. on $\pi_1^*(bX_{|\mathfrak{a}_1}) = bX_{|\mathfrak{a}}$]. **q.e.d.**

## 2.4 Functoriality

Suppose that $\pi$ is a morphism of $G$ onto another Lie group $H$. Then the dual injection $\pi^* : \mathfrak{h}^* \hookrightarrow \mathfrak{g}^*$ makes each coadjoint orbit $Y$ of $H$ into a coadjoint orbit of $G$; and we should expect a similar relation between quantum states. Indeed, assuming that $\mathrm{Ker}(\pi)$ is connected, we have:

**2.23 Proposition.** *(Notation as above.) The following are equivalent:*
  (a) *$m$ is a [0-]quantum state for $X = \pi^*(Y)$;*
  (b) *$m = n \circ \pi$ where $n$ is a [0-]quantum state for $Y$.*

**Proof.** The relation $X = \pi^*(Y)$ implies that $\operatorname{orth}(X)$ is the preimage of $\operatorname{orth}(Y)$ under the derived map $\dot\pi : \mathfrak{g} \to \mathfrak{h}$. In particular it contains $\operatorname{Ker}(\dot\pi)$. Since $\operatorname{Ker}(\pi)$ is connected, it follows by (2.10) that quantum states for $X$ are trivial on it, and hence are pulled back from states of $H$ (§ 1.2B).

So suppose that $n$ is any state of $H$ and $m = n \circ \pi$, and let $\mathfrak{a} \subset \mathfrak{g}$ project onto $\mathfrak{b} = \dot\pi(\mathfrak{a})$. From $\operatorname{orth}(X) = \dot\pi^{-1}(\operatorname{orth}(Y))$ we see that $\mathfrak{a}$ is $X$-abelian iff $\mathfrak{b}$ is $Y$-abelian. Now in this situation we have

$$m \circ \exp_{|\mathfrak{a}} = n \circ \exp_{|\mathfrak{b}} \circ \dot\pi_{|\mathfrak{a}}. \tag{2.24}$$

Writing $j$ for the injection $\hat{\mathfrak{b}} \hookrightarrow \hat{\mathfrak{a}}$ dual to $\dot\pi_{|\mathfrak{a}}$, we conclude from this and (1.30–31a) that $m \circ \exp_{|\mathfrak{a}}$ is a state concentrated on $X_{|\mathfrak{a}}$ iff $n \circ \exp_{|\mathfrak{b}}$ is a state concentrated on $j^{-1}(X_{|\mathfrak{a}}) = Y_{|\mathfrak{b}}$; and likewise with $X$, $Y$ replaced by $bX$, $bY$.    **q.e.d.**

**2.25 Remark.** This natural property would fail if Definition 2.17 used abelian subalgebras only. Indeed, assume for simplicity that the subgroup $O$ of $G$ generated by $\mathfrak{o} = \operatorname{orth}(X)$ is closed. Taking $H = G/O$ and $Y = X$ regarded as an $H$-orbit in $\operatorname{orth}(\mathfrak{o}) \simeq \mathfrak{h}^*$, we are in the setting of (2.23).

Now, revising the definition changes nothing for $Y$, since $\operatorname{orth}(Y) = \{0\}$ by construction. So if (2.23) remained true, then (2.21) would hold for non-reductive $\mathfrak{g}$ as well. But this is ruled out by the following counterexample.

**2.26 Example: Free fall orbits.** Let $G$ be the simply connected nilpotent group whose Lie algebra $\mathfrak{g}$ has a basis $(e_1, \ldots, e_5)$ with nonvanishing brackets

$$[e_5, e_4] = -e_3, \qquad [e_5, e_3] = -e_2, \qquad [e_4, e_3] = e_1. \tag{2.27}$$

Denote each element of $\mathfrak{g}^*$ by the 5-tuple of its components in the dual basis, and let $X = G(1, 1, 0, 0, c)$. This orbit is identified with $(\mathbf{R}^2, dp \wedge dq)$ by the map $\varPhi$:

$$\varPhi(p, q) = \left(1, 1, p, q, \tfrac{1}{2}p^2 - q + c\right). \tag{2.28}$$

($X$ is a homogeneous symplectic manifold of the Galilei group $G/\operatorname{center}(G)$; it describes a particle of mass 1 and rest energy $c$ in a field of force 1.) Further, let $\bar c > c$ and let $m$ be a 0-quantum state for $\bar X = G(1, 1, 0, 0, \bar c)$ (such states exist (4.34)). We claim that $m \circ \exp_{|\mathfrak{a}}$ is concentrated on $X_{|\mathfrak{a}}$ for each abelian $\mathfrak{a} \subset \mathfrak{g}$, but not for each $X$-abelian $\mathfrak{a}$.

**Proof of claim.** The first part is because we have $\bar{X}_{|\mathfrak{a}} \subset X_{|\mathfrak{a}}$ for each abelian $\mathfrak{a} \subset \mathfrak{g}$. Indeed, suppose otherwise. Then $\mathfrak{a}$ must contain a vector $e_\star$ having a nonzero component along $e_5$. As the commutant of such a vector is spanned by $e_1$, $e_2$ and itself, it suffices to consider the case $\mathfrak{a} = \mathrm{span}(e_1, e_2, e_\star)$ with $e_\star = e_5 + \alpha e_4 + \beta e_3$. But then $\bar{X}_{|\mathfrak{a}} \subset X_{|\mathfrak{a}}$ since (2.28) gives, in the dual basis of $\mathfrak{a}^*$,

$$X_{|\mathfrak{a}} = \begin{cases} \{(1,1,x_\star) : x_\star \in \mathbf{R}\} & \text{when } \alpha \neq 1, \\ \{(1,1,x_\star) : x_\star \geq c - \frac{1}{2}\beta^2\} & \text{when } \alpha = 1. \end{cases} \tag{2.29}$$

The subalgebra $\mathfrak{a} = \mathrm{span}(e_1 - e_2, e_3, e_4 + e_5)$, on the other hand, is both $X$- and $\bar{X}$-abelian; indeed $[\mathfrak{a}, \mathfrak{a}] = \mathbf{R}(e_1 - e_2)$ annihilates both orbits. Moreover (2.28) gives

$$X_{|\mathfrak{a}} = \{(0, p, \tfrac{1}{2}p^2 + c) : p \in \mathbf{R}\}. \tag{2.30}$$

So $X_{|\mathfrak{a}}$ and $\bar{X}_{|\mathfrak{a}}$ are disjoint parabolas, and $m \circ \exp_{|\mathfrak{a}}$ cannot be concentrated on both.

## 2.5 The mean value

Formal differentiation of (2.19) under the integral sign exhibits $\langle \frac{1}{i}Dm(e), Z \rangle$ as the mean value of the variable $\langle x, Z \rangle$ in the state $m$, in the sense of (2.16c). So we may think of the derivative $\frac{1}{i}Dm(e)$, when it exists, as a kind of mean value of $x$. This suggests the

**2.31 Proposition.** *Let $m$ be a 0-quantum state for $X$, twice differentiable at $e$. Then $\frac{1}{i}Dm(e)$ lies in the closed convex hull of $X$.*

**Proof.** Suppose that $X$ lies in a half-space $\{x \in \mathfrak{g}^* : \langle x, Z \rangle \leq c\}$, and put $\mathfrak{a} = \mathbf{R}Z$. Since $m \circ \exp_{|\mathfrak{a}}$ has two derivatives at the origin, we know from [L77, p. 212] that (2.19) may indeed be differentiated under the integral sign. This gives

$$\langle \tfrac{1}{i}Dm(e), Z \rangle = \int_{\mathfrak{a}^*} \langle p, Z \rangle \, d\mu(p) \leq c \tag{2.32}$$

since $\mu$ is concentrated on $X_{|\mathfrak{a}} \subset \{p \in \mathfrak{a}^* : \langle p, Z \rangle \leq c\}$ by hypothesis. Thus, any half-space which contains $X$ also contains $\frac{1}{i}Dm(e)$. But the intersection of these half-spaces is precisely the closed convex hull of $X$.    **q.e.d.**

**2.33 Remarks.** (a) The 'same' fact has been known for some time in specific cases of geometric quantization [W89, A92, N95].

(b) It would fail if $m$ was merely assumed quantum for $X$; this will follow from (4.10c) applied, for instance, to the orbits (2.26).

(c) The states $m$ subject to the hypotheses of (2.31) make a convex set of functions; so their images under $m \mapsto \frac{1}{i}Dm(e)$ make a convex set in $\mathfrak{g}^*$.

> *Any two or more states may be superposed to give a new state.* —P. A. M. Dirac

# 3. Quantum representations

## 3.1 Definition and elementary properties

We continue with $X$ a coadjoint orbit of the Lie group $G$.

**3.1 Definition.** A nonzero unitary $G$-module $\mathcal{H}$, or the representation $U$ of $G$ in it, is *quantum [resp. 0-quantum] for $X$* if, equivalently:
  (a) for each unit vector $\varphi \in \mathcal{H}$, the state $m(g) = (\varphi, g\varphi)$ is quantum [resp. 0-quantum] for $X$;
  (b) for each $X$-abelian subalgebra $\mathfrak{a}$ of $\mathfrak{g}$, $U \circ \exp_{|\mathfrak{a}}$ is a representation concentrated on $bX_{|\mathfrak{a}}$ [resp. on $X_{|\mathfrak{a}}$].

The equivalence (a) ⇔ (b) is clear from (2.9), (2.17) and our definitions in §1.2F. Likewise the properties (2.15, 2.18, 2.21, 2.23) of quantum states have immediate translations, which we leave for the reader to formulate. As a new feature we have

**3.2 Proposition.** *In (3.1) it is actually enough to verify either*
  (a) *condition (3.1a) for one cyclic $\varphi \in \mathcal{H}$; or*
  (b) *condition (3.1b) for one $X$-abelian subalgebra per conjugacy class.*

**Proof.** (a): Suppose that $m(\cdot) = (\varphi, \cdot \varphi)$ is [0-]quantum for $X$, where $\varphi$ is cyclic (i.e., its orbit $G\varphi$ spans a dense subspace of $\mathcal{H}$). Applying (2.23) to inner automorphisms of $G$, we see that

$$m_g(\cdot) = (g\varphi, \cdot g\varphi) \text{ is a [0-]quantum state for } X, \forall g \in G. \qquad (3.3)$$

In particular $m_g$ is trivial on the subgroup $O$ generated by $\mathrm{orth}(X)$ (2.10). Since $\|og\varphi - g\varphi\|^2 = 2\mathrm{Re}(1 - m_g(o))$, it follows that $O$ acts trivially on $G\varphi$ and hence on $\mathcal{H}$. From this we infer as in (2.11) that $U \circ \exp_{|\mathfrak{a}}$ is a *representation* whenever $\mathfrak{a}$ is $X$-abelian. Now (1.32), (3.3) and the cyclicity of $\varphi$ show that this representation is concentrated on $bX_{|\mathfrak{a}}$ [resp. $X_{|\mathfrak{a}}$].

(b): Suppose that the condition holds for $U \circ \exp_{|\mathfrak{a}}$ and that $\mathfrak{b} = \mathrm{Ad}(g)(\mathfrak{a})$. Putting $U = V \circ \pi$ with $\pi(h) = ghg^{-1}$, we deduce exactly as in (2.24) that $V \circ \exp_{|\mathfrak{b}}$, and hence also the unitarily equivalent representation $U \circ \exp_{|\mathfrak{b}}$, is concentrated on $bX_{|\mathfrak{b}}$ [resp. $X_{|\mathfrak{b}}$].  **q.e.d.**

**3.4 Corollary.** *A state $m$ of $G$ is [0-]quantum for $X$ iff the Gel'fand-Naĭmark-Segal module $\mathcal{H}_m$ (1.20) is.*

**3.5 Interpretation.** In the setting of (2.16c), the mean value of $H(x)$ in a state (3.1a) is $(\varphi, E(h)\varphi)$, where $E$ is the spectral measure of $U \circ \exp_{|\mathfrak{a}}$. So we may—again tentatively—regard $E(h)$ as an operator 'quantizing' the function $H = h \circ \pi$. This prompts a few remarks:

(a) In the 0-quantum case, we can replace $\hat{\mathfrak{a}}$ by $\mathfrak{a}^*$ throughout (cf. 2.19) and thus enlarge our class of 'quantizable' functions. (Only *almost periodic* functions on $\mathfrak{a}^*$ factor through $\hat{\mathfrak{a}}$.) But the class remains small if $G$ is!

(b) We always have $E(f(h_1, \ldots, h_n)) = f(E(h_1), \ldots, E(h_n))$, where the right-hand side is defined by the functional calculus of commuting operators. So question (8) of [B91] would be settled.

(c) The fact that $E(h)$ only depends on $h \circ \pi$ reflects computations actually done in physics. Thus if $G$ is the Poincaré group, $\mathfrak{a}$ the space-time translations, and $h$ the quadratic form on $\mathfrak{a}^*$ usually denoted $p_0^2 - p_1^2 - p_2^2 - p_3^2$, then $h \circ \pi$ is a constant usually denoted $m^2$ and we obtain the Klein-Gordon equation, $E(p_0^2) - E(p_1^2) - E(p_2^2) - E(p_3^2) = m^2 \mathbf{1}$; etc.

**3.6 Example: point-orbits.** Suppose that $G$ is connected and $X$ consists of a single point, $\{x\}$. Applying (3.1b) with $\mathfrak{a} = \mathfrak{g}$, we see that any quantum representation must be a multiple of a character $\chi$ given by

$$\chi(\exp(Z)) = e^{i\langle x, Z\rangle}. \tag{3.7}$$

Such a character exists iff $x$ is *integral*, in the sense that $e^{i\langle x, Z\rangle} = 1$ whenever $\exp(Z) = e$; otherwise there is no quantum representation for $X$.

## 3.2 Mackey theory

Suppose that $G$ contains a *normal $X$-abelian subgroup* $A$, i.e., a normal closed connected subgroup whose Lie algebra is $X$-abelian. Then we can draw immediate consequences from (3.1b). Indeed, we have an action of $G$ on $\mathfrak{a}^*$ which makes the projection $\mathfrak{g}^* \to \mathfrak{a}^*$ equivariant, so $X_{|\mathfrak{a}}$ is a $G$-orbit:

$$X_{|\mathfrak{a}} = G(p) = G/H, \tag{3.8}$$

where we have fixed an $x \in X$ and written $H$ for the stabilizer of $p = x_{|\mathfrak{a}}$. (Note that $A \subset H$.) If the spectral measure of $U \circ \exp_{|\mathfrak{a}}$ is concentrated on (3.8), it defines a transitive system of imprimitivity based on $G/H$; so Mackey's little group theorem [F88, p. 1284] applies[†] and gives:

**3.9 Proposition.** *(Notation as above.) If $U$ is a 0-quantum representation for $X$, then there is a unique continuous unitary representation $T$ of $H$ such that*

(a)   $U = \operatorname{Ind}_H^G T,$ (b)   $T \circ \exp_{|\mathfrak{a}} = e^{i\langle p, \cdot \rangle}\mathbf{1}.$

We shall later make an industry out of the fact that this would fail without the crucial prefix '0-'; so it is perhaps appropriate to give the simplest illustration of its role:

**3.10 Example 1: Irrational rotation orbits.** Let $G \simeq \mathbf{Z} \ltimes \mathbf{C}$ (resp. $\mathfrak{g} \simeq \mathbf{C}$) consist of all matrices of the form

$$g = \begin{pmatrix} e^{in} & c \\ 0 & 1 \end{pmatrix}, \qquad \text{resp.} \quad Z = \begin{pmatrix} 0 & \gamma \\ 0 & 0 \end{pmatrix}, \tag{3.11}$$

with $n \in \mathbf{Z}$ and $c, \gamma \in \mathbf{C}$. We identify $\mathfrak{g}^*$ with $\mathbf{C}$ so that $\langle x, Z\rangle = \operatorname{Re}(\bar{x}\gamma)$, and we consider the orbit $X = G(1) = \{e^{in} : n \in \mathbf{Z}\}$. Then (3.9), applied using $\mathfrak{a} = \mathfrak{g}$, asserts that the 0-quantum representations for $X$ are exactly all multiples of $\mathrm{Ind}_{\mathbf{C}}^{G}\bigl(e^{i\mathrm{Re}(\cdot)}\bigr)$, which is the representation acting in $\ell^2(X)$ by

$$(g\varphi)(p) = e^{i\mathrm{Re}(\bar{p}c)}\varphi(e^{-in}p). \tag{3.12}$$

---

[†]Voir aussi l'Annexe B.

Countless other possibilities arise, in contrast, if the representation is only assumed quantum for $X$; for instance (3.12) in $L^2$ of any $G$-invariant measure $\nu$ on $bX$, which is the unit circle; or the same with the right-hand side multiplied by some cocycle $s(n,p)$. Simple examples are $s(n,p) = p^n e^{-in^2/2}$ and $\nu = $ arc length, but there are many others [B84].

Of course, *another* reason why (3.9) rests on the '0-' is that (3.9a) is always continuous, whereas quantum representations need not be. Think, for instance, of the representation of the euclidean group $\mathbf{E}(2)$ in $\ell^2(\mathrm{circle})$ obtained by allowing $n \in \mathbf{R}$ in (3.11) and (3.12).

**3.13 Example 2: Galilean orbits.** Our second example really illustrates a slight extension of (3.9): even when $A$ is not normal, its normalizer, or a closed subgroup $E$ thereof, may happen to be transitive on $X_{|\mathfrak{a}}$. Then (3.1b) implies that the spectral measure of $U \circ \exp_{|\mathfrak{a}}$ gives a transitive system of imprimitivity for $U_{|E}$, which is consequently induced.

This applies when $X$ is a massive particle orbit of the extended Galilei group [S92, §3.1],[†] $A = \{\mathrm{center}\}\times\{\mathrm{boosts}\}$, $E = \{\mathrm{euclidean\ displacements}\}$. Moreover our system of imprimitivity is precisely of the kind postulated in [W62]: $X_{|\mathfrak{a}}$ is $\mathbf{R}^3$ with the usual affine action of $E$. In other words *Definition 3.1 contains Wightman's localizability axioms*, and their consequences, in this case.

It no longer applies when $X$ is a photon orbit [S92]: here $X_{|\mathfrak{a}}$ is still $\mathbf{R}^3$, but the action of $E$ reduces to the linear part, which is not transitive.

(For another example of interest, let $X$ be a twistor orbit of $\mathbf{SU}(2,2)$, $E$ the Poincaré subgroup, and $A$ the space-time translations.)

## 3.3 $\infty$-quantum representations

To proceed from (3.9), we need to know more about the representation $T$. What *ought* to be true is suggested by the following Proposition which is, in effect, the symplectic analogue of (3.9). It is proved and further discussed in the Appendix at the end of this paper.

---

[†]Voir aussi l'Annexe C.

**3.14 Proposition.** *In the setting of (3.8) there is a unique coadjoint orbit $Y$ of $H$, namely $Y = H(x_{|\mathfrak{h}})$, such that*

$$\text{(a)} \quad X = \mathrm{Ind}_H^G Y, \qquad \text{(b)} \quad Y_{|\mathfrak{a}} = \{p\}$$

*(symplectic induction). Moreover $Y$ is also the reduced space $\pi^{-1}(p)/A$, where $\pi$ is the projection $X \to \mathfrak{a}^*$.*

This—especially the last statement—suggests that the $T$ of (3.9) should not only exist, but also in turn be 0-quantum for $Y$. As this is unfortunately not automatic (see (4.38)), we are led to consider the following inductive definition, starting with $n = 0$:

**3.15 Definition.** *A representation $U$ of $G$ is $(n+1)$-quantum for $X$ if it is $n$-quantum for $X$ and, whenever $A$ is a normal $X$-abelian subgroup of $G$, the representation $T$ of (3.9) is $n$-quantum for the orbit $Y$ of (3.14).*

*If this holds for all $n$, we say that $U$ is $\infty$-quantum for $X$.*

We can, of course, introduce $\infty$-quantum states by making the analogue of (3.4) into a definition. Ideally the property should come as the consequence of some more direct definition, perhaps also having a 'compactified' version; but meanwhile (3.15) will suffice for our purposes.

**3.16 Example: Poincaré orbits.** The contents of (3.14–15) are perhaps best illustrated by the classical case of the Poincaré group. Taking for $A$ the subgroup of space-time translations, (3.14) classifies the coadjoint orbits in terms of the possible orbits $X_{|\mathfrak{a}}$ and $Y$, thus as follows:

(a) $X_{|\mathfrak{a}}$ is a timelike hyperboloid and $Y$ an orbit of $\mathbf{SO}(3)$,
(b) $X_{|\mathfrak{a}}$ is a half-cone and $Y$ an orbit of the euclidean group $\mathbf{E}(2)$,
(c) $X_{|\mathfrak{a}}$ is a spacelike hyperboloid and $Y$ an orbit of $\mathbf{SL}(2, \mathbf{R})$,
(d) $X_{|\mathfrak{a}}$ is the origin and $Y$ $(= X)$ an orbit of the Lorentz group.

In each case $Y$ may be a point-orbit: 1-quantum representations (if any) must then, by (3.6) and (3.9), be the 'expected' ones described in [S69]. But $Y$ can also be a spin sphere in case (a), a cylinder in case (b), etc. The cylinder is dealt with by going to 2-quantum if necessary; the sphere is a matter for Chapter 5.

**3.17 Technical remarks on (3.14).** The following observations will be of use in the systematic application of Proposition 3.14. Fig. 1 illustrates them well when $G = \mathbf{E}(2)$, $A =$ translation subgroup, $H = A$.

(a) (3.14a) means that the map $\Phi_{\text{ind}}$ of §1.1D is an isomorphism onto $X$; thus $X$ is the only orbit whose projection meets $Y$ (1.13). In particular, $\mathfrak{h}$ automatically satisfies the *Pukánszky condition*, $x + \text{orth}(\mathfrak{h}) \subset X$.

(b) (3.14b) implies that $A$ acts trivially on $Y$, which is therefore really a homogeneous symplectic manifold of $H/A$; we have, in fact, already taken this into account in Example 3.16.

(c) Introducing the notation $\mathfrak{e}^\perp = \text{orth}(\mathfrak{e}(x))$, we have $\mathfrak{a} \subset \mathfrak{a}^\perp = \mathfrak{h} \supset \mathfrak{h}^\perp$ and
$$\dim(G/H) = \dim(\mathfrak{a}(x)), \qquad \dim(Y) = \dim(\mathfrak{h}/\mathfrak{h}^\perp).^3 \qquad (3.18)$$

Thus, $Y$ will be zero-dimensional just when $H$ happens to be a *polarization* ($\mathfrak{h} = \mathfrak{h}^\perp$).[4] At the other extreme we see that for connected $G$, (3.14) will boil down to $X = \text{Ind}_G^G X$ just when $A$ acts trivially on $X$, i.e., when the $X$-abelian ideal $\mathfrak{a}$ is *$X$-central* (2.6); in particular, such is always the case when $\mathfrak{g}$ is reductive.

**3.19 Synopsis.** This concludes two lengthy chapters of generalities. To sum up, we now have three notions related by

$$\infty\text{-quantum} \quad \Rightarrow \quad 0\text{-quantum} \quad \Rightarrow \quad \text{quantum}. \qquad (3.20)$$

In Chapter 4 we will see that both implications are strict, and that only $\infty$-quantum gives back the Kirillov-Bernat theory of exponential groups. For compact connected groups, on the other hand, all three notions coincide by our remarks in (2.17) and (3.17c); we treat this case in Chapter 5.

---

[3]*Hint*: when $\mathfrak{f}$ normalizes $\mathfrak{e}$, the stabilizer of $x_{|\mathfrak{e}}$ in $\mathfrak{f}$ is $\mathfrak{e}^\perp \cap \mathfrak{f}$; apply this to the pairs $\mathfrak{g}, \mathfrak{a}$ and $\mathfrak{h}, \mathfrak{h}$.

[4]This is how all the polarizations of [S69] arise: see pp. 364–365 and 379, footnote.

*He compared induction to high finance. "If you don't borrow enough, you have cash flow problems. If you borrow too much, you can't pay the interest."* —R. Herb

# 4. Exponential groups

In this Chapter we assume that $G$ is an exponential Lie group. This means that $\exp : \mathfrak{g} \to G$ is a diffeomorphism, or equivalently [B60] that

(a) $G$ is connected, simply connected, and solvable; and
(b) $\mathrm{ad}(Z)$ has no purely imaginary eigenvalues for $Z \in \mathfrak{g}$.

As before $X$ denotes a coadjoint orbit of $G$. Moreover *we fix a point $x \in X$ throughout*, and we write $\xi$ for the function defined on $G$ by

$$\xi(\exp(Z)) = e^{i\langle x, Z\rangle}. \tag{4.1}$$

Following [K62], we say that a closed subgroup $H$ of $G$ is *subordinate to $x$* if $H$ is connected and $\xi_{|H}$ is a character; this means that $x_{|\mathfrak{h}}$ is a point-orbit of $H$, or in other words that $\mathfrak{h} \subset \mathfrak{h}^\perp$ (notation 3.17c).

## 4.1 The Kirillov-Bernat correspondence

**A. Orbits and representations.** Let $\hat{G}$ denote the set of equivalence classes of continuous irreducible unitary representations of $G$, and $\mathfrak{g}^*/G$ the set of coadjoint orbits. Kirillov's method, as extended to exponential groups in [B72], gives a canonical bijection

$$\mathcal{Q} : \mathfrak{g}^*/G \to \hat{G}. \tag{4.2}$$

Using Proposition 3.14 we can give the following 'classical' description of $\mathcal{Q}$. As long as $X$ is not a point-orbit, a lemma of Takenouchi [B72, p. 123] ensures that we can find an $X$-abelian ideal $\mathfrak{a}$ in $\mathfrak{g}$ which is not $X$-central. By (3.14), $X = \mathrm{Ind}_{G_1}^G X_1$ where $G_1$ is the stabilizer of $x_{|\mathfrak{a}}$ and $X_1 = G_1(x_{|\mathfrak{g}_1})$. Moreover $G_1$ has smaller dimension than $G$ (3.17c) and is again exponential [B72, p. 4]. So we may iterate the process to obtain decreasing $G_i$ such that $X = \mathrm{Ind}_{G_1}^G \cdots \mathrm{Ind}_{G_i}^{G_{i-1}} X_i = \mathrm{Ind}_{G_i}^G X_i$ where the dimension of $X_i = G_i(x_{|\mathfrak{g}_i})$ decreases at each step (1.12, 1.14). Ultimately we arrive at a point-orbit of $H = G_n$ say, so that $X = \mathrm{Ind}_H^G \{x_{|\mathfrak{h}}\}$. Then $\xi_{|H}$ is a character (4.1), and we can define

$$\mathcal{Q}(X) = \mathrm{Ind}_H^G \xi_{|H}. \tag{4.3}$$

The theory shows that the right-hand side depends indeed on $X$ only, and, 'quantizing' the above argument, that the $\mathcal{Q}(X)$ exhaust $\hat{G}$.

**4.4 Remarks.** (a) The subgroup $H$ is a polarization satisfying Pukánszky's condition. Indeed $\mathfrak{h} \subset \mathfrak{h}^\perp$ since $\{x_{|\mathfrak{h}}\}$ is a point-orbit, and this inclusion is an equality by dimension (1.12); moreover it follows inductively from (3.17a) that all $\mathfrak{g}_i$ satisfy $x + \mathrm{orth}(\mathfrak{g}_i) \subset X$.

So we have an algorithm to obtain such polarizations. Put algebraically it runs as follows, starting with $\mathfrak{g}_0 = \mathfrak{g}$: if $\mathfrak{g}_i^\perp = \mathfrak{g}_i$, we are done; else, pick an ideal $\mathfrak{a}_i$ of $\mathfrak{g}_i$ such that $\mathfrak{a}_i \subset \mathfrak{a}_i^\perp \not\supset \mathfrak{g}_i$, put $\mathfrak{g}_{i+1} = \mathfrak{g}_i \cap \mathfrak{a}_i^\perp$, and repeat.

This method gives more polarizations than M. Vergne's, but still not all Pukánszky polarizations; for instance it gives the $\mathfrak{h}$ of [B72, p. 88],[5] but not the $\mathfrak{b}_1$ of [B90, p. 313].[6]

(b) We will need the following elementary fact, found e.g. in [A71, I.5.13]. Let $\pi$ be the projection of $G$ onto a quotient group and $\pi^*$ the dual injection, as in (2.23). Then $\mathcal{Q}(\pi^*(Y)) = \mathcal{Q}(Y) \circ \pi$.

**B. Geometric quantization.** We recall how (4.3) fits into the scheme (0.1–2), following [K70, §5.7]. Using the character $\chi = \xi_{|G_x}$ of the stabilizer $G_x$, we can form the associated line bundle $L = G \times_{G_x} \mathbf{C} \xrightarrow{\pi} X$,[7] and endow it with the connection 1-form $\varpi$ whose pull-back to $G \times \mathbf{C}^\times$ is

$$[\,\cdot\,]^*\varpi = \alpha + \frac{dz}{iz} \qquad \text{where} \qquad \alpha = \langle x, g^{-1}dg\rangle. \qquad (4.5)$$

---

[5] In the notation there, choose $\mathfrak{a}_0 = \mathrm{span}(x_3, x_5)$ and then $\mathfrak{a}_1 = \mathrm{span}(x_1, x_3, x_5)$ $(= \mathfrak{h})$.

[6] Benoist's $\mathfrak{b}_1$ contains no noncentral ideal ('$\mathfrak{b}_1 \cap \mathrm{C}_2(\mathfrak{g}) \subset \mathrm{C}_1(\mathfrak{g})$'), whereas our $\mathfrak{h}$ always contains $\mathfrak{a}_0$; in fact, taking orthogonals in the relation $\mathfrak{h} \subset \mathfrak{g}_{i+1} \subset \mathfrak{a}_i^\perp$ shows that we have $\mathfrak{a}_i \subset \mathfrak{h} \subset \mathfrak{a}_i^\perp$ and $\mathfrak{g}_i^\perp \subset \mathfrak{h} \subset \mathfrak{g}_i$ for all $i$.

The curvature $\sigma$, defined by $d\varpi = \pi^*\sigma$, is precisely the orbit's 2-form. Sections of $L$ are the same thing as functions $f : G \to \mathbf{C}$ such that $f(ga) = \xi(a^{-1})f(g)$ for all $a \in G_x$. Since $\mathfrak{h} = \mathfrak{h}^\perp$, the cosets $gH$ project in $X = G/G_x$ as the leaves of a lagrangian foliation, along which a section will be constant just when the covariant derivative $\nabla f = df + if\alpha$ vanishes along the cosets:

$$df(gZ) + i\langle x, Z\rangle f(g) = 0 \qquad \forall Z \in \mathfrak{h}, \qquad (4.6)$$

i.e. just when $f(gh) = \xi(h^{-1})f(g)$ for all $h \in H$. But such functions, square integrable for a quasi-invariant measure on $G/H$, constitute precisely the space of (4.3) in Mackey's realization [F88, p. 1151].

## 4.2 Quantum representations (nilpotent groups)

In this section we determine all quantum representations for $X$ under the additional assumption that $G$ is nilpotent: $\mathrm{ad}(Z)$ is nilpotent for all $Z \in \mathfrak{g}$. Then $X$ is the image of a polynomial map $\mathfrak{g} \to \mathfrak{g}^*$, namely $Z \mapsto \exp(Z)(x)$, and therefore we know from [Z93][†] that $X$ has the same Bohr closure as its affine hull:

**4.7 Theorem.** $bX = b\widehat{X}$, *where $\widehat{X}$ is the affine hull of $X$.*

Thus the compactification can merge an orbit with many others, much as in (3.10); and this shows at once that our 'correspondence principle' (0.4) is quite a bit looser than usual. While this looseness can be remedied by simply not compactifying orbits (see §4.4), examples will show that it is in a sense quite natural.

---

[7]Quotient of $G \times \mathbf{C}$ by the equivalence relation with classes $[g, z] = \{(ga, \overline{\chi}(a)z : a \in G_x\}$; $L$ projects onto $X$ by $\pi([g, z]) = g(x)$.
†Voir aussi l'Annexe D.

Let us first detail the consequences of (4.7). We have $\widehat{X} = x + \operatorname{orth}(\mathfrak{c})$, where $\mathfrak{c}$ is the ideal defined in (1.5); writing $C = \exp(\mathfrak{c})$, we deduce (compare [S69, p. 396]!):

**4.8 Corollary 1.** *A state $m$ of $G$ is quantum for $X$ iff $m_{|C} = \xi_{|C}$.*

**Proof.** The condition means that $m \circ \exp_{|\mathfrak{c}}$ is a state concentrated on $\{x_{|\mathfrak{c}}\}$, which is just $X_{|\mathfrak{c}}$; so its necessity is clear from (2.9). To see the converse, let $\mathfrak{a}$ be any maximal $X$-abelian subalgebra of $\mathfrak{g}$. Then $\mathfrak{a}$ contains $\mathfrak{c}$, and $\widehat{X}_{|\mathfrak{a}}$ is the preimage of $x_{|\mathfrak{c}}$ under the projection $\mathfrak{a}^* \to \mathfrak{c}^*$. Therefore (4.7) gives

$$bX_{|\mathfrak{a}} = b\widehat{X}_{|\mathfrak{a}} = \pi^{-1}(x_{|\mathfrak{c}}), \qquad (4.9)$$

where $\pi$ is the projection $\hat{\mathfrak{a}} \to \hat{\mathfrak{c}}$. So we need only verify that $m \circ \exp_{|\mathfrak{a}}$ is concentrated on $\pi^{-1}(x_{|\mathfrak{c}})$, given that $m \circ \exp_{|\mathfrak{c}}$ is concentrated on $\{x_{|\mathfrak{c}}\}$. But this is clear from (1.30) and (1.31a).    **q.e.d.**

**4.10 Corollary 2.**

(a) *A unitary representation $U$ of $G$ is quantum for $X$ iff $U_{|C} = \xi_{|C}\mathbf{1}$.*
(b) *Such representations correspond one-to-one to all $\gamma$-representations of $G/C$, where $\gamma$ is the Mackey obstruction of $\xi_{|C}$.*
(c) *$\mathcal{Q}(X')$ is quantum for $X$ iff $X'$ lies in the affine hull of $X$.*

**Proof.** (a) is immediate from (4.8) and Definition 3.1a. (b) is a standard reformulation of (a), as we now recall: By definition, the Mackey obstruction of $\xi_{|C}$ is the central extension

$$\gamma : \quad 1 \to \mathbf{T} \to G \times_C \mathbf{T} \to G/C \to 1 \qquad (4.11)$$

where $\mathbf{T} = \mathbf{U}(1)$ and $G \times_C \mathbf{T}$ is the bundle associated to $G \to G/C$ by $\xi_{|C}$, with the group law for which the projection from $G \times \mathbf{T}$ is a morphism; a *$\gamma$-representation* is a projective representation, $V : G/C \to \mathbf{PU}(\mathcal{H})$, which pulls $1 \to \mathbf{T} \to \mathbf{U}(\mathcal{H}) \to \mathbf{PU}(\mathcal{H}) \to 1$ back to (4.11); and the claimed bijection between all such $V$, and all unitary representations $U$ subject to (a), obtains by requiring the commutativity of

$$\begin{array}{ccc} G & \xrightarrow{U} & \mathbf{U}(\mathcal{H}) \\ \downarrow & & \downarrow \\ G/C & \xrightarrow{V} & \mathbf{PU}(\mathcal{H}). \end{array} \qquad (4.12)$$

(c): Suppose that $X' = G(x')$ lies in the affine hull $\widehat{X} = x + \mathrm{orth}(\mathfrak{c})$. Since $G$ stabilizes $x_{|\mathfrak{c}} = x'_{|\mathfrak{c}}$ in $\mathfrak{c}^*$, conversely $C$ stabilizes $x'$ in $\mathfrak{g}^*$. Therefore we infer, *mutatis mutandis* in §4.1B, that $C$ acts in $\mathcal{Q}(X')$ as required by (a):

$$(cf)(g) = f(c^{-1}g) = f(gg^{-1}c^{-1}g) = \xi'(g^{-1}cg)f(g) = \xi(c)f(g). \qquad (4.13)$$

Conversely, suppose that $U = \mathcal{Q}(X')$ is quantum for $X$. Then $U$ is trivial on $O = \exp(\mathfrak{o})$ (2.10) and hence pulled back from $G/O$. By (4.4b) it follows that $X'$ is also pulled back and hence contained in $\mathrm{orth}(\mathfrak{o})$. Since $[\mathfrak{g}, \mathfrak{c}] \subset \mathfrak{o}$ (cf. 1.4–5) we conclude that $X'$ annihilates $[\mathfrak{g}, \mathfrak{c}]$, which is to say that $X'_{|\mathfrak{c}}$ is a single point. But now $U$ is also quantum for $X'$, by the first part of this proof. So $U \circ \exp_{|\mathfrak{c}}$ is concentrated on both points $X_{|\mathfrak{c}}$ and $X'_{|\mathfrak{c}}$, which must therefore coincide.    **q.e.d.**

Note that this Corollary reduces the classification of quantum representations to a 'hopeless' problem, since arbitrary $\gamma$-representations of $G/C$ are in general not type I [H81].

### 4.3 Discrete examples

The conclusion just drawn need not deter us from cultivating the luxuriant category of quantum representations. As they need not be continuous, a simple-minded way to capture new ones suggests itself: to strengthen the topology in Kirillov's construction, for instance by making it discrete. Thus we will replace ordinary induction by discrete induction (§ 1.2D) in (4.3), considering instead

$$\operatorname{ind}_H^G \xi_{|H}, \tag{4.14}$$

which makes sense whenever $H$ is subordinate to $x$. In other words we use the same bundle but counting measure on $G/H$, taking $\ell^2$ rather than $L^2$ sections. The resulting representations are usually discontinuous (cf. 1.24), but still measurable (cf. [H63, p. 348]).

One could as well go on to consider other invariant measures on $G/H$: Radon measures for intermediate topologies (giving e.g. $\operatorname{ind}_K^G \operatorname{Ind}_H^K \xi_{|H}$ for some $K$), Hausdorff measures, etc. But as it is, (4.14) will already make uncountably many irreducibles for most orbits, because (i) independence of polarization breaks down entirely, and (ii) $H$ need not even be a polarization in order that (4.14) be irreducible. In more detail:

**4.15 Theorem.** *Let $G$ be nilpotent, and $H$ subordinate to $x$. Then*
  (a) $\mathcal{D}_H := \operatorname{ind}_H^G \xi_{|H}$ *is quantum for $X$ iff $H$ contains $C$ (cf. 4.8).*
  (b) $\mathcal{D}_H$ *is irreducible iff $H$ is maximal subordinate to $x$.*
  (c) $\mathcal{D}_H$ *and $\mathcal{D}_K$ are inequivalent if $H$ and $K$ are different polarizations.*

**Proof.** (a): We know that $\operatorname{ind}_H^G \xi_{|H}$ is quantum for $X$ iff the state $m$ equal to $\xi_{|H}$ in $H$ and zero outside is quantum for $X$ (1.25, 3.4). According to (4.8), this happens just when $C \subset H$.

(b): Suppose that $H$ is not maximal, i.e., is strictly contained in another subordinate subgroup $K$. Since $K$ is nilpotent, the normalizer $N$ of $H$ in $K$ contains $H$ strictly. Now, given $s \in N - H$, one verifies readily that $(Jf)(g) = f(gs)$ defines a unitary intertwining operator $J : \mathcal{D}_H \to \mathcal{D}_H$ which is not scalar since $(e_e, Je_e) = 0$ (1.22–24). So $\mathcal{D}_H$ is reducible.

Conversely, suppose $\mathcal{D}_H$ is reducible. Then some double coset $D = HaH$, other than $H$, must satisfy the conditions of (1.27) with $\chi = \eta = \xi_{|H}$. But then $a$ must normalize $H$: indeed, if some $h \in H$ were outside $aHa^{-1}$, so would be $h^n$ for all $n \neq 0$; so we would have $h^p aH \neq h^q aH$ whenever $p \neq q$, and so $D/H$ would be infinite, contradicting (1.27b). Thus $a$ normalizes $H$ and stabilizes $\xi_{|H}$ (1.27a). Since the normalizer and stabilizer in question are connected [B60, B72] it follows that $A = \log(a)$ normalizes $\mathfrak{h}$ and stabilizes $x_{|\mathfrak{h}}$. Putting $\mathfrak{k} = \mathfrak{h} \oplus \mathbf{R}A$, we conclude that $\langle x, [\mathfrak{k}, \mathfrak{k}] \rangle = 0$. So $K = \exp(\mathfrak{k})$ is subordinate to $x$, and $H$ is not maximal.

(c): Let $H$ and $K$ be polarizations at $x$, and suppose there is a double coset $D = HaK$ satisfying the conditions of (1.27) with $\chi = \xi_{|H}$, $\eta = \xi_{|K}$. As above, it follows that $H = aKa^{-1}$ and $\chi(h) = \eta(a^{-1}ha)$ for all $h \in H$. Thus we have $e^{i\langle x - a(x), \mathfrak{h} \rangle} = 1$, or in other words, $a(x) \in x + \mathrm{orth}(\mathfrak{h}) = H(x)$ [B72, pp. 69–70]. Since $G_x \subset H$, this forces $a \in H$ and hence $K = H = D$. Consequently, Proposition 1.27 says that $\mathrm{Hom}_G(\mathcal{D}_H, \mathcal{D}_K)$ has dimension 1 if $H = K$, and 0 otherwise.     **q.e.d.**

**4.16 Remark.** It follows from (1.22) and (1.25a) that the cyclic vector $e_e$ of $\mathrm{ind}_H^G \xi_{|H}$ is an eigenvector under $H$: $he_e = \xi(h)e_e$. Thus the resulting state $m(g) = (e_e, ge_e)$ has $m \circ \exp_{|\mathfrak{a}}$ concentrated on $\{x_{|\mathfrak{a}}\}$ for each $X$-abelian subalgebra $\mathfrak{a}$ of $\mathfrak{h}$, so that the statistical interpretation (2.16) leads us to think of it as 'localized' on the preimage of $x_{|\mathfrak{h}}$ in $X$. Note again that this preimage is *coisotropic* because of (1.8) and since $\mathfrak{h} \subset \mathfrak{h}^\perp$.

**4.17 Example 1: Heisenberg's orbit.** The workings of 'dependance on polarization' are already instructive in the commonest case. Namely, let $G$ (resp. $\mathfrak{g}$) consist of all triangular real matrices

$$g = \begin{pmatrix} 1 & b & a \\ & 1 & c \\ & & 1 \end{pmatrix}, \quad \text{resp.} \quad Z = \begin{pmatrix} 0 & \beta & \alpha \\ & 0 & \gamma \\ & & 0 \end{pmatrix}; \tag{4.18}$$

and let $X$ be the orbit of the linear form $x: Z \mapsto -\alpha$. It is isomorphic with $(\mathbf{R}^2, dp \wedge dq)$ under the map $\Phi$ given by $\langle \Phi(p,q), Z \rangle = \begin{vmatrix} p & q \\ \beta & \gamma \end{vmatrix} - \alpha$. Here $C$ is the center $\{b = c = 0\}$, and $\xi(g) = e^{-ia}e^{ibc/2}$.

Suppose at first we use ordinary induction, following the steps of §4.1B. Trivializing the line bundle as $L = X \times \mathbf{C}$ with $\varpi = p\,dq + dz/iz$, sections are just functions $F(p,q)$ (related to the $f$'s of §4.1B by $f(g) = e^{ia}F(b,c)$), on which $G$ acts by

$$(gF)(p,q) = e^{-i(a+bq-bc)}F(p-b, q-c). \tag{4.19}$$

The possible polarizations are the $H_t = \{g : c+bt = 0\}$ or $H_\infty = \{g : b = 0\}$, thus parametrized by $t \in \mathbf{R} \cup \infty$. Inducing from $H_\infty$ means restricting (4.19) to functions of the form $F(p,q) = e^{-ipq}\varphi(p)$,[8] which gives

$$(g\varphi)(p) = e^{-ia}e^{ipc}\varphi(p-b) \qquad (4.20.\infty)$$

in $L^2(\mathbf{R})$. Likewise, inducing from $H_t$ means taking functions of the form $F(p,q) = e^{ip^2t/2}\psi(q+pt)$.[9] Computing the resulting action on $\psi$, or on its Fourier transform $\phi(p) = \sqrt{2\pi}^{-1}\int e^{ipr}\psi(r)\,dr$, we obtain

$$(g\phi)(p) = e^{-ia}e^{ipc}e^{i(bp-\frac{1}{2}b^2)t}\phi(p-b) \qquad (4.20.t)$$

in $L^2(\mathbf{R})$. Clearly (4.20.$\infty$) is equivalent to (4.20.0), and then to (4.20.$t$) as the factor $e^{i(bp-\frac{1}{2}b^2)t}$ is the coboundary, $u(p-b)/u(p)$, of $u(p) = e^{-ip^2t/2}$.

What happens if we use discrete induction instead? All of the above remains valid, provided we regard (4.20.$\infty$) as acting in $\ell^2(\mathbf{R})$ and (4.20.$t$) in $L^2$ of the Bohr compactification $\hat{\mathbf{R}}$ of $\mathbf{R}$ (i.e., in almost periodic functions), using the Fourier transform $\phi(p) = \sum e^{ipr}\psi(r)$ appropriate to $\psi \in \ell^2(\mathbf{R})$. So the inequivalences (4.15c) simply reflect the facts that the constant 1 is not in $\ell^2$, and that the factor $e^{i(bp-\frac{1}{2}b^2)t}$ is no longer a coboundary because $u$ is not almost periodic.

**4.21 Remarks.** (a) In the new version of (4.20.0), the variable $q$ is distributed according to $|\psi(q)|^2$ times counting measure on $\mathbf{R}$, and $p$ according to $|\phi(p)|^2$ times Haar measure on $\hat{\mathbf{R}}$. In particular the state can be 'localized' on a lagrangian submanifold $q = $ const., but then $p$ is uniformly distributed: $|\phi(p)|^2 \equiv 1$. As we shall see in §6.2, this compliance with Heisenberg's uncertainty principle is no accident; it makes, in any case, the appearance of Bohr compactification quite natural.

(b) It may come as a surprise to some readers that the reducible representation (4.19) in $L^2(X)$ is quantum for $X$; indeed the statements can be found in the literature that irreducibility is 'demanded by the uncertainty principle', and that (4.19) 'would violate the uncertainty principle since square integrable sections of $L$ can have arbitrarily small support'. The answer is, of course, that $|F(p,q)|^2$ has nothing to do with the probability densities of $p$ or $q$, which must be computed by the rules (3.5). The reader may then verify that the product of variances $\Delta p\,\Delta q$ *increases* from $\frac{1}{2}$ to $\infty$ when $F$ shrinks to the origin.[10]

---

[8] $\nabla F = dF + iFp\,dq$ vanishes along the leaves $\{p = $ const.$\}$ iff $\partial F/\partial q + ipF = 0$.
[9] $\nabla F$ vanishes along the leaves $\{q + pt = $ const.$\}$ iff $\partial F/\partial p = t(\partial F/\partial q + ipF)$.

(c) Discretely induced representations of this $G$ have already appeared in the literature: When regarded as a projective representation of $G/C = \mathbf{R}^2$ (4.10b), the $\ell^2$ version of (4.19) is just the regular $\gamma$-representation, much-studied as perhaps the simplest example of a type II representation [B73, p. 311]. The $\ell^2$ version of (4.20.$\infty$) was considered in [E81], and is the restriction of a projective representation of $\mathbf{R}\times\hat{\mathbf{R}}$ discussed by Mackey [F88, p. 1202]. Finally (4.20.$t$) appears up to notation in [B84, p. 265], with the sole difference that $G$ there is over the integers, so that $t$ must be a rational (written $d_1/d_2$ in lowest terms) and the representation acts in $2\pi d_2$-periodic instead of almost periodic functions.

**4.22 Example 1 (continued).** We can get yet another uncountable family of irreducible quantum representations by 'discretizing' the Weil-Cartier-Zak models of the Schrödinger representation [C66, Z67]. To describe these, let us fix numbers $B$ and $C$ with $BC = 2\pi$, and apply the Mackey analysis [F88] to the cocompact normal subgroup $\Gamma = \{g \in G : (b,c) \in B\mathbf{Z} \times C\mathbf{Z}\}$. There results that (4.20.$\infty$) is equivalent to $\mathrm{Ind}_\Gamma^G \chi$, where $\chi(g) = e^{-ia}$ ($g \in \Gamma$; note that $\chi \neq \xi_{|\Gamma}$). This representation is naturally realized by restricting (4.19) to the functions which satisfy the Bloch conditions

$$F(p+B, q) = F(p, q), \qquad F(p, q+C) = e^{-ipC}F(p,q) \qquad (4.23)$$

for all $p, q$ and are square integrable over any rectangle of size $B \times C$; the intertwining operator from (4.20.$\infty$) is then $F(p,q) = \sum_{k \in p+B\mathbf{Z}} e^{-ikq}\varphi(k)$.

---

[10]For instance, $F(p,q) = \frac{1}{\sqrt{2\pi\varepsilon}}\, e^{-(p^2+q^2)/4\varepsilon}$ gives $\Delta p\,\Delta q = \frac{1}{2}\sqrt{1+\frac{1}{4\varepsilon^2}}$.

Now consider instead
$$U = \mathrm{ind}_\Gamma^G \chi, \tag{4.24}$$

thus replacing integration by summation. An easy application of (1.27) shows that $U$ is irreducible; moreover, we claim that if $U_M := U \circ \Theta_M$ where $\Theta_M$ is the automorphism of $G$ defined by $M \in \mathbf{SL}(2,\mathbf{R})$:

$$\Theta_M\left(\exp\begin{pmatrix}0 & \beta & \alpha \\ & 0 & \gamma \\ & & 0\end{pmatrix}\right) = \exp\begin{pmatrix}0 & \beta' & \alpha \\ & 0 & \gamma' \\ & & 0\end{pmatrix} \quad \text{with} \quad \begin{pmatrix}\beta' \\ \gamma'\end{pmatrix} = M\begin{pmatrix}\beta \\ \gamma\end{pmatrix}, \tag{4.25}$$

then $U_M \simeq U_N$ iff $M^{-1}N$ has rational coefficients, so that we get a family of irreducible quantum representations parametrized by $\mathbf{SL}(2,\mathbf{R})/\mathbf{SL}(2,\mathbf{Q})$.

**Proof of claim.** We have $U_M = \mathrm{ind}_{\Gamma_M}^G \chi_M$, where

$$\Gamma_M = \Theta_{M^{-1}}(\Gamma) \quad \text{and} \quad \chi_M = \chi \circ \Theta_M \tag{4.26}$$

[F88, p. 1172]. If $M^{-1}N$ is irrational, then every $\Gamma_M$-orbit on the torus $G/\Gamma_N$ is infinite, and so no double coset $\Gamma_M a \Gamma_N$ can satisfy the criterion (1.27b). If it is rational, on the other hand, then the intersection $\Delta$ of $\Gamma_M$ and $\Gamma_N$ has finite index in both, and, as observed in [C66, eq. 73], there is an $s \in G$ such that $\chi_N(\delta) = \chi_M(s^{-1}\delta s)$ for all $\delta \in \Delta$. The double coset $\Gamma_M s^{-1} \Gamma_N$ then gives rise to the explicit intertwining operator

$$(Jf)(g) = \sum_{\gamma\Delta \in \Gamma_N/\Delta} \chi_N(\gamma) f(g\gamma s) \tag{4.27}$$

from $U_M$ to $U_N$ in the realization (1.25). Here as usual the notation implies that each summand only depends on $\gamma$ via its coset $\gamma\Delta$.

**4.28 Example 2: Bargmann's orbit.** The effects of compactification in the previous example were still rather mild, insofar as $X$ was equal to its affine hull (cf. 4.7). So we move on to the next simplest example, where $G$ (resp. $\mathfrak{g}$) consists of all real matrices of the form

$$g = \begin{pmatrix} 1 & b & \frac{1}{2}b^2 & a \\ & 1 & b & c \\ & & 1 & e \\ & & & 1 \end{pmatrix}, \quad \text{resp.} \quad Z = \begin{pmatrix} 0 & \beta & 0 & \alpha \\ & 0 & \beta & \gamma \\ & & 0 & \varepsilon \\ & & & 0 \end{pmatrix}. \tag{4.29}$$

We denote elements of $\mathfrak{g}^*$ as 4-tuples $(m,p,q,h)$, paired to $\mathfrak{g}$ according to $\langle x, Z\rangle = \left|\begin{smallmatrix} p & q \\ \beta & \gamma \end{smallmatrix}\right| - h\varepsilon - m\alpha$, and we consider the orbit of $x = (1,0,0,0)$; it is isomorphic with $(\mathbf{R}^2, dp \wedge dq)$ under the map $\Phi$:

$$\Phi(p,q) = (1, p, q, \tfrac{1}{2}p^2). \tag{4.30}$$

There is now a unique polarization, $H_\infty = \{g : b = 0\}$, giving rise to the quantum representation

$$(g\varphi)(p) = e^{-ia} e^{i(pc - \frac{1}{2}p^2 e)} \varphi(p - b) \tag{4.31}$$

in $L^2(\mathbf{R})$. The statistical interpretation (3.5) is the familiar one: $p$ is distributed according to the measure $|\varphi(p)|^2\, dp$, the pair $(p, h)$ according to its image under $p \mapsto (p, \tfrac{1}{2}p^2)$, and the variable $r := q + pt$ according to $|\psi(r,t)|^2\, dr$ where $\psi(r,t) = \sqrt{2\pi}^{-1} \int e^{-i(pr - \frac{1}{2}p^2 t)} \varphi(p)\, dp$. The functions $\psi$ thus obtained make what has been called *the* Hilbert space of solutions of the Schrödinger equation, $i\psi_t = \tfrac{1}{2}\psi_{rr}$, and $G$ acts on them by

$$(g\psi)(r,t) = e^{-ia} e^{-i\{b(r-c) - \frac{1}{2}b^2(t-e)\}} \psi((r-c) - b(t-e),\, t - e). \tag{4.32}$$

As in (4.17, 4.21a), using discrete induction instead gives us the same picture except that $dp$ is to be replaced throughout by counting measure on $\mathbf{R}$ and $dr$ by Haar measure on $\hat{\mathbf{R}}$. This allows states which are localized in $p$ but then uniformly distributed in $q$, and produces *another* space of solutions of Schrödinger's equation, including e.g. $\psi \equiv 1$.

Things become more interesting if we induce from the subgroup $H_0 = \{g : c = e = 0\}$, which is maximal subordinate but not a polarization. Whereas ordinary induction would then give the representation (4.32) but in $L^2(\mathbf{R}^2)$, where it is is reducible, discrete induction gives it in $\ell^2(\mathbf{R}^2)$, where it is irreducible by (4.15b). Moreover the statistical interpretation of $\psi$ is now quite different: for instance if $\psi$ is the characteristic function of the origin $(0,0)$, so that

$$m(g) = (\psi, g\psi) = \begin{cases} e^{-ia} & \text{if } g \in H_0 \\ 0 & \text{otherwise,} \end{cases} \tag{4.33}$$

one finds that $q$ is distributed according to Dirac measure at the origin, $r = q+pt$ according to Haar measure on $\hat{\mathbf{R}}$ as soon as $t$ is nonzero, and $(p, h)$ according to Haar measure on $\hat{\mathbf{R}}^2$. This is allowed because the parabola $\{(p, \frac{1}{2}p^2) : p \in \mathbf{R}\}$ has the same Bohr closure as the whole plane (cf. 4.7), and corresponds to the fact that the Schrödinger equation has evaporated from (4.32).

Again we will see in Chapter 6 that this uniform distribution of $(p, h)$ is inescapable if we insist on fixing $q$. We can already make this plausible by observing that the state $(\varphi, g\varphi)$ defined by $\varphi(p) = (\frac{\alpha}{2\pi})^{1/4} e^{-\alpha p^2/4}$ in (4.31), gives $\Delta q = \frac{\sqrt{\alpha}}{2}$ and tends to (4.33) pointwise as $\alpha$ tends to zero.

## 4.4 ∞-quantum representations

We now return to general exponential groups, and investigate what happens if we forget about the compactification in (0.4) and consider 0-quantum representations. The cleanest statement is in fact obtained if we go over directly to ∞-quantum (cf. 3.20), for we recover then the Kirillov-Bernat correspondence (4.2):

**4.34 Theorem.** *Let $X$ be a coadjoint orbit of the exponential group $G$, and $U$ a unitary representation of $G$. Then*

$$U \text{ is } \infty\text{-quantum for } X \quad \Leftrightarrow \quad U \text{ is a multiple of } \mathcal{Q}(X). \qquad (4.35)$$

**Proof.** A preliminary observation: Suppose that $\mathfrak{a}$ is an $X$-abelian ideal in $\mathfrak{g}$, and let $Y$ be the unique coadjoint orbit of the stabilizer $H$ of $x_{|\mathfrak{a}}$ such that $X = \mathrm{Ind}_H^G Y$ (3.14). Then it follows from induction in stages and from the *definition* of $\mathcal{Q}$ as we described it, that

$$\mathcal{Q}(X) = \mathrm{Ind}_H^G \mathcal{Q}(Y). \qquad (4.36)$$

$\Leftarrow$: We verify that $\mathcal{Q}(X)$ is $n$-quantum for $X$ for all $n$. First let $n = 0$, and let $\mathfrak{a} \subset \mathfrak{g}$ be $X$-abelian. Then $X_{|\mathfrak{a}}$ consists entirely of point-orbits of the subgroup $A = \exp(\mathfrak{a})$, and H. Fujiwara [F91] has shown that $\mathcal{Q}(X)_{|A}$ is the direct integral of the associated characters, relative to the projection of (the class of) Liouville measure on $X$. Composing with $\exp_{|\mathfrak{a}}$, we conclude that $\mathcal{Q}(X) \circ \exp_{|\mathfrak{a}}$ is concentrated on $X_{|\mathfrak{a}}$.

Now assume inductively that the assertion is true for $n$ for all orbits of all exponential groups. The first condition in the definition of $n+1$-quantum (3.15) is then trivially fulfilled. Moreover, given an $X$-abelian ideal $\mathfrak{a}$ in $\mathfrak{g}$, (4.36) and the uniqueness statement in (3.9) show that with $U = \mathcal{Q}(X)$, the $T$ of (3.9) must be $\mathcal{Q}(Y)$. As this is $n$-quantum for $Y$ by the induction hypothesis, the second condition in (3.15) is also fulfilled.

$\Rightarrow$: We actually prove the stronger assertion in which $\infty$ is replaced by the number $n = \frac{1}{2}\dim(X)$ on the l.h.s. When $n = 0$, it is clear by (3.6). Now assume it inductively for all orbits of dimension less than $2n$ of all exponential groups. Pick an $X$-abelian, non $X$-central ideal $\mathfrak{a} \subset \mathfrak{g}$, and let $Y$ and $T$ be the unique orbit and representation of the stabilizer of $x_{|\mathfrak{a}}$ from which $X$ and $U$ are induced (3.9, 3.14). Since $Y$ has lower dimension than $X$ (3.17c), it follows by the definition of $n$-quantum and the induction hypothesis that $T$ is a multiple of $\mathcal{Q}(Y)$. Therefore we conclude by (4.36) that $U$ is a multiple of $\mathcal{Q}(X)$.   **q.e.d.**

Our next objective is to show that $\infty$ cannot always be replaced by $0$ in (4.35), even for nilpotent groups. The first step however is to observe that it can often be replaced not only by $\frac{1}{2}\dim(X)$ as we have just seen, but in fact by much less. More precisely, let us say that $X$ is *flat* if it is an affine subspace of $\mathfrak{g}^*$, and that $X$ *satisfies Corwin's condition* if there is an $X$-abelian ideal $\mathfrak{a}$ such that $\mathfrak{a}^\perp$ normalizes $\mathfrak{a}^{\perp\perp}$. (This appears as condition (25) in [C90, p. 108].) We have then:

**4.37 Proposition.** *Suppose that $G$ is nilpotent. Then*
  (a) *if $X$ is flat, we can replace $\infty$ by $0$ in (4.35);*
  (b) *if $X$ satisfies Corwin's condition, we can replace $\infty$ by $1$.*

**Proof.** (a) is immediate from (4.10c). Next we recall from [C90, p. 99] that $X$ is flat iff $G$ normalizes the stabilizer $\mathfrak{g}_x = \mathfrak{g}^\perp$. (Indeed this means that $\mathfrak{g}_{x'}$, and hence the tangent space $\mathfrak{g}(x') = \mathrm{orth}(\mathfrak{g}_{x'})$, is constant for $x' \in X$.) Likewise, the orbit $Y$ of (3.14) will be flat just when $H$ normalizes $\mathfrak{h}^\perp$, or, referring to (3.17c), when $\mathfrak{a}^\perp$ normalizes $\mathfrak{a}^{\perp\perp}$.

Thus, assuming Corwin's condition, we can get a flat $Y$ in (3.14). If $U$

is 1-quantum for $X$, we conclude then by (a) that the $T$ of (3.9) is multiple of $\mathcal{Q}(Y)$, and hence by (4.36) that $U$ is a multiple of $\mathcal{Q}(X)$.    **q.e.d.**

The main interest of Corwin's condition is that, as Corwin himself stated in Ref. 16 of [C90], 'there does not appear to be an example in the literature where it is not satisfied.' Thus (4.37b) explains why we know of no nilpotent example where 1-quantum would not imply $\infty$-quantum—although such examples most probably exist. On the other hand we have:

**4.38 Example: 0-quantum $\not\Rightarrow$ 1-quantum.** We consider the simply connected nilpotent group $G$ whose Lie algebra $\mathfrak{g}$ is defined on a basis $(e_1, \ldots, e_6)$ by the nonvanishing brackets

$$
\begin{aligned}
&[e_6, e_5] = -e_4, &&[e_5, e_2] = -e_1, \\
&[e_6, e_4] = -e_3, &&[e_4, e_3] = e_1. \\
&[e_6, e_3] = -e_2,
\end{aligned}
\qquad (4.39)
$$

Let us denote each element of $\mathfrak{g}^*$ by the 6-tuple of its components in the dual basis, and let $X$ be the orbit of $x = (1, 0, 0, 0, 0, c)$. This is identified with $(\mathbf{R}^4, dp \wedge dq + de \wedge df)$ by the map $\Phi$:

$$\Phi(p, q, e, f) = (1, f, p, q, e, \tfrac{1}{2}p^2 - fq + c). \qquad (4.40)$$

Likewise let $\bar{X} = G(1, 0, 0, 0, 0, \bar{c})$ where $\bar{c} > c$. We claim that $\mathcal{Q}(\bar{X})$ is 0-quantum, but not 1-quantum, for $X$.

**Proof of claim.** That $\mathcal{Q}(\bar{X})$ cannot be 1-quantum follows from (4.37b) because $X$ satisfies Corwin's condition, as one sees by taking $\mathfrak{a} = \mathrm{span}(e_1, e_2, e_3)$, getting $\mathfrak{a}^\perp = \mathfrak{a}^{\perp\perp} = \mathrm{span}(e_1, e_2, e_3, e_6)$. On the other hand it is 0-quantum for $X$, because we have $\bar{X}_{|\mathfrak{a}} \subset X_{|\mathfrak{a}}$ for each abelian subalgebra $\mathfrak{a}$ of $\mathfrak{g}$. (Note that $\mathrm{orth}(X) = \mathrm{orth}(\bar{X}) = \{0\}$, so that $X$-abelian means abelian here.) Indeed, suppose that $\mathfrak{a}$ is a maximal abelian subalgebra with $\bar{X}_{|\mathfrak{a}} \neq X_{|\mathfrak{a}}$. Then $\mathfrak{a}$ must contain the center $\mathbf{R}e_1$, and also some vector $e_\star$ having a nonzero component along $e_6$. In other words $\mathfrak{a}$ contains a subalgebra $\mathfrak{a}_0 = \mathrm{span}(e_1, e_\star)$ with, say,

$$e_\star = e_6 + \alpha e_5 + \beta e_4 + \gamma e_3 + \delta e_2. \qquad (4.41)$$

If $\alpha \neq 0$, then $\mathfrak{a}_0$ is its own centralizer, so that we have $\mathfrak{a} = \mathfrak{a}_0$; but then $\bar{X}_{|\mathfrak{a}} = X_{|\mathfrak{a}}$, a straight line in $\mathfrak{a}^*$. Next if $\alpha = 0$, then $\mathfrak{a}_0$ has centralizer $\mathfrak{a}_1 = \mathrm{span}(e_1, e_2, e_\star)$, and we have $\mathfrak{a} = \mathfrak{a}_1$. Now we can assume that $\delta = 0$, and conjugate by $\exp(\beta e_5 + \gamma e_4)$ to reduce to the case $e_\star = e_6$. Then (4.40) gives

$$X_{|\mathfrak{a}} = \{(1, f, \tfrac{1}{2}p^2 - fq + c) : f, p, q \in \mathbf{R}\}. \tag{4.42}$$

This is the plane $\{(1, f, h) : f, h \in \mathbf{R}\}$ with the half-line $\{(1, 0, h) : h < c\}$ removed. In particular if $\bar{c} > c$ then $\bar{X}_{|\mathfrak{a}} \subset X_{|\mathfrak{a}}$, as claimed.

## 4.5 Two digressions

Before we end this Chapter, we digress to discuss two (nonexponential) examples which throw some further light on §4.2. The first one suggests that the phenomenon (4.7) and its corollaries might be true in much greater generality; the second one balances this by showing that condition (0.4) becomes much more stringent for certain 'large' groups.

**A. The metaplectic orbit.** Let $G$ be any covering of $\mathbf{Sp}(2n, \mathbf{R})$, acting on $W = \mathbf{R}^{2n} \setminus \{0\}$ with its standard 2-form $\sigma$, and let $X$ be the image of the resulting momentum map $\Phi$:

$$\langle \Phi(w), Z \rangle = \tfrac{1}{2}\sigma(w, Zw). \tag{4.43}$$

The affine hull of this orbit is all of $\mathfrak{g}^*$; so the analogue of (4.7) reads:

**4.44 Proposition.** $bX = \hat{\mathfrak{g}}$.

**Proof.** We could rely on [Z93], but prefer to give a simplified proof based on the same idea. We consider Gibbs measures $\mu_\mathrm{T} = \Phi(\nu_\mathrm{T})$ on $X$, where

$$\nu_\mathrm{T}(f) = \left(\tfrac{1}{2\pi\mathrm{T}}\right)^n \int_W f(w) e^{-\|w\|^2/2\mathrm{T}} dw, \tag{4.45}$$

and let $\mathrm{T} \to \infty$. Then although $\mu_\mathrm{T}$ has no limit as a measure on $X$ or $\mathfrak{g}^*$, we claim that it does tend vaguely to *Haar measure* $\eta$ on $\hat{\mathfrak{g}}$. Note that the claim implies the Proposition in view of (1.31a) and (1.31c).

But to prove the claim it suffices to observe that the Fourier transform (2.2) of $\mu_\mathrm{T}$:

$$\hat{\mu}_\mathrm{T}(Z) = \prod_{j=1}^{2n} \frac{1}{\sqrt{1 - 2i\mathrm{T}\zeta_j}} \tag{4.46}$$

($\zeta_j$ = coefficients of the quadratic form $\langle \Phi(\cdot), Z \rangle$ in principal axes) tends pointwise to the characteristic function of $0 \in \mathfrak{g}$, which is $\hat{\eta}$. Thus we have $\mu_\mathrm{T}(f) \to \eta(f)$ whenever $f$ is a character of $\hat{\mathfrak{g}}$, and hence also for any continuous $f$ since linear combinations of characters are uniformly dense in the continuous functions on $\hat{\mathfrak{g}}$.     **q.e.d.**

It is tautological from (4.44) that *any* probability measure on $\hat{\mathfrak{g}}$ is concentrated on $bX$, and by the same token, that *any* unitary representation of $G$ is quantum for $X$ (3.1b). But at the same time the proof illustrates in what sense even Haar measure on $\hat{\mathfrak{g}}$ can arise from 'physics on $X$': namely, for instance as an 'infinite temperature state'. This is the kind of objects that the compactification in (0.4) allows.

It would be interesting to know if all noncompact orbits of simple Lie groups enjoy the density property (4.44); this can be verified for $\mathbf{SL}(2,\mathbf{R})$.

**B. The prequantum representation in $L^2(\mathbf{R}^2)$.** Whenever (4.7) holds, our principle (0.4) will degenerate to its mere application to *constant* hamiltonians (cf. 4.10a). Yet, we should emphasize, this does not always happen. Indeed counterexamples occur in compact groups (cf. 2.18), or e.g. in $\mathbf{E}(2)$ (Fig. 1). Moreover they will occur in suitably chosen infinite-dimensional groups, such as the following.

To describe it, we return to the plane $X = (\mathbf{R}^2, dp \wedge dq)$ of (4.17) and let $G$ be the full automorphism group of the line bundle $L = X \times \mathbf{C}$ over it. Explicitly $G$ consists of all diffeomorphisms $g$ of the form

$$g(p,q,z) = (s(p,q), z e^{iS(p,q)}) \qquad (4.47)$$

with $dS = p\,dq - s^*(p\,dq)$. Its 'Lie algebra' of vector fields, $\mathfrak{g}$, is isomorphic with $C^\infty(X)$ (1.2) under $H \mapsto \hat{\eta}$:

$$\hat{\eta}(p,q,z) = (\eta(p,q), iz\ell(p,q)) \qquad (4.48)$$

where $\eta = \mathrm{drag}\,H = (-\partial H/\partial q, \partial H/\partial p)$ and $\ell = H - p\partial H/\partial p$ [L90, p. 270]. We can regard $X$ as one of its coadjoint orbits by letting $(p,q)$ stand for the linear form $H \mapsto H(p,q)$, and make sense of Definition 2.3 here by restricting attention to *complete* commuting $\hat{\eta}_j$. (Such infinite-dimensional cases are, in fact, included in the original definition of [S88].)

Consider now the natural representation of $G$ in square integrable sections of $L$. Regarding sections as functions $F(p,q)$ as in (4.19), it is given in $L^2(X)$ by

$$(gF)(p,q) = e^{iS(s^{-1}(p,q))}F(s^{-1}(p,q)). \tag{4.49}$$

It is true by construction that this represents the center correctly, i.e., the flow $e^{t\hat{\eta}}(p,q,z) = (p,q,ze^{ita})$ of a constant $H \equiv a$ is represented by $e^{ita}\mathbf{1}$. Nevertheless,

**4.50 Proposition.** *The representation (4.49) of $G$ is not quantum for $X$.*

**Proof.** We consider the function $H(p,q) = \sin p$. It gives rise to an infinitesimal automorphism (4.48) whose flow writes

$$e^{t\hat{\eta}}(p,q,z) = (p, q+t\cos p, ze^{it(\sin p - p\cos p)}). \tag{4.51}$$

Computing the resulting action (4.49) of this one-parameter subgroup on $F$ and on its partial Fourier transform $G(p,k) = \sqrt{2\pi}^{-1}\int e^{ikq}F(p,q)\,dq$, we obtain

$$(e^{t\hat{\eta}}G)(p,k) = e^{it\{\sin p + (k-p)\cos p\}}G(p,k) \tag{4.52}$$

in $L^2(\mathbf{R}^2)$. This shows that the spectral measure of $(F, e^{t\hat{\eta}}F)$ is the image of $|G(p,k)|^2\,dp\,dk$ under $(p,k) \mapsto \sin p + (k-p)\cos p$. Now if (4.49) were quantum for $X$, this measure would be always concentrated on the range $[-1,1]$ of $H$ (2.9); but this is clearly not the case. **q.e.d.**

**4.53 Remark.** It may seem bold to wonder about quantum representations of $G$ for $X$. Yet we should emphasize that Van Hove's Theorem [V51, p. 86] does not forbid the existence of such things: rather it says that—under certain technical hypotheses—they must be *reducible* when restricted to the Heisenberg subgroup defined by the action $g(p,q,z) = (p+b, q+c, ze^{-i(a+bq)})$ of (4.18) on $L$. But this need not worry us, since we have seen that quantum representations of the Heisenberg group could very well be reducible (4.21b); we understand the words 'too large' in (0.1) as referring to spectra rather than reducibility. (Of course, such remarks will remain rhetorical until someone comes up with a quantum representation of $G$ for $X$.)

*Three quarks for muster Mark!*
—J. Joyce

# 5. Compact groups

In this Chapter, $G$ is a compact connected Lie group. We fix a maximal torus $T \subset G$ throughout, and we write $W$ for the resulting Weyl group, $W = \mathrm{Normalizer}(T)/T$. It is finite and acts on $\mathfrak{t}$ and $\mathfrak{t}^*$ by conjugation.

In view of (2.18), we restrict attention to (unitary) $G$-modules $\mathcal{H}$ in which the action of $G$ is *continuous*.

## 5.1 The Borel-Weil correspondence

**A. Orbits and representations.** Let us establish some standard notation. We will say that a form $\mu \in \mathfrak{t}^*$ is *integral* if the character $e^{i\langle\mu,\cdot\rangle}$ of $\mathfrak{t}$ factors through the covering $\exp : \mathfrak{t} \to T$, and is a *weight of* $\mathcal{H}$ if the resulting character of $T$ occurs in the decomposition of $\mathcal{H}$ as a $T$-module; we write this decomposition

$$\mathcal{H} = \bigoplus\nolimits_{\mu \,:\, \mathrm{weight}} \mathcal{H}^\mu. \qquad (5.1)$$

Moreover, following [B60, p. IX.66], we let

(5.2) R consist of the nonzero weights of $\mathfrak{g}_{\mathbf{C}}$ (adjoint action);
(5.3) C be a chosen connected component of $\mathfrak{t} \setminus \bigcup_{\alpha \in \mathrm{R}} \mathrm{Ker}(\alpha)$;
(5.4) $\leq$ be defined on $\mathfrak{t}^*$ by: $\lambda \leq \mu$ iff $\mu - \lambda$ is nonnegative on $\overline{\mathrm{C}}$;
(5.5) D consist of those $\mu \in \mathfrak{t}^*$ such that $w(\mu) \leq \mu$ for all $w \in W$.

(R, C, D stand for *roots, chamber, dominant*.) Then $\overline{\mathrm{C}}$ and D are fundamental domains for the action of $W$ on $\mathfrak{t}$ and $\mathfrak{t}^*$, and we have (ibid., pp. 67 and 14):

**5.6 Theorem (Cartan-Weyl).**

(a) *Every irreducible G-module $\mathcal{H}$ has a unique $\leq$-maximal weight $\lambda$, which characterizes $\mathcal{H}$ and can be any integral point of D.*
(b) *Every coadjoint orbit $X$ of $G$ intersects $\mathfrak{t}^*$ in a $W$-orbit, and hence D in a unique point $\mu$.*

In (b), we have assumed the canonical inclusion $\mathfrak{t}^* \hookrightarrow \mathfrak{g}^*$ which arises as follows: being maximal abelian, $\mathfrak{t}$ coincides with the space of all $T$-fixed points in $\mathfrak{g}$; whence a canonical projection, $\int_T \mathrm{Ad}(t)\, dt : \mathfrak{g} \to \mathfrak{t}$, whose transpose identifies $\mathfrak{t}^*$ with the $T$-fixed points of $\mathfrak{g}^*$.

Together (a) and (b) give a bijection

$$\mathcal{Q} : (\mathfrak{g}^*/G)_{\mathrm{int}} \to \hat{G} \qquad (5.7)$$

between all *integral* coadjoint orbits (those through integral points of D) and all equivalence classes of irreducible $G$-modules. This gives an intrinsic reformulation of (5.6a), independant of the choice of $T$ and C.

We recall below how geometric quantization realizes the map $\mathcal{Q}$ (without, however, giving an independant proof of its existence).

**B. Prequantization.** Let $X$ and $\mathcal{H}$ correspond under (5.7): $\lambda = \mu$. Following [B54, T55] we describe a projective embedding $X \to \mathbf{P}(\mathcal{H})$ which will exhibit $X$ both as a complex manifold and as the base of a line bundle (0.1). To this end we consider the point $\boldsymbol{\lambda} \in \mathbf{P}(\mathcal{H})$ defined by the maximal weight space $\mathcal{H}^\lambda$ (which is 1-dimensional [B60]) and recall the

**5.8 Theorem.** *The orbit $\boldsymbol{X} := G(\boldsymbol{\lambda})$ is a complex submanifold of $\mathbf{P}(\mathcal{H})$.*

**Proof.** We can replace $G$ and $\mathfrak{g}$ by their images in $\mathrm{End}(\mathcal{H})$, and consider the action of $G_{\mathbf{C}} = \exp(\mathfrak{g} \oplus i\mathfrak{g})$ on $\mathbf{P}(\mathcal{H})$. Let us show that $\mathfrak{g}_{\mathbf{C}}(\boldsymbol{\lambda}) = \mathfrak{g}(\boldsymbol{\lambda})$, so that our compact orbit is both open and closed in the complex orbit.

We have $\mathfrak{g}_{\mathbf{C}} = \bigoplus_{\alpha \in \mathrm{R} \cup \{0\}} \mathfrak{g}^\alpha$ with $\mathfrak{g}^\alpha = \{Z : [A, Z] = i\langle \alpha, A\rangle Z$ for $A \in \mathfrak{t}\}$, and each $\alpha \in \mathrm{R}$ is either $> 0$ or $< 0$ since it cannot change sign on C (5.3). If $\alpha > 0$, then $\mathfrak{g}^\alpha(\boldsymbol{\lambda}) = 0$; indeed the relation $AZ\varphi = [A, Z]\varphi + ZA\varphi$ shows that $\mathfrak{g}^\alpha \mathcal{H}^\lambda \subset \mathcal{H}^{\lambda+\alpha}$, which is zero since $\lambda$ is maximal. Likewise $\mathfrak{g}^0(\boldsymbol{\lambda}) = 0$ since $\mathfrak{g}^0$ stabilizes $\mathcal{H}^\lambda$. Finally, suppose that $Z \in \mathfrak{g}^\alpha$, $\alpha < 0$; then one verifies that the adjoint $Z^*$ is in $\mathfrak{g}^{-\alpha}$ and so acts trivially; thus $Z$ has the same effect on $\boldsymbol{\lambda}$ as $Z - Z^*$ which, being skew-adjoint, is in $\mathfrak{g}$.    **q.e.d.**

The orbit $\boldsymbol{X}$ is known as the manifold of *coherent states* [O75].

To complete the construction, we regard $\mathbf{P}(\mathcal{H})$ as consisting of all rank 1 self-adjoint projectors $\boldsymbol{p}$ in $\mathcal{H}$, and we consider over it the 'tautological' line bundle $\mathcal{L} = \{(\boldsymbol{p}, \varphi) : \varphi \in \mathrm{Im}(\boldsymbol{p})\}$. It is associated with $\mathcal{L}^\times = \mathcal{H} - \{0\}$, and admits the connection

$$T_\varphi \mathcal{L}^\times = \mathcal{H} = \mathrm{Ker}(\boldsymbol{p}) \oplus \mathrm{Im}(\boldsymbol{p}). \qquad (5.9)$$

Writing 'hor' and 'vert' for the projections on these summands, $\pi$ for the projection $\mathcal{L}^\times \to \mathbf{P}(\mathcal{H})$, and $\varpi$ and $\sigma$ for the connection and curvature forms (defined by $\mathrm{vert}(\delta\varphi) = i\varphi\varpi(\delta\varphi)$ and $d\varpi = \pi^*\sigma$), we find

$$\varpi(\delta\varphi) = \frac{(\varphi, \delta\varphi)}{i\|\varphi\|^2}, \qquad \sigma(\delta\boldsymbol{p}, \delta'\boldsymbol{p}) = \mathrm{Tr}(\delta'\boldsymbol{p} J \delta\boldsymbol{p}), \qquad (5.10)$$

where $J$ is the complex structure of $\mathbf{P}(\mathcal{H})$: $J\delta\boldsymbol{p} = \frac{1}{i}[\boldsymbol{p}, \delta\boldsymbol{p}]$. Note that $\sigma$ is nondegenerate on any complex submanifold, for $J$ is *negative* relative to $\sigma$ in the sense that $\sigma(J\delta\boldsymbol{p}, \delta\boldsymbol{p}) < 0$ for any nonzero tangent vector $\delta\boldsymbol{p}$. In particular $\boldsymbol{X}$ is homogeneous symplectic, so (1.7) asserts that the resulting momentum map, computed as in (1.6):

$$\langle \Phi(\boldsymbol{x}), Z \rangle = \varpi(Z\xi) = \frac{(\xi, Z\xi)}{i\|\xi\|^2}, \qquad \xi \in \mathrm{Im}(\boldsymbol{x}), \qquad (5.11)$$

covers a coadjoint orbit. But this orbit is $X$ because we have $\Phi(\boldsymbol{\lambda}) = \lambda$ (equality in $\mathfrak{t}^* \subset \mathfrak{g}^*$), and is simply connected because (co)adjoint stabilizers are connected [B60, p. IX.10]. So $\Phi$ is a symplectic diffeomorphism $\boldsymbol{X} \to X$, whose inverse is the sought embedding, $x \mapsto \boldsymbol{x}$.

Pulling $\mathcal{L}$ and $J$ back to $X$ now provides, canonically, a line bundle with connection $L \to X$ whose curvature is the orbit's 2-form, and a negative $G$-invariant complex structure on $X$. We continue to write $J$ for it, and $\varpi$ for the connection form (5.10) on $L^\times$, which is thus naturally a complex submanifold of $\mathcal{H} - \{0\}$.

**C. Polarization.** Sections of $L$ can be viewed either as maps $s: X \to L$ such that $s(x) \in L_x$, or as functions $f: L^\times \to \mathbf{C}$ such that $f(z\xi) = \overline{z}f(\xi)$, the correspondence between the two being fixed by

$$s(x) = \frac{\xi f(\xi)}{\|\xi\|^2}, \qquad \xi \in L_x^\times. \qquad (5.12)$$

Square integrable sections for Liouville measure on $X$ make a unitary $G$-module $H(L)$ with $(gf)(\xi) = f(g^{-1}\xi)$. Writing $H^\bullet(L)$ for the submodule in which, equivalently,

(a) $f$ is antiholomorphic,
(b) $\nabla f = df \circ \text{hor} = df + if\overline{\varpi}$ is complex antilinear,
(c) $\nabla s = \text{vert} \circ Ds$ vanishes on all $+i$-eigenvectors of $J$ in $TX^{\mathbf{C}}$,

we have the

**5.13 Theorem (Borel-Weil).** $H^\bullet(L)$ *is isomorphic to* $\mathcal{H}$.

**Proof.** Sending $\varphi \in \mathcal{H}$ to the section given by $f(\xi) = (\xi, \varphi)$ gives a nonzero intertwining operator; so we need only see that $H^\bullet(L)$ is irreducible.

Now [B48, p. 117] shows that $H^\bullet(L)$ is a closed subspace of $H(L)$ and that the evaluation functionals $f \mapsto f(\xi)$ are bounded on it. So the reproducing kernel argument of [K61] applies and gives the conclusion.     **q.e.d.**

**5.14 Remark.** Using (5.13) to actually construct $\mathcal{H}$ requires an independent description of the ingredients $(L, J)$ which define $H^\bullet(L)$.

(a) Concerning $L$, one observes that the isotropy action of $G_\lambda$ in $L_\lambda$ defines a character $\chi$ of $G_\lambda$ whose differential is $i\lambda_{|\mathfrak{g}_\lambda}$, as (5.11) shows. Thus we know $\chi$ and can form the associated bundle with connection $G \times_{G_\lambda} \mathbf{C}$ as in (4.5); and this is isomorphic to $L$ under $[g, z] \mapsto (g(\boldsymbol{\lambda}), gz\xi_0)$, where $\xi_0$ is any given unit vector in $L_\lambda$.

(b) As for $J$, it is determined by its restriction to $T_\lambda X$, which is $T$-invariant and so must preserve the decomposition of $T_\lambda X$ as a $T$-module. Now the proof of (5.8) exhibits this as $\cong \bigoplus_{\alpha<0} \mathfrak{g}^\alpha(\boldsymbol{\lambda})$, where the nonzero summands are 2-planes because $\dim_{\mathbf{C}}(\mathfrak{g}^\alpha) = 1$ [B60]. So we need only pick in each plane, of the two existing rotation-invariant complex structures, the one which is negative relative to the orbit's 2-form. (Cf. [B58, p. 498].)

## 5.2 Quantum representations

Let $\mathcal{H}(\lambda)$ denote the irreducible module with maximal weight $\lambda$, and $X(\mu)$ the coadjoint orbit with dominant element $\mu$ (5.6). We have:

**5.16 Theorem.** $\mathcal{H}(\lambda)$ *is quantum for* $X(\mu)$ *iff* $\lambda \leq \mu$.

**Proof.** Since all maximal tori of $G$ are conjugate [B60], it suffices to apply our criterion (3.1b) with $\mathfrak{a} = \mathfrak{t}$; cf. (3.2b, 2.21, 2.15). Moreover, we know from Kostant's convexity theorem that $X(\mu)_{|\mathfrak{t}}$ is the convex hull of the Weyl group orbit $W(\mu)$ [Z92].[†] Thus $\mathcal{H}(\lambda)$ will be quantum for $X(\mu)$ iff

$$\{\text{weights of } \mathcal{H}(\lambda)\} \subset \text{conv}(W(\mu)). \tag{5.17}$$

If (5.17) fails, then some weight $\kappa$ is separated from $W(\mu)$ by a hyperplane; i.e., we can find $Z \in \mathfrak{t}$ such that $\langle w^{-1}(\mu) - \kappa, Z\rangle = \langle \mu - w(\kappa), w(Z)\rangle < 0$ for all $w \in W$. Fixing $w$ so that $w(Z) \in \overline{C}$ this says that the weight $w(\kappa)$ is not $\leq \mu$, whence $\lambda \not\leq \mu$.

Conversely if $\lambda \not\leq \mu$, then some $Z \in \overline{C}$ satisfies $\langle \mu - \lambda, Z\rangle < 0$, whence also $\langle w(\mu) - \lambda, Z\rangle < 0$ for all $w \in W$ (5.4, 5.5). So $\lambda$ is separated from $W(\mu)$ by a hyperplane, and (5.17) fails.     **q.e.d.**

**5.18 Remark.** One can also prove that $\mathcal{H}(\lambda)$ is quantum for $X(\lambda)$ without using the convexity theorem: Given a weight $\kappa$, with eigenprojector $\boldsymbol{E}$ say, we must find a point of $X(\lambda)$ above $\kappa$; and this can be done by maximizing the transition probability $\varrho(\boldsymbol{x}) = \text{Tr}(\boldsymbol{Ex})$ over the orbit (5.8).

This Rayleigh-Ritz method (adapted from [G82, thm 5.3]) can actually be used to prove the convexity theorem itself; see [Z92] for more details.

**5.19 Example 1: Fundamental orbits of U(n).** Let $\mathcal{E}$ denote the Hilbert space $\mathbf{C}^n$, $G$ its unitary group, and identify $\mathfrak{g}^*$ with the self-adjoint matrices by writing $\langle x, Z \rangle = \frac{1}{i}\text{Tr}(xZ)$. Thus an orbit is a conjugacy class of such matrices; it is integral iff their spectrum consists of integers. We consider the Grassmannian of all rank $k$ self-adjoint projectors, i.e. the orbit $X = X(\pi_k)$ with

$$\pi_k = \begin{pmatrix} \mathbf{1} & \mathbf{0} \\ \mathbf{0} & \mathbf{0} \end{pmatrix} \begin{matrix} \}k \\ \}n-k \end{matrix}. \tag{5.20}$$

---

[†]Voir aussi l'Annexe E.

This turns out to be a *minimal* dominant integral form in the order (5.4).[11] So (5.16) asserts that there is a unique irreducible quantum module for $X$, namely $\mathcal{H} = \mathcal{H}(\pi_k)$; this, as one knows, is just the $k$th exterior power $\wedge^k \mathcal{E}$, the embedding $X \to \mathbf{P}(\mathcal{H})$ being the Plücker map $x \mapsto \wedge^k(x)$.

In particular for $k = 1$ we see that the only irreducible quantum module for $\mathbf{P}(\mathcal{E})$ is $\mathcal{E}$ itself—in keeping with the idea that 'quantum mechanics in $\mathcal{E}$' should correspond to 'classical mechanics on $\mathbf{P}(\mathcal{E})$'.

**5.21 Example 2: Hadron orbits.** Specializing to $\mathbf{U}(3)$, Fig. 2 confirms pictorially that the triangle with vertices

$$u = \begin{pmatrix} 1 & & \\ & 0 & \\ & & 0 \end{pmatrix}, \qquad d = \begin{pmatrix} 0 & & \\ & 1 & \\ & & 0 \end{pmatrix}, \qquad s = \begin{pmatrix} 0 & & \\ & 0 & \\ & & 1 \end{pmatrix}, \tag{5.22}$$

encloses no other integral points, so that $\mathbf{C}^3$ itself is the only irreducible quantum module for $X(u)$. On the other hand the projection of $X(2u - s)$ contains exactly three integral $W$-orbits, giving rise to three quantum modules: $\mathcal{H}(2u - s)$, but also $\mathcal{H}(u)$ and $\mathcal{H}(u + d - s)$.

---

[11]In $\mathfrak{u}(n)$, we let $\mathfrak{t}$ be the diagonal matrices and $C = \{\text{diag}(i\vartheta_1, \ldots, i\vartheta_n) : \vartheta_1 > \ldots > \vartheta_n\}$; the Weyl group is the symmetric group acting by permutation of the diagonal coefficients and we have $D = \{\text{diag}(\lambda_1, \ldots, \lambda_n) : \lambda_1 \geq \ldots \geq \lambda_n\}$, the order (5.4) being given by

$$\lambda \leq \mu \quad \Leftrightarrow \quad \begin{cases} \lambda_1 + \cdots + \lambda_j \leq \mu_1 + \cdots + \mu_j & \forall j, \\ \lambda_1 + \cdots + \lambda_n = \mu_1 + \cdots + \mu_n. \end{cases}$$

When $\mu = \pi_k$ this implies $\lambda_1 \leq 1$, $\lambda_n \geq 0$ and $\sum \lambda_i = k$, whence necessarily $\lambda = \pi_k$ if the $\lambda_i$ are to be decreasing integers.

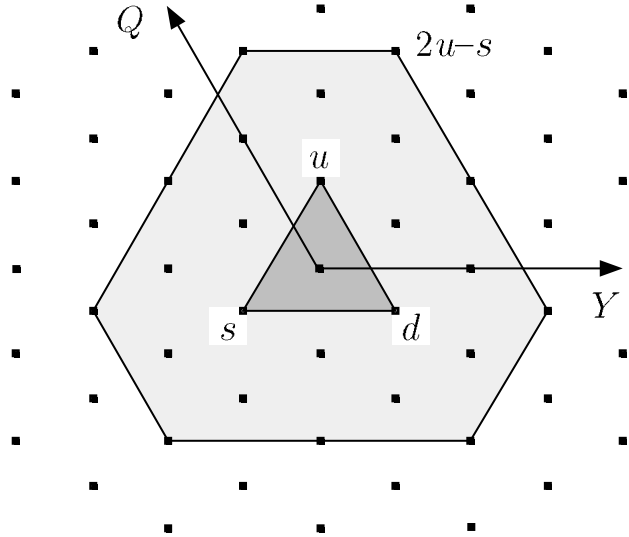

**Fig. 2** U(3)-orbits projected in the dual of the diagonal subalgebra.

## 5.3 Quantum representations within $H(L)$

The foregoing shows that our principle (0.4) is not enough to fix a single representation for each orbit $X$, except in minimal cases like (5.19); but the picture changes drastically if we borrow an idea from (0.1–2) and, assuming that $X$ is integral, restrict attention to submodules of $H(L)$. Indeed:

**5.23 Theorem.** *The only quantum submodule of $H(L)$ is $H^{\bullet}(L)$. In other words, the state $m(g) = (f, gf)$ defined by a normalized section (5.12) of $L$ is quantum for $X$ iff $f$ is antiholomorphic.*

**Proof.** Write $\lambda$ for the dominant element of $X$, and recall the identification (5.14a) of $L$ with the associated bundle $G \times_{G_\lambda} \mathbf{C}$. In terms of sections this means that
$$H(L) = \mathrm{Ind}_{G_\lambda}^G \chi \tag{5.24}$$
where $\chi(\exp(Z)) = e^{i\langle \lambda, Z\rangle}$. By Frobenius reciprocity, [F88], the irreducible module with maximal weight $\nu$ occurs in (5.24) with multiplicity
$$\text{mult. of } \mathcal{H}(\nu) \text{ in } H(L) = \text{mult. of } \chi \text{ in } \mathcal{H}(\nu)_{|G_\lambda}. \tag{5.25}$$

Since $G_\lambda$ contains $T$, the right-hand side is bounded above by the multiplicity of $\lambda$ as a weight of $\mathcal{H}(\nu)$, which is certainly zero unless $\lambda \leq \nu$. On the other hand we know that $\mathcal{H}(\nu)$ is not quantum for $X$ unless $\nu \leq \lambda$. So the only quantum summands in the decomposition of (5.24) are those isomorphic to $\mathcal{H}(\lambda)$; but since the multiplicity of $\lambda$ in $\mathcal{H}(\lambda)$ is 1 there is just one such summand, which (5.13) identifies as $H^{\bullet}(L)$.    **q.e.d.**

The question of why one should reject all summands of (5.24) except the lowest one was already raised in [S66, pp. 333, 337] in the case of $\mathbf{SU}(2)$; see also [S69, pp. 374–375]. We now see that (0.4) provides just the right principle for doing this; it explains, in a way, the polarization step of geometric quantization.

## 5.4 Another characterization of $\mathcal{Q}(X)$

While (5.23) nicely motivates the restriction from $H(L)$ to $H^{\bullet}(L)$, it does not say why we should be looking at sections of $L$ in the first place. This idea is quite heterogeneous to what we do, and one might prefer to deduce the bound $\lambda \geq \mu$ opposite to (5.16) from something 'opposite to (0.3)'.

We note here—as a curiosity—that an inequality considered by B. Simon [S80] precisely fits this purpose. Assuming $\mathcal{H}$ finite-dimensional it reads

$$\tfrac{1}{\dim(\mathcal{H})}\mathrm{Tr}(e^{-\boldsymbol{H}_Z}) \;\geq\; \tfrac{1}{\mathrm{vol}(X)}\int_X e^{-H_Z(x)}\,dx \qquad (5.26)$$

where $\boldsymbol{H}_Z$ is the self-adjoint generator of $\exp(tZ)$ acting on $\mathcal{H}$, $H_Z$ is the hamiltonian function (2.1), and $dx$ is Liouville measure. We have then:

**5.27 Proposition.** *If $X$ is integral, then there is a unique irreducible quantum module for $X$ which satisfies (5.26) for all $Z \in \mathfrak{g}$, namely $\mathcal{H} = \mathcal{Q}(X)$. Otherwise, there is no such module.*

**Proof.** Write $\mathcal{H} = \mathcal{H}(\lambda)$ and $X = X(\mu)$ as before. When $\lambda = \mu$, the truth of (5.26) was proved by Simon using the coherent state embedding $x \mapsto \boldsymbol{x}$,

as follows:

$$\tfrac{1}{\text{vol}(X)} \int_X e^{-\langle x, Z\rangle}\, dx = \tfrac{1}{\text{vol}(X)} \int_X e^{-\text{Tr}(\boldsymbol{x}\boldsymbol{H}_Z)}\, dx \tag{5.28a}$$

$$\leq \tfrac{1}{\text{vol}(X)} \int_X \text{Tr}(\boldsymbol{x} e^{-\boldsymbol{H}_Z})\, dx \tag{5.28b}$$

$$= \tfrac{1}{\dim(\mathcal{H})} \text{Tr}(e^{-\boldsymbol{H}_Z}). \tag{5.28c}$$

Here (a) is just (5.11); (b) follows from the convexity of the exponential map:

$$e^{\text{Tr}(\boldsymbol{x}\boldsymbol{H}_Z)} = \exp(\sum c_i h_i) \leq \sum c_i \exp(h_i) = \text{Tr}(\boldsymbol{x} e^{\boldsymbol{H}_Z}) \tag{5.29}$$

where $\boldsymbol{H}_Z = \sum h_i \boldsymbol{E}_i$ is the spectral decomposition of $\boldsymbol{H}_Z$ and $c_i = \text{Tr}(\boldsymbol{x}\boldsymbol{E}_i)$; and (c) is by Schur's lemma: $\int_X \boldsymbol{x}\, dx = \tfrac{\text{vol}(X)}{\dim(\mathcal{H})} \mathbf{1}$. To discuss the case $\lambda < \mu$, on the other hand, we note that for $Z \in \mathfrak{t}$ one has

$$\text{Tr}(e^{-\boldsymbol{H}_Z}) = \sum_{\nu:\,\text{weight}} \dim(\mathcal{H}^\nu) e^{-\langle \nu, Z\rangle} \tag{5.30}$$

while Harish-Chandra's formula [B92] gives, with $dx$ normalized as $\tfrac{1}{n!}\left(\tfrac{\sigma}{2\pi}\right)^n$ ($n = \tfrac{1}{2}\dim(X)$),

$$\int_X e^{-\langle x, Z\rangle}\, dx = \sum_{p \in W(\mu)} \frac{e^{-\langle p, Z\rangle}}{\prod_{\alpha \in \text{R}_p} \langle \alpha, Z\rangle} \tag{5.31}$$

where $\text{R}_p$ is the set of weights of the isotropy representation of $T$ in $T_p X$.[12] Now if $\lambda < \mu$, there is an $A \in \text{C}$ such that $\langle \lambda, A\rangle < \langle \mu, A\rangle$; putting $Z = -At$, we see that (5.30) grows like $e^{\langle \lambda, A\rangle t}$ and (5.31) like $t^{-n} e^{\langle \mu, A\rangle t}$, so that (5.26) fails for sufficiently large $t$.   **q.e.d.**

One may ask if (5.26), thought of as an inequality between the classical and quantum partition functions at cotemperature $Z \in \mathfrak{g}$ (cf. [S66, p. 326]), is true in any generality. We wish only to observe here that it requires the

---

[12]Explicitly $\text{R}_p = \{\alpha \in \text{R} : \langle p, \alpha^\vee\rangle < 0\}$ where $\alpha^\vee$ denotes the unique element of $[\mathfrak{g}^\alpha, \mathfrak{g}^{-\alpha}]$ such that $\langle \alpha, \alpha^\vee\rangle = 2$. (Cf. 5.14b.)

*non*-inclusion of 'half-forms', 'rho-shift', 'zero-point energy', etc. Thus for a harmonic oscillator the analogue of (5.26) would be

$$\sum_{n=0}^{\infty} e^{-\beta n} \geq \tfrac{1}{2\pi} \int_{\mathbf{R}^2} e^{-\beta(p^2+q^2)/2} dp\, dq. \qquad (5.32)$$

This inequality is indeed true for all $\beta$; but it would be reversed if $n$ were replaced by $n + \tfrac{1}{2}$ inside the sum. (Compare [S80], ineqs. (3.8) and (6.1).)

> *An eigenstate belonging to an eigenvalue in a range is a mathematical idealization of what can be attained in practice. All the same such eigenstates play a very useful role in the theory and one could not very well do without them.* —P. A. M. Dirac

# 6. Localized states

We have seen after (4.10) that quantum representations for a noncompact orbit make a rather untractable class. On the other hand, we described in (4.15–16) certain quantum representations in which states could be 'localized' on coisotropic submanifolds of the orbit, thus solving in some cases, in a rather unexpected way, what A. Weinstein [W82] once called the 'fundamental quantization problem': to attach quantum states to lagrangian submanifolds.

It is natural to ask if such solutions are unique. We will show that this is so in a number of cases, and that in fact the property of admitting an eigenvector under a sufficiently large subgroup often characterizes a representation entirely. This is somewhat reminiscent of the representation theory of Lie algebras, where the classification of modules only becomes manageable after one imposes the presence of eigenvectors [B90].

$X$ continues to denote a coadjoint orbit of the Lie group $G$.

## 6.1 Definition

A reasonable definition in general would seem to be the following. Consider a closed subgroup $H \subset G$, and a coadjoint orbit $Y$ of $H$ contained in $X_{|\mathfrak{h}}$.

**6.1 Definition.** A quantum state for $X$ is *localized* at $Y \subset \mathfrak{h}^*$ if the restriction $m_{|H}$ is a quantum state for $Y$.

We may also say by abuse that such a state is localized on $\pi^{-1}(Y)$, where $\pi$ denotes the projection $X \to \mathfrak{h}^*$. We recall from [K78] that this set is generically a coisotropic submanifold of $X$. Note also from (2.23) that if $g \in G$ normalizes $H$, then its natural actions on $\mathfrak{h}^*$ and on quantum states (cf. 3.3) are such that

$$m \text{ localized at } Y \quad \Leftrightarrow \quad m_g \text{ localized at } g(Y). \tag{6.2}$$

We will mostly apply (6.1) when $H$ is connected and $Y$ is a point-orbit; then $m_{|H}$ will be a character $\chi$ (by 3.6), and the cyclic vector $e_e = \overline{m}$ of the Gel'fand-Naĭmark-Segal module $\mathcal{H}_m$ will be an eigenvector of type $\chi$ under $H$ (by 1.19).

## 6.2 Nilpotent groups

Here $G$ is assumed nilpotent, connected, simply connected, $H$ connected; we fix $x \in X$ and $\xi$ as in (4.1), and recall that $\{x_{|\mathfrak{h}}\}$ will be a point-orbit just when $H$ is subordinate to $x$.

**6.3 Theorem.** *Let $H$ be maximal subordinate to $x$. Then there is a unique quantum state for $X$ localized at $\{x_{|\mathfrak{h}}\}$, namely*

$$m(g) = \begin{cases} \xi(g) & \text{if } g \in H, \\ 0 & \text{otherwise.} \end{cases} \tag{6.4}$$

**Proof.** The fact that $m$ must coincide with $\xi$ in $H$ is just (3.6). To see that it must vanish outside $H$, we consider the sequence $H = H_0 \subset H_1 \subset H_2 \ldots$ where $H_{i+1}$ is the normalizer of $H_i$ in $G$. Since $G$ is nilpotent, the $H_i$ are connected and all equal to $G$ after finitely many steps [B60, ch. III, p. 236]; so it is enough to show inductively that $m$ vanishes in $H_{i+1} - H_i$ for all $i$.

*Case $i = 0$.* Let $a \in H_1 - H$. Applying (1.19) twice, we get

$$\xi(h)m(a) = m(ha) = m(aa^{-1}ha) = m(a)\xi(a^{-1}ha) \tag{6.5}$$

for all $h \in H$. Thus, if $m(a)$ was nonzero, $a$ would both normalize $H$ and stabilize $\xi_{|H}$, and we would conclude just as in the proof of (4.15b) that $H$ is not maximal subordinate to $x$.

*Case $i > 0$.* Let $a \in H_{i+1} - H_i$. Then $a$ normalizes $H_i$ but not $H$, so we can fix an $h \in H$ such that $a^{-1}ha \in H_i - H$. Putting $g_n = h^n a$ it follows that $g_p^{-1} g_q \in H_i - H$ whenever $p \neq q$. The induction hypothesis then shows that $m(g_p^{-1} g_q) = 0$, which is to say that the $\delta^{g_n}$ make an orthonormal set relative to the sesquilinear form (1.16). Therefore Bessel's inequality gives

$$\sum_n |m(g_n)|^2 = \sum_n |(\delta^e, \delta^{g_n})_m|^2 \leq (\delta^e, \delta^e)_m = 1. \tag{6.6}$$

Now this forces $m(a) = 0$, because we have $|m(g_n)| = |\xi(h^n) m(a)| = |m(a)|$ for all $n$. **q.e.d.**

**6.7 Remarks.** (a) Note that the condition that $m$ be quantum for $X$ has not been used in the proof. Rather it comes as a bonus from (4.8).

(b) Referring back to Example 4.28, we now see that (4.33) is the *unique* state localized at $q = 0$. In particular the uniform spread of $(p, h)$ over the compactified energy-momentum plane, irrespective of the relation $h = \frac{1}{2}p^2$, is a necessary consequence of this localization. Quite to the point, a main ingredient in the proof was (1.19), which is really a version of the Heisenberg inequalities: compare [F88, p. 529].

(c) One can read (6.3) as a Frobenius reciprocity theorem. In fact, given an irreducible unitary $G$-module $\mathcal{H}$, write $\mathcal{H}^{H,\xi}$ for the subspace of $H$-eigenvectors of type $\xi_{|H}$. One has then a linear map

$$\mathrm{Hom}_G(\mathcal{H}, \mathrm{ind}_H^G \xi_{|H}) \longrightarrow \mathcal{H}^{H,\xi} \tag{6.8}$$

sending $J$ to the vector $\vartheta = J^* e_e$, in the notation of (4.16). Frobenius' Theorem asserts that similar maps are isomorphisms for finite groups, but in general this need not a priori be so [M51]. Nevertheless, here

**6.9 Corollary.** *(6.8) is an isomorphism.*

**Proof.** Injectiveness is because one recovers $J$ from $\vartheta$ by $(J\varphi)(g) = (g\vartheta, \varphi)$ as in (1.23); this is perfectly general. To prove surjectiveness, suppose $\psi$ is a unit vector in $\mathcal{H}^{H,\xi}$. By (6.3) and (6.7a), the state $m(g) = (\psi, g\psi)$ must be given by (6.4). Since $\mathcal{H}$ is irreducible, $\psi$ is cyclic and so [F88,VI.19.8] asserts that $\mathcal{H}$ coincides with the Gel'fand-Naĭmark-Segal module $\mathcal{H}_m$, which is just $\mathrm{ind}_H^G \xi_{|H}$ (1.25).     **q.e.d.**

**Note added.** D. Testard has pointed out that Theorem 6.3 was known in the case of the Heisenberg group: [B74, thm 2.18]. This paper predates [E81] (cf. 4.21c), and also has some overlap with Example 4.22.

## 6.3 Compact groups

Definition 6.1 turns out to be productive also in singular cases where the preimage of $Y$ in $X$ is not coisotropic. For instance if $G$ is compact connected and $H = T$ a maximal torus, the useful case to consider is when $Y$ is an extreme point of the polytope $X_{|\mathfrak{t}}$, say the dominant element $\mu$ of $X$ (5.6b). We have then:

**6.10 Theorem.** *If $\mu$ is integral, there is a unique quantum state for $X(\mu)$ localized at $\{\mu\} \subset \mathfrak{t}^*$, namely $m(g) = (\varphi, g\varphi)$ where $\varphi$ is a maximal weight vector in $\mathcal{H}(\mu)$. Otherwise there is no such state.*

**Proof.** Let $\mathcal{H}_m = \bigoplus_j \mathcal{H}(\lambda_j)$ be a decomposition of the Gel'fand-Naĭmark-Segal module generated by $m$ into irreducibles. As $\mathcal{H}_m$ is quantum for $X$ (3.4), all $\lambda_j$ are $\leq \mu$ (5.16). Moreover, as we saw after (6.2), the cyclic vector $e_e$ is a weight vector of weight $\mu$. So $\mu$ must be integral, and $e_e$ is orthogonal to all summands of maximal weight $\lambda_j < \mu$, which must therefore vanish since $e_e$ is cyclic. Also it is orthogonal to all except the maximal weight space in each remaining summand; so its decomposition writes $e_e = \sum_j c_j \varphi_j$ where $\varphi_j$ is a normalized maximal weight vector in $\mathcal{H}(\lambda_j) \simeq \mathcal{H}(\mu)$. Now we have $(\varphi_j, g\varphi_k) = \delta_{jk}(\varphi, g\varphi)$ where $\varphi$ is as in the statement of the Theorem, and hence

$$m(g) = (e_e, ge_e) = \sum_{j,k} \bar{c}_j c_k (\varphi_j, g\varphi_k) = (\varphi, g\varphi). \tag{6.11}$$

(Of course, it follows a posteriori that there was only one summand.)     **q.e.d.**

**6.12 Remarks.** (a) Thus we get a third characterization, after (5.23) and (5.27), of $\mathcal{H}(\mu)$ among irreducible quantum modules for $X(\mu)$: namely, it is the only one to accomodate a state localized at $\{\mu\}$.

(b) Letting the Weyl group act, (6.2) will give a unique quantum state localized at any other extreme point of $X_{|\mathfrak{t}}$.

**6.13 Examples.** Referring back to Example 5.19, the state defined by the maximal weight vector $\varphi_k = e_1 \wedge \cdots \wedge e_k$ in $\mathcal{H}(\pi_k)$ is the minor

$$m_k(g) = (\varphi_k, g\varphi_k) = \begin{vmatrix} g_{11} & \cdots & g_{1k} \\ \cdots\cdots\cdots\cdots \\ g_{k1} & \cdots & g_{kk} \end{vmatrix}. \tag{6.14}$$

(Note that this illustrates (1.26b) with $H = G_{\pi_k} = \mathbf{U}(k) \times \mathbf{U}(n-k)$.) Likewise the state $m_1^{\ell_1} \cdots m_n^{\ell_n}$ belongs to the module with maximal weight $\ell_1\pi_1 + \cdots + \ell_n\pi_n$; cf. [C13]. Specializing to $\mathbf{U}(3)$, (6.10) and (6.12b) give for instance the following table (cf. Fig. 2):

| name | weight $\mu$ | state $m(g)$ |
| --- | --- | --- |
| up quark | $u$ | $g_{11}$ |
| strange antiquark | $-s$ | $(g_{11}g_{22} - g_{12}g_{21})/\det(g)$ |
| proton | $2u + d$ | $g_{11}(g_{11}g_{22} - g_{12}g_{21})$ |
| xi$^-$ hyperon | $d + 2s$ | $g_{33}(g_{33}g_{22} - g_{32}g_{23})$ |
| pi$^+$ meson | $u - d$ | $g_{11}(g_{11}g_{33} - g_{13}g_{31})/\det(g)$ |
| omega$^-$ | $3s$ | $g_{33}^3.$ |

## Appendix: Proof of (3.14)

Proposition 3.14 is a corollary of [Z96], but does not require the full generality of that paper. We give here a direct proof along the lines of [G83], where the case of certain semidirect products was already treated; cf. also the related work in [R75] and [D92]. As we have learned to expect, the main point will reside in the following:

**A.1 Lemma.** *In the setting of (3.8) one automatically has $x+\mathrm{orth}(\mathfrak{h}) \subset X$ ('Pukánszky condition'). In more detail one has*

 (i) $\mathfrak{a}(x) = \mathrm{orth}(\mathfrak{h})$;
 (ii) $A(x) = x + \mathrm{orth}(\mathfrak{h})$;
 (iii) $H(x) = \eta^{-1}(H(x_{|\mathfrak{h}}))$, *where $\eta$ is the projection $\mathfrak{g}^* \to \mathfrak{h}^*$.*

**Proof.** (i) We have $\langle \mathfrak{a}(x), Z \rangle = \langle x, [\mathfrak{a}, Z] \rangle = \langle p, [\mathfrak{a}, Z] \rangle = \langle Z(p), \mathfrak{a} \rangle$ for all $Z \in \mathfrak{g}$. Consequently, $\mathrm{orth}(\mathfrak{a}(x)) = \mathfrak{g}_p = \mathfrak{h}$.

(ii) The relation $\langle \mathfrak{a}(x), \mathfrak{h} \rangle = 0$ just established shows that $\mathfrak{a}$ stabilizes $x_{|\mathfrak{h}}$; therefore so does $A$, which is to say that $A(x) \subset x + \mathrm{orth}(\mathfrak{h})$. For the reverse inclusion we note that since $\mathfrak{a}$ is an ideal we have $\mathrm{ad}(Z)^n(\mathfrak{g}) \subset [\mathfrak{a},\mathfrak{a}]$ for all $Z \in \mathfrak{a}$ and all $n \geq 2$. As $\mathfrak{a}$ is $X$-abelian this implies that

$$\langle \exp(Z)(x), Z' \rangle = \Big\langle x, \sum_{n=0}^{\infty} \frac{(-)^n}{n!} \mathrm{ad}(Z)^n(Z') \Big\rangle$$
$$= \langle x, Z' - [Z, Z'] \rangle$$
$$= \langle x + Z(x), Z' \rangle \qquad (A.2)$$

for all $Z \in \mathfrak{a}$, $Z' \in \mathfrak{g}$. Thus we see that $A(x)$ contains $x + \mathfrak{a}(x)$, as desired. Finally (iii) follows by writing (ii) in the form $A(x) = \eta^{-1}(x_{|\mathfrak{h}})$ and letting $H$ act on both sides. **q.e.d.**

**Proof of (3.14).** We want to show that $H(x_{|\mathfrak{h}})$ is the unique coadjoint orbit $Y$ of $H$ such that, using the notation of §1.1D,

 (a) $\Phi_{\mathrm{ind}}$ is a symplectic diffeomorphism $\mathrm{Ind}_H^G Y \to X$;
 (b) $Y_{|\mathfrak{a}} = \{p\}$.

Assume that $Y = H(x_{|\mathfrak{h}})$. Then (b) is clear, and A.1(iii) says that $X$ is the only coadjoint orbit of $G$ whose projection meets $Y$; hence $\Phi_{\mathrm{ind}}$ is onto $X$, by (1.13). To see that it is one-to-one, it is enough by equivariance to show that $\Phi_{\mathrm{ind}}^{-1}(x)$ is one point, or in other words, that $\phi^{-1}(x) \cap \psi^{-1}(0)$ is one $H$-orbit. Now the formulas for $\phi$ and $\psi$ in §1.1D exhibit the latter set as

$$\{(xq, q^{-1}(x)_{|\mathfrak{h}}) : q \in Q\}, \qquad (A.3)$$

where $Q$ consists of all $q \in G$ such that $q^{-1}(x) \in \eta^{-1}(Y)$. But then A.1(iii) shows that $Q = H$, so that (A.3) is indeed a single $H$-orbit. Thus $\Phi_{\text{ind}}$ is a bijection onto $X$ and hence, by (1.7), a symplectic diffeomorphism.

Conversely assume that $Y$ satisfies (a, b). Then $X$ lies in the image of $\Phi_{\text{ind}}$, and therefore $Y$ lies in $X_{|\mathfrak{h}}$ (1.13). So given $y \in Y$ there will be some $g \in G$ such that $y = g(x)_{|\mathfrak{h}}$. Projecting this relation in $\mathfrak{a}^*$ and using (b) we deduce that $p = g(p)$. Therefore $g \in H$, and $Y = H(x_{|\mathfrak{h}})$.

To prove the remaining statement, we note that $\pi^{-1}(p) = H(x)$ since $\pi$ is $G$-equivariant, and that an $A$-orbit in this set is just a fiber of $X \to \mathfrak{h}^*$: A.1(ii) says this for one fiber, but it then follows for all because $A$ is normal. Thus we obtain $\pi^{-1}(p)/A = H(x)_{|\mathfrak{h}} = Y$, as claimed. **q.e.d.**

> ... to fight the bourgeois concept of representation? But I wonder if this is not setting the problem in academic terms.   —J.-L. Godard

# References

(The labelling is one-to-one.)

# THÉORIE DE MACKEY SYMPLECTIQUE
## avec une application
## aux
# REPRÉSENTATIONS QUANTIQUES

François Ziegler
CNRS-CPT Luminy, Case 907
13288 Marseille Cedex 09
ziegler@cpt.univ-mrs.fr

*Table*

# 0. Introduction

KAZHDAN, KOSTANT ET STERNBERG ont introduit dans [K78] une construction destinée à jouer, en géométrie symplectique, le rôle que tient l'induction unitaire en théorie des représentations. Comme on sait, plusieurs méthodes dont le prototype est celle de Kirillov, associent des variétés symplectiques à des représentations unitaires; et l'idée que la correspondance puisse entrelacer induction symplectique et induction unitaire est manifestement une des motivations de [K78].

Ceci suggère aussitôt l'existence d'une version symplectique de la théorie de Mackey ('normal subgroup analysis')—ce que nous développerons ici. L'entreprise n'est certes pas très originale: sa possibilité transparaît clairement dans la plupart des articles cités en bibliographie, et de fait, la page 1 de [K66] mentionne déjà, parmi huit projets, celui de traiter symplectiquement 'la théorie de Mackey pour le cas des produits semidirects'; projet dont il revient à Rawnsley [R75] d'avoir donné une première réalisation. Quoique négligée jusqu'ici, l'extension à des cas plus généraux peut se justifier au moins par la vision *synoptique* qu'on y gagne, chaque fois qu'on applique la théorie de Mackey aux représentations associées à des variétés symplectiques.

Mais notre motivation véritable réside ailleurs. C'est que la *définition même* de ces représentations semble souvent dictée, a posteriori, par l'idée d'entrelacement mentionnée plus haut. Il est alors naturel de se demander si on ne pourrait faire résulter cela de quelque principe général—comme celui de [S88, Z96], qui exige en substance que 'le spectre quantique d'observables qui commutent soit concentré sur l'ensemble de valeurs classique'. En effet, nous avons vu dans [Z96, 3.9] des cas où ce principe implique qu'à une variété induite correspondent nécessairement des représentations induites; et nous verrons au §6 comment cela peut se généraliser à toute variété induite.

De ce point de vue, notre principale observation sera que l'*imprimitivité* se manifeste, en géométrie symplectique, par la présence d'un certain type d'algèbre de fonctions en involution, auxquelles on pourra précisément imposer le principe ci-dessus; et les §§ 1 à 5 sont surtout destinés à montrer l'étendue du domaine où ceci peut s'appliquer.

$$* \quad * \quad * \quad *$$

**Notation.** Toutes les variétés (et groupes de Lie) de cet article sont supposées séparées et dénombrables à l'infini. Si $G$ est un groupe de Lie, alors $\mathfrak{g}$ est son algèbre de Lie. Si $G$ agit sur une variété $X$, de sorte qu'on a un morphisme $g \mapsto g_X$ de $G$ dans les difféomorphismes de $X$, avec $g_X(x)$ fonction $C^\infty$ de $(g,x)$, on définit l'action infinitésimale

$$Z \mapsto Z_X, \qquad \mathfrak{g} \to \text{champs de vecteurs sur } X \qquad (0.1)$$

par $Z_X(x) = \frac{d}{dt}\exp(tZ)_X(x)\big|_{t=0}$. C'est un morphisme d'algèbres de Lie, si on définit le crochet des champs de vecteurs avec l'opposé du signe habituel. Chaque fois que ce sera possible, nous omettrons les indices et écrirons $g(x)$ et $Z(x)$, au lieu de $g_X(x)$ et $Z_X(x)$. Les notations

$$G(x), \qquad \mathfrak{g}(x), \qquad G_x, \qquad \mathfrak{g}_x, \qquad (0.2)$$

désigneront respectivement la $G$-orbite de $x$, son espace tangent en $x$, le stabilisateur de $x$ dans $G$, et le stabilisateur de $x$ dans $\mathfrak{g}$.

Il sera commode d'avoir une notation concise pour l'action des translations sur les vecteurs tangents et cotangents au groupe; ainsi, si $g, q \in G$ nous écrirons

$$\begin{array}{cc} T_q G \to T_{gq} G \\ v \mapsto gv, \end{array} \qquad \text{resp.} \qquad \begin{array}{cc} T_q^* G \to T_{gq}^* G \\ p \mapsto gp \end{array} \qquad (0.3)$$

pour la dérivée de $q \mapsto gq$, resp. pour l'application transposée définie par $\langle gp, v \rangle = \langle p, g^{-1}v \rangle$. On définit de même $vg$ et $pg$ avec $\langle pg, v \rangle = \langle p, vg^{-1} \rangle$. Rappelons enfin l'action coadjointe, $g(x) = gxg^{-1}$, sur $\mathfrak{g}^* = T_e^* G$; infinitésimalement cela donne: $Z(x) = \langle x, [\,\cdot\,, Z] \rangle$.

# 1. Rappels symplectiques

**A. Champs de vecteurs hamiltoniens.** Soit $X$ une variété symplectique —une variété munie d'une 2-forme non dégénérée et fermée $\sigma$. Un champ de vecteurs $\eta$ sur $X$ est dit *symplectique* si son flot préserve la 2-forme: $\pounds(\eta)\sigma = 0$. D'après la formule de Cartan pour la dérivée de Lie, $\pounds(\eta)\sigma = i(\eta)d\sigma + di(\eta)\sigma$, ceci a lieu ssi la 1-forme $i(\eta)\sigma = \sigma(\eta,\cdot)$ est fermée. Si de plus cette 1-forme est exacte, de sorte qu'on a

$$i(\eta)\sigma = -dH, \tag{1.1}$$

pour une fonction $H$ sur $X$, nous dirons que $\eta$ est *hamiltonien*, $\eta \in \mathrm{Ham}(X)$, et nous écrirons $\eta = \mathrm{drag}\, H$ ('gradient symplectique'). Le champ $\eta$ détermine $H$ à une constante additive locale près, de sorte que si $X$ a $c$ composantes connexes on a une suite exacte

$$0 \longrightarrow \mathbf{R}^c \longrightarrow \mathrm{C}^\infty(X) \overset{\mathrm{drag}}{\longrightarrow} \mathrm{Ham}(X) \longrightarrow 0. \tag{1.2}$$

C'est une extension centrale d'algèbres de Lie si l'on munit $\mathrm{C}^\infty(X)$ du *crochet de Poisson*: $\{H, H'\} = \sigma(\mathrm{drag}\, H', \mathrm{drag}\, H)$.

**B. $G$-espaces hamiltoniens.** L'action sur $X$ d'un groupe de Lie $G$, préservant $\sigma$, est dite hamiltonienne si (0.1) est à valeurs dans $\mathrm{Ham}(X)$; ou en d'autres termes, s'il existe une *application moment*, $\varPhi : X \to \mathfrak{g}^*$, telle que

$$i(Z_X)\sigma = -dH_Z \qquad \text{où} \qquad H_Z = \langle \varPhi(\cdot), Z \rangle. \tag{1.3}$$

On dit que l'application moment est *équivariante*, et que le couple $(X, \varPhi)$ est un *$G$-espace hamiltonien*, si $\varPhi$ entrelace l'action de $G$ sur $X$ et l'action coadjointe sur $\mathfrak{g}^*$. Infinitésimalement cela signifie que $Z \mapsto H_Z$ est un morphisme:

$$\{H_Z, H_{Z'}\} = H_{[Z,Z']}. \tag{1.4}$$

(Comme on sait, cette condition n'est pas vraiment restrictive; on peut toujours la satisfaire en passant à une extension centrale de $\mathfrak{g}$ par $\mathbf{R}^c$ (1.2).)

La notion d'*isomorphisme* de $G$-espaces hamiltoniens est claire:

$$\operatorname{Hom}_G(X_1, X_2) \tag{1.5}$$

désignera l'ensemble des difféomorphismes $X_1 \to X_2$ qui transforment $\sigma_1$ en $\sigma_2$ et $\Phi_1$ en $\Phi_2$.

**C. Exemples de base.** (a) Si $G$ agit sur une variété $Q$, on obtient une action sur $X = T^*Q$ qui préserve la 1-forme canonique $\varpi = \text{``}\langle p, dq \rangle\text{''}$ et $\sigma = d\varpi$. Dans ce cas, $i(Z_X)d\varpi + di(Z_X)\varpi = 0$ montre qu'une application moment (équivariante) est donnée par la formule d'E. Nœther: $H_Z = i(Z_X)\varpi$, i.e.,

$$\langle \Phi(p), Z \rangle = \langle p, Z(q) \rangle, \qquad p \in T_q^*Q. \tag{1.6}$$

(b) Si $X$ est une orbite dans $\mathfrak{g}^*$ pous l'action coadjointe de $G$, alors la 2-forme définie sur $X$ par $\sigma(Z(x), Z'(x)) = \langle Z(x), Z' \rangle$ fait de $(X, X \hookrightarrow \mathfrak{g}^*)$ un $G$-espace hamiltonien homogène. Réciproquement (Kirillov-Kostant-Souriau), un tel espace est toujours un revêtement d'orbite coadjointe:

**1.7 Théorème.** *Soit $(X, \Phi)$ un $G$-espace hamiltonien, et supposons que $G$ agisse transitivement sur $X$. Alors $\Phi$ est un revêtement symplectique de son image, qui est une orbite coadjointe de $G$.*

Par revêtement *symplectique* on entend bien sûr que $\Phi$ lie la 2-forme de l'orbite à celle de $X$; c'est une simple conséquence de (1.4). Que $\Phi$ soit un revêtement (ait des fibres discrètes) résulte de la première de deux conséquences instructives de (1.3), valables pour toute application moment:

$$\operatorname{Ker}(D\Phi(x)) = \mathfrak{g}(x)^\sigma, \qquad \operatorname{Im}(D\Phi(x)) = \operatorname{orth}(\mathfrak{g}_x). \tag{1.8}$$

Dans ces formules, l'exposant $^\sigma$ signifie 'sous-espace orthogonal relativement à $\sigma$'; et si $(\cdot)$ est une partie de $\mathfrak{g}$ ou de $\mathfrak{g}^*$, $\operatorname{orth}(\cdot)$ désigne son annulateur dans l'autre.

**D. Induction symplectique.** Nous pouvons maintenant rappeler la construction de [K78]. Etant donné un sous-groupe fermé $H$ de $G$, et un $H$-espace hamiltonien $(Y, \Psi)$, elle produit un $G$-espace hamiltonien $(\operatorname{Ind}_H^G Y, \Phi_{\text{ind}})$ comme suit. *(Nous utilisons la notation (0.3).)*

Munissons $N = T^*G \times Y$ de la forme symplectique $\omega = d\varpi + \tau$, où $\varpi$ est la 1-forme canonique de $T^*G$ et $\tau$ la 2-forme donnée sur $Y$; et faisons agir $H$ 'diagonalement' sur $N$ par $h(p, y) = (ph^{-1}, h(y))$. C'est une action hamiltonienne, d'application moment $\psi$:

$$\psi(p, y) = \Psi(y) - q^{-1}p_{|\mathfrak{h}} \tag{1.9}$$

pour $p \in T_q^*G$. Le second terme désigne ici la restriction à $\mathfrak{h}$ de $q^{-1}p \in \mathfrak{g}^*$; il vient de (1.6) avec $Z(q) = -qZ$. La variété induite est alors définie comme la réduite de Marsden-Weinstein de $N$ en zéro [M74], i.e.,

$$\operatorname{Ind}_H^G Y := \psi^{-1}(0)/H. \tag{1.10}$$

C'est une variété de dimension $2\dim(G/H) + \dim(Y)$, qui admet une 2-forme naturelle $\sigma_{\text{ind}}$ définie par la condition que son image réciproque dans $\psi^{-1}(0)$ soit la restriction de $\omega$. Pour en faire un $G$-espace, faisons agir $G$ sur $N$ par $g(p, y) = (gp, y)$. Cette action commute avec l'action diagonale de $H$, et préserve $\psi^{-1}(0)$. De plus elle est hamiltonienne et son application moment $\phi : N \to \mathfrak{g}^*$, donnée par (1.6) avec maintenant $Z(q) = Zq$:

$$\phi(p, y) = pq^{-1}, \qquad p \in T_q^*G, \tag{1.11}$$

est constante sur les $H$-orbites. Passant au quotient, on obtient l'action induite de $G$ sur $\operatorname{Ind}_H^G Y$, et son application moment $\Phi_{\text{ind}} : \operatorname{Ind}_H^G Y \to \mathfrak{g}^*$.

## 2. Le théorème d'imprimitivité

**A. Systèmes d'imprimitivité.** Soit $X$ un $G$-espace hamiltonien, dont on notera $\sigma$ la 2-forme et $\Phi : X \to \mathfrak{g}^*$ l'application moment. L'action naturelle de $G$ sur l'algèbre de Lie $\mathrm{C}^\infty(X)$ sera notée: $g(f) = f(g^{-1}(\cdot))$; elle préserve le crochet de Poisson.

**2.1 Définition.** Un *système d'imprimitivité* sur $X$ est une sous-algèbre abélienne $G$-invariante $\mathfrak{f}$ de $\mathrm{C}^\infty(X)$, telle que le champ hamiltonien drag $f$ soit complet pour toute $f \in \mathfrak{f}$.

La terminologie peut se justifier ainsi: selon des principes consacrés, toute 'quantification' de $X$ doit fournir une représentation $U$ de $G$ (resp. $E$ de $\mathfrak{f}$) par des opérateurs unitaires (resp. autoadjoints), avec naturellement la relation d'équivariance $E(g(f)) = U(g)E(f)U(g^{-1})$; d'où un système d'imprimitivité au sens classique de Mackey.

Si $\mathfrak{f}$ est un système d'imprimitivité, on notera $\mathfrak{f}^*$ son dual algébrique, et $\mathcal{F}$ l'ensemble $\mathfrak{f}$ vu comme groupe additif. Alors $G$ agit sur $\mathfrak{f}^*$ par contragrédience, $\langle g(b), f \rangle = \langle b, g^{-1}(f) \rangle$, et $\mathcal{F}$ agit sur $X$ par exponentiation des champs drag $f$. Bien que $\mathcal{F}$ ne soit pas en général un groupe de Lie, on pourra encore noter $\mathcal{F}(x)$ la $\mathcal{F}$-orbite de $x$, $\mathfrak{f}(x) = \{(\mathrm{drag}\, f)(x) : f \in \mathfrak{f}\}$, et considérer comme *moment* de cette action l'application

$$\pi : X \to \mathfrak{f}^*, \qquad \langle \pi(x), f \rangle = f(x). \tag{2.2}$$

L'ensemble $B = \pi(X)$ s'appellera *base* de $\mathfrak{f}$. Comme (2.2) est visiblement $G$-équivariante, $B$ est en général une réunion de $G$-orbites. D'où:

**2.3 Définition.** Le système d'imprimitivité $\mathfrak{f}$ sera dit *transitif* si
  (i) l'action de $G$ sur la base $B = \pi(X)$ est transitive;
  (ii) $\pi : X \to B$ est $C^\infty$ pour la structure de variété $G$-homogène de $B$.

**2.4 Remarques.** (a) C'est par exemple automatiquement le cas lorsque l'action de $G$ sur $X$ est elle-même transitive.

(b) La structure de variété dont il est question est bien définie, car le stabilisateur de n'importe quel point $b = \pi(x)$ de $B$ est fermé. En effet, $g \in G_b$ signifie qu'on a $\langle g(b), f \rangle = \langle b, f \rangle$, i.e. $f(g(x)) = f(x)$, pour toute $f \in \mathfrak{f}$; or cette condition est fermée par continuité des applications $g \mapsto f(g(x))$.

(c) Nous appellerons encore *base de* $\mathfrak{f}$ tout $G$-ensemble muni d'une bijection $G$-équivariante *fixée* avec $B$. Ceci permet de parler de systèmes d'imprimitivité *de même base*.[1]

---

[1] Une différence formelle entre (2.1) et Mackey est que nous définissons la base en termes de $\mathfrak{f}$, tandis que Mackey se la donne d'avance. Mais j'imagine, sans trop y avoir réfléchi, que pour un $G$-module unitaire $\mathcal{H}$ on pourrait définir d'abord un système d'imprimitivité comme une sous-algèbre commutative $G$-invariante de $\mathrm{End}(\mathcal{H})$, et ensuite sa base comme le spectre de cette sous-algèbre.

**B. Le système d'imprimitivité associé à une variété induite.** Supposons que $X = \operatorname{Ind}_H^G Y$, où $H$ est un sous-groupe fermé de $G$ et $Y$ un $H$-espace hamiltonien (§ 1D). On obtient alors une projection $G$-équivariante

$$\pi_{\text{ind}} : \operatorname{Ind}_H^G Y \to G/H \tag{2.5}$$

en remarquant que l'application $T^*G \times Y \to G/H$ qui envoie $T_q^*G \times Y$ sur $qH$ passe au quotient (1.10) car elle est constante sur les $H$-orbites. Il en résulte un système d'imprimitivité canonique sur $X$:

**2.6 Proposition.** *Si $X = \operatorname{Ind}_H^G Y$, alors $\mathfrak{f}_{\text{ind}} := \pi_{\text{ind}}^*(C^\infty(G/H))$ est un système d'imprimitivité transitif sur $X$, de base $G/H$.*

**Preuve.** (On reprend les notations du § 1D.) Soit $f \in C^\infty(G/H)$, et notons encore $f$ sa relevée à $X$, à $G$, ou à $T^*G \times Y$. Alors son gradient symplectique sur ce dernier espace est le champ $\eta$ de flot

$$e^{t\eta}(p, y) = (p - tDf(q), y) \tag{2.7}$$

($p \in T_q^*G$). En effet, dérivant (2.7) en $t = 0$ dans une carte standard $(p_i, q_i)$ du cotangent, on obtient en coordonnées $\eta = (\delta p_i, \delta q_i, \delta y) = (-\partial f/\partial q_i, 0, 0)$ et donc

$$\omega(\eta, \cdot) = \delta p_i dq_i - \delta q_i dp_i + \tau(\delta y, \cdot) = -df, \tag{2.8}$$

i.e. $\eta = \operatorname{drag} f$. Le flot (2.7) est complet; comme d'autre part $f$ est $H$-invariante, on sait [M74, cor. 3] que ce flot passe au quotient $X = \psi^{-1}(0)/H$ (1.10), où il est encore celui de $\operatorname{drag} f$ (calculé sur $X$). Ce dernier est donc aussi complet. Si par ailleurs $f'$ est une autre fonction de $G/H$, on voit sur (2.7) qu'elle est constante le long du flot, de sorte qu'on a $\{f, f'\} = 0$. Enfin l'équivariance de (2.5) montre que ces fonctions constituent un espace $G$-invariant. On a donc bien un système d'imprimitivité, dont la base $B$ s'identifie de façon évidente à $G/H$. **q.e.d.**

**C. Le théorème d'imprimitivité.** Le théorème de Mackey [F88, p. 1194] affirme que l'existence d'un système d'imprimitivité transitif (au sens classique) caractérise les représentations induites. Son analogue symplectique consistera donc à compléter (2.6) par une réciproque:

**2.9 Théorème.** *Soit $(X,\Phi)$ un $G$-espace hamiltonien admettant un système d'imprimitivité transitif $\mathfrak{f}$, et notons $H$ le stabilisateur d'un $b \in \pi(X)$. Alors il existe un unique $H$-espace hamiltonien $(Y,\Psi)$ tel que*

$$\text{(a)} \quad X = \text{Ind}_H^G Y, \qquad \text{(b)} \quad \pi^{-1}(b) = \pi_{\text{ind}}^{-1}(eH). \qquad (2.10)$$

*Explicitement $Y = \pi^{-1}(b)/\mathcal{F}$, i.e. $Y$ est l'espace réduit de $X$ en $b \in \mathfrak{f}^*$.*

**Remarques.** La condition (b) ne sert qu'à assurer l'unicité de $Y$, qui s'entend à isomorphisme près. De même par (2.10) il faut bien sûr entendre "il existe un $J \in \text{Hom}_G(X, \text{Ind}_H^G Y)$ qui envoie $\pi^{-1}(b)$ sur $\pi_{\text{ind}}^{-1}(eH)$."

La preuve détaillera la structure de $H$-espace hamiltonien de $\pi^{-1}(b)/\mathcal{F}$.

**Preuve.** 1. Le niveau $X_b := \pi^{-1}(b)$ est une sous-variété de $X$. En effet, comme $\pi : X \to B$ est équivariante, sa dérivée en $x \in X_b$ (2.3ii) applique $\mathfrak{g}(x)$ sur $\mathfrak{g}(b) = T_b B$. Donc $\pi$ est une submersion, d'où notre assertion.

2. Cette sous-variété est *coisotrope*, et plus précisément, l'orthogonal symplectique de $T_x X_b$ est donné par

$$(T_x X_b)^\sigma = \mathfrak{f}(x) \qquad (2.11)$$

(qui est isotrope puisque $\mathfrak{f}$ est abélienne). En effet, la suite exacte transposée de $0 \to T_x X_b \to T_x X \to T_b B \to 0$ montre d'abord que $(T_x X_b)^\sigma$ est l'ensemble de valeurs de l'injection

$$j_x : T_b^* B \hookrightarrow T_x X \qquad (2.12)$$

obtenue en composant $T_b^* B \hookrightarrow T_x^* X$ avec l'isomorphisme $T_x^* X \to T_x X$ fourni par la structure symplectique. D'autre part, chaque $f \in \mathfrak{f}$ est par construction la relevée à $X$ d'une fonction $\dot f$ sur $B$, qui est aussi $C^\infty$ puisque $\pi$ est une submersion. Mais alors (1.1) dit que $(\text{drag } f)(x) = j_x(-D\dot f(b))$, de sorte que prouver (2.11) revient à voir que l'application

$$\mathfrak{f} \to T_b^* B, \qquad f \mapsto -D\dot f(b) \qquad (2.13)$$

est surjective. Or si $Z \in \mathfrak{g}$ et que $g_t = \exp(tZ)$, on a

$$
\begin{aligned}
Z \in \mathfrak{h} &\Leftrightarrow \langle g_t(b), f \rangle = \langle b, f \rangle && \forall f \in \mathfrak{f},\ \forall t \\
&\Leftrightarrow \dot{f}(g_t(b)) = \dot{f}(b) && \forall f \in \mathfrak{f},\ \forall t \\
&\Leftrightarrow \tfrac{d}{dt}\, \dot{f}(g_t(b)) = 0 && \forall f \in \mathfrak{f},\ \forall t \\
&\Leftrightarrow D(\dot{f} \circ g_t)(b)(Z(b)) = 0 && \forall f \in \mathfrak{f},\ \forall t \\
&\Leftrightarrow D\dot{f}(b)(Z(b)) = 0 && \forall f \in \mathfrak{f} \quad (2.14)
\end{aligned}
$$

puisque $\mathfrak{f}$ est $G$-invariante. Comme $Z \in \mathfrak{h}$ équivaut aussi à $Z(b) = 0$, ceci montre que les $D\dot{f}(b)$ séparent $T_b B$; donc (2.13) est surjective, et (2.11) est démontré.

3. L'espace d'orbites $Y = X_b/\mathcal{F}$ admet une unique structure de variété telle que $X_b \to Y$ soit une submersion. (Remarquons que (2.11) entraîne que $\mathfrak{f}(x) \subset \mathfrak{f}(x)^\sigma = T_x X_b$, de sorte que l'action de $\mathcal{F}$ préserve bien $X_b$.) En effet, il résulte de ce qu'on a dit avant (2.13) que l'action de $\mathcal{F}$ sur $X_b$ se factorise en

$$
\begin{array}{ccccc}
\mathcal{F} & \xrightarrow{(2.13)} & T_b^* B & \longrightarrow & \mathrm{Diff}(X_b) \\
f & \longmapsto & a & \longmapsto & e^{\hat{a}}.
\end{array} \qquad (2.15)
$$

où $\hat{a}$ désigne le champ de vecteurs défini sur $X_b$ par $\hat{a}(x) = j_x(a)$ (2.12). Comme (2.13) est surjective, ceci montre que les $\mathcal{F}$-orbites dans $X_b$ sont en fait les orbites d'une action du *groupe de Lie* additif $T_b^* B$. D'autre part, (1.3) et la définition de $\hat{a}$ donnent, pour tout $Z \in \mathfrak{g}$,

$$
\langle D\Phi(x)(\hat{a}(x)), Z \rangle = \sigma(\hat{a}(x), Z(x)) = \langle a, Z(b) \rangle = \langle \check{a}, Z \rangle \qquad (2.16)
$$

où $a \mapsto \check{a}$ est la transposée de $Z \mapsto Z(b)$, donc une bijection $T_b^* B \to \mathrm{orth}(\mathfrak{h})$. Ainsi $\Phi$ lie $\hat{a}$ au champ *constant* $\check{a}$, et entrelace donc l'action de $T_b^* B$ sur $X_b$ avec une simple action par translations:

$$
\Phi(e^{\hat{a}}(x)) = \Phi(x) + \check{a}. \qquad (2.17)
$$

Celle-ci étant libre et propre, celle-là l'est aussi par [B60, ch. III, §4, prop. 5], d'où notre assertion par [B90, ch. III, §1, prop. 10].

4. $Y$ est naturellement un $H$-espace hamiltonien. En effet $\sigma_{|X_b}$ s'annule précisément le long des $\mathcal{F}$-orbites (2.11), et provient donc d'une 2-forme symplectique $\tau$ de $Y$ [S69, p. 85]. De même l'action de $H$, qui préserve $X_b$, passe au quotient car elle normalise l'image de (2.15), $\mathfrak{f}$ étant $H$-invariante. Enfin (2.17) montre que l'action de $H$ sur $Y$ ainsi obtenue admet un moment $\Psi$, défini par la commutativité de

$$\begin{array}{ccc} X_b & \xrightarrow{\Phi} & \mathfrak{g}^* \\ \mathcal{F}(\cdot)\downarrow & & \downarrow \\ Y & \xrightarrow{\Psi} & \mathfrak{h}^*. \end{array} \qquad (2.18)$$

5. L'induite $\mathrm{Ind}_H^G Y$ est isomorphe à $X$ en tant que $G$-espace hamiltonien. Reprenons en effet les notations du §1D, et considérons les applications $\varepsilon$ (resp. $j$) de $G \times X_b$ dans $X$ (resp. dans $T^*G \times Y$):

$$\varepsilon(q, x) = q(x), \qquad \text{resp.} \qquad j(q, x) = (q\Phi(x), \mathcal{F}(x)). \qquad (2.19)$$

Il résulte de (1.9) et de (2.15–17) que $j$ est un difféomorphisme de $G \times X_b$ sur $\psi^{-1}(0)$. De plus, on vérifie sans peine que $j$ envoie les fibres de $\varepsilon$ sur les $H$-orbites dans $\psi^{-1}(0)$; passant au quotient, on obtient donc un difféomorphisme $J$:

$$\begin{array}{ccc} G \times X_b & \xrightarrow{j} & \psi^{-1}(0) \subset T^*G \times Y \\ \varepsilon\downarrow & & \downarrow (1.10) \\ X & \xrightarrow{J} & \mathrm{Ind}_H^G Y \end{array} \qquad (2.20)$$

qui visiblement est $G$-équivariant et envoie $X_b$ sur $\pi_{\mathrm{ind}}^{-1}(eH)$. Pour voir que $J$ est symplectique, considérons les variables qui apparaissent dans (2.19):

$$\tilde{x} = q(x), \qquad y = \mathcal{F}(x), \qquad m = \Phi(x), \qquad n = (qm, y) \qquad (2.21)$$

comme fonctions de $(q, x) \in G \times X_b$. Chaque vecteur $(\delta q, \delta x) \in T_q G \times T_x X_b$ donne alors par dérivation un vecteur $(\delta \tilde{x}, \delta y, \delta m, \delta n)$, ainsi qu'un élément $Z = q^{-1}\delta q$ de $\mathfrak{g}$. Cela étant, la définition de la 2-forme $\omega = d\varpi + \tau$ de

$T^*G \times Y$, celle de $\tau$ ci-dessus, et les formules (1.3–4), donnent

$$\begin{aligned}
\omega(\delta n, \delta' n) &= \langle \delta m, Z' \rangle - \langle \delta' m, Z \rangle + \langle m, [Z', Z] \rangle + \tau(\delta y, \delta' y) \\
&= \sigma(\delta x, Z'(x)) - \sigma(\delta' x, Z(x)) + \sigma(Z(x), Z'(x)) + \sigma(\delta x, \delta' x) \\
&= \sigma(\delta x + Z(x), \delta' x + Z'(x)) \\
&= \sigma(\delta \tilde{x}, \delta' \tilde{x}).
\end{aligned} \qquad (2.22)$$

Ceci montre que $j^*\omega = \varepsilon^*\sigma$, d'où $J^*\sigma_{\text{ind}} = \sigma$ par (2.20) et par définition de $\sigma_{\text{ind}}$. Enfin il est clair sur (1.11) et (2.19) que $J^*\Phi_{\text{ind}} = \Phi$.

6. Reste à établir la propriété d'unicité. Soit donc $Y$ une solution quelconque du problème, de sorte qu'on a un diagramme commutatif $G$-équivariant

$$\begin{array}{ccc}
X & \xrightarrow{J} & \text{Ind}_H^G Y \\
\pi \downarrow & & \downarrow \pi_{\text{ind}} \\
B & \longrightarrow & G/H,
\end{array} \qquad (2.23)$$

où les flèches horizontales sont des isomorphismes, celle du bas envoyant $b$ sur $eH$. Alors $T_b^*B$ se trouve identifié à $(\mathfrak{g}/\mathfrak{h})^* \simeq \text{orth}(\mathfrak{h})$, et son action (2.15) sur $X_b$ à une action de $\text{orth}(\mathfrak{h})$ sur $\pi_{\text{ind}}^{-1}(eH)$. Or les éléments de cet ensemble, vus comme $H$-orbites dans $T^*G \times Y$ (1.10), ont chacun un unique représentant dans $T_e^*G \times Y$, de sorte qu'on peut le voir comme la partie de $\mathfrak{g}^* \times Y$ constituée des $(m, y)$ tels que $m_{|\mathfrak{h}} = \Psi(y)$ (1.9). Alors l'action de $\text{orth}(\mathfrak{h})$ se fait par translation de $m$ et la 2-forme se réduit à celle de $Y$, qui s'identifie donc bien au quotient $X_b/\mathcal{F}$.   **q.e.d.**

**2.24 Exemples.** (a) Si $X$ est une orbite coadjointe de $G$ et $A$ un sous-groupe abélien invariant fermé, les hamiltoniens $H_Z(x) = \langle x, Z \rangle$ ($Z \in \mathfrak{a}$) constituent un système d'imprimitivité transitif sur $X$, et on retrouve la proposition 3.14 de [Z96]. (Cf. aussi [K68, §2].)

(b) Le théorème contient aussi l'égalité $T^*(G/H) = \text{Ind}_H^G\{0\}$, vraie pour tout sous-groupe fermé $H$ de $G$.

(c) Si par contre on prive $T^*(G/H)$ de sa section nulle, le théorème ne s'applique plus car les gradients symplectiques des fonctions relevées de la base ne sont plus complets (cf. 2.7).

**2.25 Remarque.** Lorsque $G$ et $H$ sont connexes et simplement connexes, on peut aussi déduire (2.9) du théorème de Xu [X91, thm. 3.3] (démontré aussi par Landsman [L95, thm. 5]). Dans ce cas en effet, la double fibration

$$\begin{array}{c} \phantom{\mathfrak{h}^*} \overset{\alpha}{\swarrow} \; T^*G \; \overset{\beta}{\searrow} \phantom{\mathfrak{g}^* \times G/H} \\ \mathfrak{h}^* \phantom{\overset{\alpha}{\swarrow} \; T^*G \; \overset{\beta}{\searrow}} \mathfrak{g}^* \times G/H \end{array} \qquad (2.26)$$

($\alpha(p) = q^{-1}p_{|\mathfrak{h}}$, $\beta(p) = (pq^{-1}, qH)$, $p \in T_q^*G$) est un *bimodule d'équivalence* au sens de Xu; d'autre part, un système d'imprimitivité de base $G/H$ fournit une 'réalisation symplectique complète' $\Phi \times \pi : X \to \mathfrak{g}^* \times G/H$; enfin on peut vérifier que l'équivalence de presque-catégories construite par Xu est ici, comme il se doit, l'induction symplectique.

**D. Le cas homogène.** Supposons que $X$ soit *homogène*, ou en d'autres termes, un revêtement d'orbite coadjointe (1.7). Le critère d'inductibilité (2.6, 2.9) peut alors se traduire en termes essentiellement infinitésimaux, comme suit:

**2.27 Proposition.** *Soit* $(X, \Phi)$ *un $G$-espace hamiltonien homogène, et $H$ un sous-groupe fermé de $G$. Alors $X$ est induit par un $H$-espace hamiltonien ss'il existe $x \in X$ tel que, notant $\dot{x} = \Phi(x)$,*

(a) $H$ *contient* $G_x$;
(b) $\mathfrak{h}$ *est coisotrope en* $\dot{x}$: $\operatorname{orth}(\mathfrak{h}(\dot{x})) \subset \mathfrak{h}$;
(c) $\mathfrak{h}$ *satisfait la condition de Pukánszky en* $\dot{x}$: $\dot{x} + \operatorname{orth}(\mathfrak{h}) \subset G(\dot{x})$.

*De plus, dans ces conditions, l'espace induisant* $(Y, \Psi)$ *est aussi homogène; et si le revêtement* $\Phi$ *est injectif, alors* $\Psi$ *aussi.*

**Preuve.** Supposons que $X$ soit induite de $H$, et notons $B = G/H$. Alors on a un diagramme (2.23) et, d'après la description de la fibre en $b = eH$ qui le suit,

$$\Phi(X_b) = \{m \in \mathfrak{g}^* : m_{|\mathfrak{h}} \in \Psi(Y)\}. \qquad (2.28)$$

Soit $x \in X_b$. Comme $G$ est transitif sur $X$, on a $H(x) = X_b$. Donc l'orbite $H(\dot{x})$ égale (2.28) et contient par conséquent $\dot{x} + \operatorname{orth}(\mathfrak{h})$, d'où (c). De même son espace tangent $\mathfrak{h}(\dot{x})$ contient $\operatorname{orth}(\mathfrak{h})$, d'où (b). Enfin l'équivariance de $\pi$ assure (a).

Réciproquement, soit $x$ vérifiant (a,b,c) et montrons que $X$ admet un système d'imprimitivité de base $B = G/H$. Par (a), $g(x) \mapsto gH$ définit-bien une submersion équivariante $\pi$ de $X$ sur $B$, dont la fibre en $b = eH$ est l'orbite $X_b = H(x)$. Posons $\mathfrak{f} = \pi^*(C^\infty(B))$. Alors on a les relations

$$\mathfrak{h}(x)^\sigma \subset \mathfrak{h}(x), \qquad \mathfrak{h}(x)^\sigma = \mathfrak{f}(x), \qquad \mathfrak{f}(x) \subset \mathfrak{f}(x)^\sigma. \tag{2.29}$$

La première résulte de (b) et de ce que $\sigma(\mathfrak{h}(x), Z(x)) = \langle \mathfrak{h}(\dot{x}), Z \rangle$ (cf. 1.7); la seconde se démontre comme (2.11), et la troisième s'ensuit. Cela dit, il est clair par 'transport de structure' que la troisième relation vaut encore pour tout $x' \in X$, et les deux autres pour tout $x'$ dans la $H$-orbite $X_b$. Ceci montre que $\mathfrak{f}$ est une sous-algèbre abélienne de $C^\infty(X)$, dont les champs hamiltoniens $\eta$ sont tangents à $X_b$. Reste à voir alors que $e^{t\eta}(x)$ existe pour tout $t$. Or (2.16) montre que $\Phi$ lie $\eta_{|X_b}$ à un champ constant $\check{a} \in \text{orth}(\mathfrak{h})$, dont par (c) la courbe intégrale est dans $\Phi(X)$. Comme $\Phi$ est un revêtement cette courbe se relève à $X$, d'où la conclusion. Enfin les deux dernières affirmations résultent du lemme ci-dessous.   **q.e.d.**

**2.30 Lemme.** *Soit $(Y, \Psi)$ un $H$-espace hamiltonien.*
  (a) *Si $G$ est transitif sur $\text{Ind}_H^G Y$, alors $H$ l'est sur $Y$.*
  (b) *Si $\Phi_{\text{ind}}$ est injective, alors $\Psi$ l'est aussi.*

**Preuve.** Soient $y_1, y_2 \in Y$, et choisissons des $m_i \in \mathfrak{g}^*$ tels que $\Psi(y_i) = m_{i|\mathfrak{h}}$. Alors les orbites $[m_i, y_i] := H(m_i, y_i)$ sont dans la variété induite (1.10).

(a) Si $G$ est transitif sur celle-ci, alors $[m_1, y_1] = [gm_2, y_2]$ pour un $g \in G$; et cela signifie en particulier que $y_1$ est dans la $H$-orbite de $y_2$.

(b) Si $\Psi(y_1) = \Psi(y_2)$, on peut prendre $m_1 = m_2$. Comme $\Phi_{\text{ind}}([m_i, y_i]) = m_i$, l'injectivité de $\Phi_{\text{ind}}$ entraîne alors que $(m_1, y_1) = (m_2 h^{-1}, h(y_2))$ pour un $h \in H$. Or cette égalité exige $h = e$, et donc $y_1 = y_2$.   **q.e.d.**

**2.31 Exemple.** Tout ceci peut s'appliquer lorsque $G$ est fini ou discret. Alors les conditions (2.27b, c) sont automatiquement vérifiées, de sorte qu'un espace $G/K$ est induit de tout sous-groupe $H$ contenant $K$: (2.27) dit simplement que $G/K = \text{Ind}_H^G(H/K)$. En d'autres termes, on retrouve la notion *classique* d'imprimitivité. Cf. [H78, p. 151; Z78, ex. 2.3].

**2.32 Remarques.** (a) La proposition 2.27 est implicite dans [D84], IV(17). Sous ses hypothèses, Duflo montre en IV(11) que les représentations associées à $X$ par sa méthode sont précisément induites de celles associées à $Y$; dans certains cas, c'est même vrai *par définition*.

(b) L'inclusion (2.27b) est une égalité ssi $\mathfrak{h}$ est une *polarisation réelle*, ssi $Y$ est de dimension zéro. Plus généralement, les deux termes de (2.27b) peuvent être les sous-algèbres '$\mathfrak{d}$' et '$\mathfrak{e}$' associées à une polarisation complexe; dans ce cas (2.27) est essentiellement la proposition 3.12 de [D92].

(c) La variété $Y$ construite ici est le quotient de $H(x)$ par son feuilletage caractéristique (2.11), et un revêtement de l'orbite coadjointe $H(\dot{x}_{|\mathfrak{h}})$ (cf. 2.18). Notons cependant que la propriété $X = \operatorname{Ind}_H^G Y$ ne suffit pas à caractériser $(Y, \Psi)$: c'est le rôle de la condition (2.10b); cf. (2.36a).

**E. Morphismes.** On peut préciser (2.30a) en montrant que l'homogénéité de $Y$ sous $H$ équivaut à celle de $\operatorname{Ind}_H^G Y$ sous l'action conjointe de $G$ et $\mathcal{F}_{\text{ind}}$ (2.6). C'est en un sens l'analogue d'un autre théorème de Mackey [B62]— mais en un sens limité, car l'homogénéité n'est qu'un analogue très-imparfait de l'irréductibilité (cf. 2.36b).

La proposition suivante offre un meilleur analogue, et éclaire du même coup le rôle de la condition (2.10b). Étant donnés des $G$-espaces hamiltoniens $X_1, X_2$ munis de systèmes d'imprimitivité de même base $B$ (cf. 2.4c), convenons d'appeler $\operatorname{Hom}_B(X_1, X_2)$ l'ensemble des $J \in \operatorname{Hom}_G(X_1, X_2)$ tels que

$$\begin{array}{c} X_1 \xrightarrow{J} X_2 \\ {}_{\pi_1}\searrow \quad \swarrow{}_{\pi_2} \\ B \end{array} \qquad (2.33)$$

commute.

**2.34 Proposition.** $\operatorname{Hom}_{G/H}(\operatorname{Ind}_H^G Y_1, \operatorname{Ind}_H^G Y_2) = \operatorname{Hom}_H(Y_1, Y_2)$.

**Preuve.** Notons $X_1, X_2$ les variétés induites par $Y_1, Y_2$, et $\pi_1, \pi_2$ leurs projections (2.5) sur $G/H$. Si $j \in \operatorname{Hom}_H(Y_1, Y_2)$, alors visiblement le difféomorphisme symplectique

$$\operatorname{id} \times j : T^*G \times Y_1 \to T^*G \times Y_2 \qquad (2.35)$$

passe aux quotients (1.10), et y définit un élément $J$ de $\operatorname{Hom}_{G/H}(X_1, X_2)$. Réciproquement, tout $J \in \operatorname{Hom}_{G/H}(X_1, X_2)$ applique $\pi_1^{-1}(eH)$ sur $\pi_2^{-1}(eH)$ en respectant les 2-formes et donc leurs feuilletages caractéristiques. Il en résulte un isomorphisme entre les espaces quotient qui, d'après (2.9), ne sont autres que $Y_1$ et $Y_2$. Enfin on vérifie sans peine que les correspondances ainsi définies sont inverses l'une de l'autre.    ***q.e.d.***

**2.36 Exemples.** La proposition autorise $\mathrm{Hom}_G(\mathrm{Ind}_H^G Y_1, \mathrm{Ind}_H^G Y_2)$ à être strictement plus grand que $\mathrm{Hom}_H(Y_1, Y_2)$. Ainsi:

(a) *Des espaces non isomorphes peuvent induire des espaces isomorphes.* Soit $G = \mathbf{U}(1) \ltimes \mathbf{C}$ le groupe des déplacements du plan, et $\Re, \Im$ les formes linéaires 'partie réelle' et 'partie imaginaire' considérées comme orbites coadjointes du sous-groupe $\mathbf{C}$. Alors $\mathrm{Ind}_{\mathbf{C}}^G \Re$ et $\mathrm{Ind}_{\mathbf{C}}^G \Im$ sont une seule et même orbite coadjointe de $G$ (un cylindre).

Dans cet exemple, $Y_1$ et $Y_2$ sont au moins isomorphes comme $G$-variétés (en oubliant les applications moment); mais cela n'est pas même nécessaire. En effet, (2.31) fournira un contre-exemple chaque fois que $G$ est un groupe fini, $H$ un sous-groupe invariant, et $K_1, K_2$ des sous-groupes de $H$ conjugués sous $G$ mais pas sous $H$. Ceci se réalise facilement dans les matrices triangulaires strictes sur $\mathbf{Z}/2\mathbf{Z}$.

(b) *Un espace sans automorphismes peut induire un espace qui en a.* Soit $G_1 = \mathbf{R} \ltimes \mathbf{C}$ le revêtement universel du précédent. Alors $\mathrm{Ind}_{\mathbf{C}}^{G_1} \Re$ est le revêtement universel du cylindre précédent, dont l'homotopie fournit des automorphismes non triviaux. Leur présence reflète le fait que la représentation $\mathrm{Ind}_{\mathbf{C}}^{G_1}(e^{i\Re})$ est réductible [B72, p. 189], bien que $\mathrm{Ind}_{\mathbf{C}}^{G_1} \Re$ soit homogène sous $G_1$.

## 3. Le théorème du petit groupe

Pour appliquer le théorème d'imprimitivité, il nous faut des systèmes d'imprimitivité (transitifs). En théorie des représentations, l'observation fondamentale de Frobenius et Mackey est que la décomposition de la restriction d'une représentation (irréductible) de $G$ à un sous-groupe fermé normal, $N$, fournit un tel système.

L'analogue pour un $G$-espace hamiltonien $X$ s'impose de lui-même. Puisque $N$ est normal, $G$ agit naturellement dans $\mathfrak{n}$ et $\mathfrak{n}^*$, et respecte la partition de $\mathfrak{n}^*$ en $N$-orbites. On a donc une action de $G$ dans l'espace d'orbites $\mathfrak{n}^*/N$, et une suite d'applications $G$-équivariantes

$$X \xrightarrow{\Phi} \mathfrak{g}^* \longrightarrow \mathfrak{n}^* \longrightarrow \mathfrak{n}^*/N. \tag{3.1}$$

On peut alors s'attendre à ce que l'image de $X$ tout à droite soit la base d'un système d'imprimitivité, et à ce que l'application de (2.9) à ce système donne l'analogue suivant du théorème de Mackey [F88, p. 1284]. Dans cet énoncé, $\Phi_N$ désigne la composée des deux premières applications (3.1).

**3.2 Théorème.** *Soit $(X, \Phi)$ un $G$-espace hamiltonien homogène, et $N \subset G$ un sous-groupe fermé normal. Soit $W \in \mathfrak{n}^*/N$ une des orbites contenues dans $\Phi(X)_{|\mathfrak{n}}$, et supposons que son stabilisateur $H := G_W$ est fermé dans $G$. Alors $N \subset H$, et il existe un unique $H$-espace hamiltonien $(Y, \Psi)$ tel que*

$$\text{(a)} \quad X = \operatorname{Ind}_H^G Y, \qquad \text{(b)} \quad \Psi(Y)_{|\mathfrak{n}} = W. \tag{3.3}$$

*Explicitement $Y$ est le quotient de $\Phi_N^{-1}(W)$ par son feuilletage caractéristique, i.e., c'est l'espace réduit de $X$ en $W$ au sens de KKS* [K78].

**Preuve.** L'inclusion $N \subset H$ est claire, car $N$ agit trivialement dans $\mathfrak{n}^*/N$.

Soit $B$ l'image de $X$ par (3.1). C'est une $G$-orbite dans $\mathfrak{n}^*/N$, et comme le stabilisateur $H$ du point $b := W$ est fermé, c'est naturellement une variété. Notant $\pi$ la projection $X \to B$, nous vérifierons ci-dessous que

$$\pi^*(C^\infty(B)) \text{ est un système d'imprimitivité sur } X. \tag{3.4}$$

Montrons d'abord comment le théorème en découle. Pour l'*existence*, soit $(Y, \Psi)$ le $H$-espace hamiltonien que nous fournit le théorème 2.9. Il satisfait (a), et nous savons qu'il est le quotient de $X_b := \pi^{-1}(b) = \Phi_N^{-1}(W)$ par son feuilletage caractéristique (2.11). De plus, le diagramme (2.18) montre que $\Psi(Y)_{|\mathfrak{n}} = \Phi(X_b)_{|\mathfrak{h}|\mathfrak{n}} = W$, d'où (b).

Réciproquement, soit $(Y, \Psi)$ satisfaisant (a,b), et $J \in \operatorname{Hom}_G(X, \operatorname{Ind}_H^G Y)$. Alors la partie *unicité* résultera aussi du théorème 2.9, si nous montrons que $J$ envoie nécessairement $X_b$ dans $\pi_{\operatorname{ind}}^{-1}(eH)$. Soit donc $x \in X_b$, et notons $\dot{x}$, $w$ et (donc) $W$ ses images successives par (3.1), de sorte qu'on a

$$x \in X, \qquad \dot{x} = \Phi(x), \qquad w = \dot{x}_{|\mathfrak{n}}, \qquad W = N(w). \tag{3.5}$$

Soit d'autre part $(p, y) \in T_q^*G \times Y$ un point de la $H$-orbite qu'est $J(x)$ (1.10). Comme $J^*\Phi_{\text{ind}} = \Phi$, les applications moment (1.11) et (1.9) prennent en ce point les valeurs $\dot{x}$ et 0, ce qui s'écrit: $pq^{-1} = \dot{x}$ et $q^{-1}p_{|\mathfrak{h}} = \Psi(y)$. Substituant la première de ces relations dans la seconde, il vient $q^{-1}(\dot{x})_{|\mathfrak{h}} = \Psi(y)$. Comme par hypothèse $\dot{x}$ et $\Psi(y)$ se projettent tous deux sur $W$, on en déduit que $q^{-1}(W) = W$ et donc $q \in H$, ce qu'il fallait démontrer (cf. 2.5).

Reste donc à vérifier (3.4). Mais cela résultera de la preuve de (2.27), si nous montrons que $H$ satisfait les trois conditions de cette proposition; or c'est ce qu'assure le lemme suivant. **q.e.d.**

**3.6 Lemme [P78, Lemma 2].** *Gardons les hypothèses du théorème, et soient $x$, $\dot{x}$, $w$ comme en (3.5). Alors*

 (a) *$H$ contient $G_x$,*
 (b) *$\mathfrak{n}_w(\dot{x}) = \text{orth}(\mathfrak{h})$,*
 (c) *$N_w^0(\dot{x}) = \dot{x} + \text{orth}(\mathfrak{h})$,*

*où $N_w^0$ est la composante neutre du stabilisateur $N_w$.*

**Preuve.** (a) est clair par équivariance de (3.1). (b): Il résulte de la définition des stabilisateurs $H$, $G_w$, $N_w$, et de calculs ensemblistes élémentaires, que

  i. $H = NG_w$,
  ii. $\mathfrak{g}_w = \text{orth}(\mathfrak{n}(\dot{x}))$,
  iii. $\mathfrak{n}_w = \{Z \in \mathfrak{n} : Z(\dot{x}) \in \text{orth}(\mathfrak{n})\}$,

d'où $\text{orth}(\mathfrak{h}) \stackrel{\text{i}}{=} \text{orth}(\mathfrak{n} + \mathfrak{g}_w) = \text{orth}(\mathfrak{g}_w) \cap \text{orth}(\mathfrak{n}) \stackrel{\text{ii}}{=} \mathfrak{n}(\dot{x}) \cap \text{orth}(\mathfrak{n}) \stackrel{\text{iii}}{=} \mathfrak{n}_w(\dot{x})$.
(c): La relation (b) dit que $\mathfrak{n}_w$ stabilise $\dot{x}_{|\mathfrak{h}}$. Donc $N_w^0$ aussi, ce qui signifie que $N_w^0(\dot{x}) \subset \dot{x} + \text{orth}(\mathfrak{h})$. Pour l'inclusion réciproque, remarquons que comme $\mathfrak{n}$ est un idéal on a $\text{ad}(Z)^n(\mathfrak{g}) \subset [\mathfrak{n}_w, \mathfrak{n}]$ pour tout $Z \in \mathfrak{n}_w$ et tout $n \geq 2$. Puisque $w = \dot{x}_{|\mathfrak{n}}$ s'annule sur $[\mathfrak{n}_w, \mathfrak{n}]$, ceci entraîne la formule

$$\begin{aligned}\langle \exp(Z)(\dot{x}), Z'\rangle &= \Big\langle \dot{x}, \sum_{n=0}^{\infty} \frac{(-)^n}{n!}\text{ad}(Z)^n(Z')\Big\rangle \\ &= \langle \dot{x}, Z' - [Z, Z']\rangle \\ &= \langle \dot{x} + Z(\dot{x}), Z'\rangle \end{aligned} \quad (3.7)$$

pour tous $Z \in \mathfrak{n}_w$, $Z' \in \mathfrak{g}$. En particulier $\dot{x} + \mathfrak{n}_w(\dot{x})$ est bien contenu dans $N_w^0(\dot{x})$, ce qu'il fallait démontrer. **q.e.d.**

Le théorème ci-dessus n'est pas tout à fait aussi détaillé que le modèle de [F88] que nous avons invoqué. Pour compléter le parallèle, on peut lui adjoindre l'énoncé suivant:

**3.8 Proposition.** *Soit $N \subset G$ un sous-groupe fermé normal, et $W \in \mathfrak{n}^*/N$ une orbite telle que $H := G_W$ soit fermé. Alors $Y \mapsto X = \mathrm{Ind}_H^G Y$ définit une bijection entre (classes d'isomorphisme de)*

  (a) *$G$-espaces hamiltoniens homogènes $(X, \Phi)$ tels que $W \subset \Phi(X)_{|\mathfrak{n}}$;*
  (b) *$H$-espaces hamiltoniens homogènes $(Y, \Psi)$ tels que $W = \Psi(Y)_{|\mathfrak{n}}$.*

*De plus, dans ces conditions, le revêtement $\Phi$ est trivial ssi $\Psi$ l'est.*

**Preuve.** Le théorème 3.2 fournit un 'inverse-à-droite' (a) → (b) de l'induction symplectique, avec $\Phi$ trivial ⇒ $\Psi$ trivial (2.30). Reste à voir que cet inverse est surjectif, et ⇐.

Soit donc $(Y, \Psi)$ du type (b), et montrons que $G$ est transitif sur $\mathrm{Ind}_H^G Y$. Pour cela, choisissons $y \in Y$, posons $w = \Psi(y)_{|\mathfrak{n}}$, et fixons un $m \in \mathfrak{g}^*$ tel que $m_{|\mathfrak{h}} = \Psi(y)$. Appliquons alors le lemme 3.6 en y remplaçant la donnée $(G, X, N, W)$ par $(H, Y, N, W)$, resp. $(G, G(m), N, W)$. En tenant compte du fait que $\Psi$ est un revêtement, on obtient les relations

$$N_w^0(y) = \{y\}, \qquad \text{resp.} \qquad N_w^0(m) = m + \mathrm{orth}(\mathfrak{h}). \tag{3.9}$$

Cela étant, soit $(qm_1, y_1) \in T_q^*G \times Y$ un représentant d'un élément quelconque de la variété induite (1.10). Choisissons un $h \in H$ tel que $y_1 = h(y)$, et posons $m_2 = h^{-1}(m_1)$. Alors $(qhm_2, y)$ est un autre représentant du même élément. De plus, la nullité de (1.9) entraîne que $m_{2|\mathfrak{h}} = \Psi(y) = m_{|\mathfrak{h}}$, et (3.9) montre alors que $m_2 = n(m)$ pour un $n \in N_w^0 \subset H_y$. On en déduit que $(qhnm, y)$ est un troisième représentant du même élément. Celui-ci est donc l'image par $g = qhn$ de celui défini par $(m, y)$, d'où la transitivité.

Supposons maintenant que le revêtement $\Psi$ soit trivial, i.e., que l'inclusion $H_y \subset H_{\Psi(y)}$ soit une égalité. Notant $[m, y]$ l'élément de $\mathrm{Ind}_H^G Y$ dont on vient de parler (et que $\Phi_{\mathrm{ind}}$ envoie visiblement sur $m$), il nous reste à montrer que l'inclusion $G_{[m,y]} \subset G_m$ est alors elle aussi une égalité. Or $m$ a été choisi tel que

$$m_{|\mathfrak{h}} = \Psi(y) \qquad \text{et} \qquad m_{|\mathfrak{n}} = w. \tag{3.10}$$

La deuxième relation implique que $G_m \subset G_w \subset H$, et la première que de plus $G_m \subset H_{\Psi(y)} = H_y$. On en conclut que si $g \in G_m$ alors $g([m,y]) = [gm, y] = [mg, y] = [m, g(y)] = [m, y]$, d'où $G_m \subset G_{[m,y]}$. **q.e.d.**

**3.11 Exemples.** (a) Lorsque $N$ est abélien, ou que $W$ est un point, (3.2) se réduit encore à la proposition rappelée en (2.24a). L'application aux cas classiques des groupes d'Euclide et Poincaré est détaillée dans [Z96, §3.3].

(b) Si $G$ est résoluble exponentiel, nous avons vu dans [Z96, §4.1] comment on peut itérer (3.2) jusqu'à ce que $X$ apparaisse comme induite d'un point; c'est l'analogue symplectique du théorème de Takenouchi selon lequel les représentations unitaires irréductibles de $G$ sont monomiales [B72, p. 140]. Cela redémontre en somme l'existence, dans un tel groupe, de polarisations réelles vérifiant la condition de Pukánszky (cf. 2.32b).

(c) L'hypothèse de (3.2) et (3.8) selon laquelle $G_W$ est fermé n'est pas gratuite, car $\mathfrak{n}^*/N$ n'est pas forcément séparé. De fait, elle est en défaut dans des cas bien connus comme le suivant: Soit $G$ le groupe des matrices de la forme
$$g = \begin{pmatrix} a & b \\ 0 & c \end{pmatrix} \tag{3.12}$$
où $a$ est prise dans le tore diagonal $T \subset \mathbf{U}(2)$, $b$ diagonale dans $\mathfrak{gl}(2, \mathbf{C})$, et $c$ dans une droite irrationnelle $S \subset T$. Si $N$ désigne le sous-groupe $a = \mathbf{1}$, et $W$ l'orbite de la valeur en l'élément neutre de la 1-forme $\mathrm{Re}(\mathrm{Tr}(db))$, alors on trouve $G_W = \{g \in G : a \in S\}$, sous-groupe dense de $G$.

**3.13 Remarques.** Lorsque $N$ est le groupe additif d'un espace vectoriel, et $G$ un produit semidirect $K \ltimes N$ où $K$ agit linéairement dans $N$, le théorème 3.2 est dû à V. Guillemin et S. Sternberg [G83, §4].

Notons que le théorème d'imprimitivité (2.9)—dont nous avons déduit (3.2)—rentre lui-même formellement dans ce cadre si on considère le produit semidirect $G \ltimes \mathcal{F}$ défini par l'action naturelle de $G$ sur $\mathcal{F}$ (§2A). De ce point de vue, *toute* variété induite l'est pour des raisons de 'normal subgroup analysis'; comparer à ce propos la question posée dans [F88, p. 1371].

Par ailleurs, le lecteur pourra constater que des variantes du lemme 3.6 interviennent *très* souvent dans la littérature consacrée à la méthode des orbites; et qu'avoir alors à l'esprit le théorème du petit groupe est rarement sans éclairer la situation.

## 4. Questions de décomposition

Les représentations d'un groupe $H$, dont la restriction à un sous-groupe normal $N$ est multiple d'une irréductible donnée $S$, s'obtiennent comme suit: on étend $S$ en une représentation projective de $H$, et on tensorise avec une représentation projective de $H/N$ qui 'tue le cocycle'. C'est l''étape III' de l'analyse de Mackey [F88, p. 1263].

La situation analogue se présente en (3.8b): un $H$-espace hamiltonien $Y$, dont l'image dans $\mathfrak{n}^*$ se réduit à une seule orbite $W$. Peut-on en tirer un théorème de structure pour $Y$? Nous nous contenterons ici de quelques remarques. La première est que $H$ agit naturellement sur $W$ (puisqu'il la stabilise dans $\mathfrak{n}^*/N$), et que $H$ et même $H/N$ agissent sur l'espace d'orbites $Y/N$ (puisque $H$ normalise $N$). On a donc un diagramme $H$-équivariant

$$\begin{array}{ccc} & Y & \\ {}^{\alpha}\swarrow & & \searrow{}^{\beta} \\ Y/N & & W. \end{array} \qquad (4.1)$$

D'autre part, on déduit sans peine de (1.8a) et (3.9a) que les fibres de $\alpha$ et $\beta$ sont des sous-variétés symplectiques de $Y$, supplémentaires et orthogonales l'une de l'autre relativement à la 2-forme de $Y$. *Si $Y/N$ est une variété*, elle admettra une structure symplectique liée à celles de $Y$ et $W$ par

$$\sigma_Y = \alpha^* \sigma_{Y/N} + \beta^* \sigma_W, \qquad (4.2)$$

et $\alpha \times \beta$ fera de $Y$ un revêtement symplectique du produit $Y/N \times W$, chaque fibre de $\alpha$ étant un revêtement de $W$ et vice versa. Enfin, on peut vérifier que l'action de $H$ sur $W$ préserve $\sigma_W$. *Si elle est hamiltonienne*, il en résultera selon (1.2) une extension centrale de $\mathfrak{h}$ (en fait de $\mathfrak{h}/\mathfrak{n}$), nécessairement compensée par $Y/N$ dans (4.2). C'est, infinitésimalement, l'obstruction de Mackey.

La *décomposition barycentrique* de [S69] est de ce type, en prenant pour $N$ un groupe de Heisenberg. Notons aussi que lorsque $Y$ provient d'un $G$-espace hamiltonien comme en (3.8), le fibration de $Y$ (réduite KKS) en revêtements de $Y/N$ (réduite de Marsden-Weinstein) a déjà été observée dans [M76], [K78], et [G84, 26.6].

## 5. Groupes réductifs

L'analyse de Mackey ne s'applique guère aux groupes simples, qui n'ont pas de sous-groupes normaux, ni plus généralement aux groupes semisimples ou réductifs; et pourtant l'induction ('parabolique') joue un rôle dans la construction de leurs représentations.

Voyons comment ceci se traduit géométriquement. Suivant Vogan [V87], nous appellerons *G réductif* s'il y a un homomorphisme $\eta : G \to \mathbf{GL}(n, \mathbf{R})$ de noyau fini et d'image $\Theta$-stable, où $\Theta(g) = {}^t g^{-1}$. Alors $\mathfrak{g}$ s'identifie à une algèbre de Lie de matrices, $\mathfrak{g}^*$ s'identifie à $\mathfrak{g}$ au moyen de la forme trace $\langle Z, Z' \rangle = \mathrm{Tr}(ZZ')$, et toute $x \in \mathfrak{g}^*$ a une unique décomposition de Jordan

$$x = x_\mathrm{h} + x_\mathrm{e} + x_\mathrm{n}, \qquad x_\mathrm{h}, x_\mathrm{e}, x_\mathrm{n} \in \mathfrak{g}^*, \tag{5.1}$$

où la matrice $x_\mathrm{h}$ est hyperbolique (diagonalisable à valeurs propres réelles), $x_\mathrm{e}$ elliptique (diagonalisable à valeurs propres imaginaires) et $x_\mathrm{n}$ nilpotente, et où $x_\mathrm{h}$, $x_\mathrm{e}$ et $x_\mathrm{n}$ commutent entre elles. (Pour tout ceci, voir [V87].)

Cela étant, on a le

**5.2 Théorème.** *Soit $X = G(x)$ une orbite coadjointe du groupe réductif $G$. Notons $\mathfrak{n}$ la somme des espaces propres associés aux valeurs propres positives de $\mathrm{ad}(x_\mathrm{h})$, et $P$ le normalisateur de $\mathfrak{n}$ dans $G$. Alors on a*

$$X = \mathrm{Ind}_P^G Y \qquad \text{où} \qquad Y = P(x_{|\mathfrak{p}}). \tag{5.3}$$

*De plus on a $Y_{|\mathfrak{n}} = \{0\}$, de sorte que $Y$ provient d'une orbite du quotient $P/\exp(\mathfrak{n})$.*

**Preuve.** Ceci résultera de la proposition 2.27, et de la remarque 2.32c qui la suit, si nous montrons que

(a) $P$ contient $G_x$,
(b) $\mathfrak{n}(x) = \mathrm{orth}(\mathfrak{p})$,
(c) $N(x) = x + \mathrm{orth}(\mathfrak{p})$, où $N := \exp(\mathfrak{n})$.

La première relation est claire: si $g \in G_x$, alors $\eta(g)$ commute avec $x$ et donc avec $x_{\mathrm{h}}$, par une propriété bien connue de la décomposition de Jordan; en particulier $g$ préserve les espaces propres de $\mathrm{ad}(x_{\mathrm{h}})$ et donc $\mathfrak{n}$, d'où $g \in P$.

Pour les deux autres, il nous faut rappeler quelques éléments de [V87]. Notons $\mathfrak{g}^a$ le sous-espace propre de $\mathfrak{g}$ pour la valeur propre $a$ de $\mathrm{ad}(x_{\mathrm{h}})$. Alors $\mathfrak{g} = \bigoplus_{a \in \mathbf{R}} \mathfrak{g}^a$, et par définition $\mathfrak{n} = \bigoplus_{a > 0} \mathfrak{g}^a$. L'identité de Jacobi montre que
$$[\mathfrak{g}^a, \mathfrak{g}^b] \subset \mathfrak{g}^{a+b}, \tag{5.4}$$
et l'invariance de la forme trace montre que $0 = \langle [x_{\mathrm{h}}, Z], Z' \rangle + \langle Z, [x_{\mathrm{h}}, Z'] \rangle = (a+b)\langle Z, Z' \rangle$ pour $(Z, Z') \in \mathfrak{g}^a \times \mathfrak{g}^b$. Donc
$$\langle \cdot, \cdot \rangle_{|\mathfrak{g}^a \times \mathfrak{g}^b} \text{ est } \begin{cases} \text{nulle} & \text{si } a + b \neq 0, \\ \text{non dégénérée} & \text{si } a + b = 0. \end{cases} \tag{5.5}$$

Cela dit, il résulte de (5.4) que $\mathfrak{p}$, normalisateur de $\mathfrak{n}$, égale $\mathfrak{g}^0 \oplus \mathfrak{n}$, et (5.5) montre alors que $\mathrm{orth}(\mathfrak{p}) = \mathfrak{n}$. Donc (b) revient à voir que $\mathfrak{n}(x) = \mathfrak{n}$, ou en d'autres termes, que $\mathrm{ad}(x)$ envoie $\mathfrak{n}$ sur $\mathfrak{n}$. En effet il l'envoie *dans* $\mathfrak{n}$ par (5.4) car $x \in \mathfrak{g}^0$, et *sur* $\mathfrak{n}$ car $\mathrm{Ker}(\mathrm{ad}(x)) = \mathfrak{g}_x \subset \mathfrak{g}_{x_{\mathrm{h}}} = \mathfrak{g}^0$.

De (b) on déduit d'abord que $\mathfrak{n}$ stabilise $x_{|\mathfrak{p}}$, d'où $N(x) \subset x + \mathrm{orth}(\mathfrak{p})$, et ensuite que $N(x)$ est une sous-variété de dimension maximale et donc ouverte de $x + \mathrm{orth}(\mathfrak{p})$. Mais $N$ est un groupe nilpotent, par (5.4); donc $N(x)$ est aussi fermée [B72, p. 7], d'où (c).

Quant à la dernière assertion de l'énoncé, il suffit de noter que $Y_{|\mathfrak{n}}$ est la $P$-orbite de $x_{|\mathfrak{n}}$ dans $\mathfrak{n}^*$, et que $x_{|\mathfrak{n}} = 0$ par (5.5) puisque $x \in \mathfrak{g}^0$.　　**q.e.d.**

**5.6 Remarques.** (a) $P$ est un sous-groupe parabolique de $G$. Lorsqu'il est *minimal* (sous-groupe de Borel), et que $Y$ satisfait une certaine condition de régularité, notre théorème est dû à Guillemin et Sternberg [G83, §3], et l'égalité (c) ci-dessus à Harish-Chandra [H54, Lemma 8].

(b) Le programme de la méthode des orbites est, ici encore, de définir les représentations associées à $X$ comme induites de celles associées à $Y$. Cependant ces dernières ne sont actuellement connues, grâce au foncteur de Zuckerman, que dans le cas où $x_{\mathrm{n}} = 0$.

## 6. Application aux représentations quantiques

> *"Comme le Dieu des philosophes, l'opération de Hilbert est incompréhensible et ne se voit pas; mais elle gouverne tout, et ses manifestations sensibles éclatent partout."* —R. Godement

**A. Prologue.** Autorisons-nous une minute de spéculation.

Tout se passe comme s'il existait un foncteur, qui à toute variété symplectique associerait une famille de représentations projectives de son groupe d'automorphismes. Le foncteur se manifesterait à nous au travers de représentations restreintes, qui seraient des multiples de celles qu'ont découvertes Maxwell, Dirac, Kirillov ...

Personne n'a jamais vu ce foncteur, et certains ont même prétendu qu'il n'existait pas. (C'est, selon d'autres, qu'ils exigeaient de lui des vertus un peu exagérées.)

Mais supposons un instant qu'il existe, et qu'il vérifie l'axiome proposé par J. M. Souriau dans [S88]—sous la forme légèrement renforcée de [Z96]. Prenons pour variété symplectique un $G$-espace hamiltonien, et supposons qu'il soit induit d'un sous-groupe fermé, $H$. Notons $U$ la restriction à $G$ d'une des représentations obligeamment fournies par notre foncteur.

Alors $U$ est aussi induite de $H$, et le but du chapitre est de le démontrer. La méthode consistera à adjoindre à $G$ le groupe $\mathcal{F}$ défini par le système d'imprimitivité associé (2.6), et à tirer les conséquences de l'axiome pour la représentation de $\mathcal{F} \rtimes G$ dont on peut supposer—si foncteur il y a—que $U$ est aussi la restriction.

**B. Représentations 0-quantiques.** L'axiome en question s'exprime au travers de la définition suivante. Soit $(X, \Phi)$ un $G$-espace hamiltonien, et convenons d'appeler *$X$-abéliennes* les sous-algèbres $\mathfrak{a}$ de $\mathfrak{g}$ telles que $\{H_Z, H_{Z'}\} = 0$ chaque fois que $Z, Z' \in \mathfrak{a}$ (cf. 1.3).

**6.1 Définition.** Une *représentation 0-quantique* pour $(X, \Phi)$ est une représentation unitaire $U$ de $G$, telle que pour toute sous-algèbre $X$-abélienne $\mathfrak{a}$ de $\mathfrak{g}$, $U \circ \exp_{|\mathfrak{a}}$ soit une représentation de $(\mathfrak{a}, +)$ concentrée sur $\Phi(X)_{|\mathfrak{a}}$.

*Concentrée* est à prendre au sens suivant. On regarde le dual $\mathfrak{a}^*$ comme une partie du dual de Pontrjagin $\hat{\mathfrak{a}}$ du groupe additif *discret* $\mathfrak{a}$, en identifiant $p \in \mathfrak{a}^*$ au caractère $e^{i\langle p,\cdot\rangle}$ de $\mathfrak{a}$; et on demande que la mesure spectrale $E$ définie sur $\hat{\mathfrak{a}}$ par le théorème de Stone soit concentrée sur $\Phi(X)_{|\mathfrak{a}}$, ce qui veut dire que le complémentaire de cet ensemble est négligeable pour la mesure de Radon $(\varphi, E(\cdot)\varphi)$, quel que soit $\varphi$ dans l'espace de la représentation. (Pour tout ceci, voir [Z96, §1.1F].)

Remarquons que nous ne supposons pas $U$ continue a priori. Quant au préfixe '0-', il vient du fait que nous utilisons ici la définition de [Z96], qui renforce un peu celle de Souriau. Cette dernière reviendrait en effet, comme on l'a vu dans [Z96], à supposer $U \circ \exp_{|\mathfrak{a}}$ concentrée non pas sur $\Phi(X)_{|\mathfrak{a}}$ mais sur sa fermeture dans $\hat{\mathfrak{a}}$.

**C. Représentations 0-quantiques de $\mathcal{F} \rtimes G$.** Soit maintenant $H$ un sous-groupe fermé de $G$, et supposons que $X$ soit induite d'un $H$-espace hamiltonien: $X = \mathrm{Ind}_H^G Y$. On a donc sur $X$ un système d'imprimitivité $\mathfrak{f} = \mathfrak{f}_{\mathrm{ind}}$ (2.6), dont nous avons noté $\mathcal{F}$ le groupe additif sous-jacent (§2A); et on peut former le produit semidirect $\mathcal{F} \rtimes G$:

$$(f,g)(f',g') = (f + g(f'), gg'). \tag{6.2}$$

Bien que $\mathcal{F} \rtimes G$ ne soit pas en général un groupe de Lie, tous les ingrédients nécessaires pour appliquer la définition 6.1 sont réunis: on a une action de $\mathcal{F} \rtimes G$ sur $X$ (§2E), et une application moment équivariante

$$\pi \times \Phi : X \to \mathfrak{f}^* \times \mathfrak{g}^*. \tag{6.3}$$

Par ailleurs, toute représentation de $\mathcal{F} \rtimes G$ s'écrit $(f,g) \mapsto D(f)U(g)$ pour une représentation $D$ de $\mathcal{F}$ et une représentation $U$ de $G$. Cela étant, on peut énoncer le

**6.4 Théorème.** *Soit $(D, U)$ une représentation 0-quantique de $\mathcal{F} \rtimes G$ pour $X = \mathrm{Ind}_H^G Y$, et fixons $x \in X$. Alors il existe une unique représentation unitaire continue $T$ de $H$, telle que*

$$\text{(a)} \quad U = \mathrm{Ind}_H^G T, \qquad \text{(b)} \quad (D(f)\varphi)(g) = e^{if(g(x))}\varphi(g). \tag{6.5}$$

**Remarque.** L'unicité de $T$ s'entend à équivalence unitaire près. De même, par (6.5) il faut bien sûr entendre que le couple $(D, U)$ est unitairement équivalent à celui défini par (a, b); rappelons à ce propos que la représentation induite agit dans un espace de (classes de) fonctions $\varphi$ sur $G$.

**Preuve.** La représentation $U$ est continue par [Z96, 2.18]. Appliquons d'autre part (6.1) à la sous-algèbre abélienne $\mathfrak{f}$ (sur laquelle exp sera, bien entendu, l'identité $\mathfrak{f} \to \mathcal{F}$). La mesure spectrale $E$ de $D \circ \exp_{|\mathfrak{f}} = D$ est définie par $E(f^\vee) = D(f)$, où $f^\vee$ désigne le caractère continu $\chi \mapsto \chi(f)$ du groupe abélien compact $\hat{\mathfrak{f}}$. Comme $(D, U)$ est une représentation, on a $D(g(f)) = U(g)D(f)U(g^{-1})$ et donc

$$E(g(f^\vee)) = U(g)E(f^\vee)U(g^{-1}). \tag{6.6}$$

De plus, $E$ est par hypothèse concentrée sur la base $\pi(X) \subset \mathfrak{f}^* \subset \hat{\mathfrak{f}}$ de $\mathfrak{f}$. Par conséquent, le théorème d'imprimitivité de Mackey [F88, p. 1194] s'applique, et nous fournit une unique représentation continue $T$ de $H$ telle que

$$\text{(a)} \quad U = \text{Ind}_H^G T, \qquad \text{(b)} \quad (E(f^\vee)\varphi)(g) = f^\vee(g(b))\varphi(g), \tag{6.7}$$

où $b = \pi(x)$ vu comme élément de $\hat{\mathfrak{f}}$. Or ceci n'est autre que (6.5). **q.e.d.**

# Références

# ANNEXES

Table

---

N.B. Tous les renvois faits dans les annexes le sont au premier article, pp. 1–66.

# Une propriété des mesures de Radon

Nous démontrons ci-dessous la proposition 1.31b (p. 13). En effet, elle ne semble se trouver que dans les exercices de Bourbaki [k], avec de surcroît un énoncé erroné. (Il est trivialement en défaut quand l'application $\pi$ considérée n'est pas surjective.) Voici un énoncé correct:

**Proposition.** *Soit $\pi : X \to Y$ une application continue entre espaces localement compacts, où $X$ est dénombrable à l'infini, et soit $\nu$ une mesure positive sur $Y$. Si $\nu$ est concentrée sur $\pi(X)$, alors il existe une mesure positive $\mu$ sur $X$ telle que $\nu = \pi(\mu)$.*

On a besoin de cette proposition en (2.18), et dans l'annexe B ci-dessous. L'hypothèse de dénombrabilité est essentielle, comme on le voit en prenant $X = \mathbf{R}_{\text{discret}}$, $Y = \mathbf{R}$, $\nu =$ Lebesgue.

**Preuve.** Par hypothèse, $X$ est la réunion d'une famille dénombrable de compacts $V_n$. Remplaçant chaque $V_n$ par un voisinage compact $W_n$, puis par $X_n = W_1 \cup \cdots \cup W_n$, on obtient une suite croissante de compacts $X_n$ dont les intérieurs recouvrent $X$. Posons $Y_n = \pi(X_n)$, et

$$\pi_n = \text{l'application } X_n \to Y_n \text{ déduite de } \pi,$$
$$\nu_n = \nu|Y_n, \text{ la mesure induite par } \nu \text{ sur } Y_n \text{ [c]}.$$

**Lemme 1.** *Il existe une mesure positive $\mu_n$ sur $X_n$, telle que $\nu_n = \pi_n(\mu_n)$.*

En effet, soit $F_n = \{\, g \circ \pi_n : g \in \mathrm{C}(Y_n) \}$, où $\mathrm{C}(Y_n)$ désigne les fonctions continues sur $Y_n$. Comme $g \circ \pi_n$ détermine $g$, il y a une unique forme linéaire (positive) $\lambda_n$ sur $F_n$ telle que $\nu_n(g) = \lambda_n(g \circ \pi_n)$. Par Hahn-Banach, ou [b], on peut étendre $\lambda_n$ par une forme linéaire positive $\mu_n$ sur $\mathrm{C}(X_n)$.

**Lemme 2.** *On peut choisir les $\mu_n$ de sorte que $\mu_n|X_m = \mu_m$ pour $n \geq m$.*

Par récurrence, il suffit de voir que $\mu_1$ étant donnée, on peut modifier $\mu_2$ de sorte qu'elle induise $\mu_1$ sur $X_1$. Pour cela, il suffit de remplacer $\mu_2$ par

$$\tilde{\mu}_2 = i(\mu_1) + \phi.\mu_2,$$

où $i$ est l'injection $X_1 \hookrightarrow X_2$, et $\phi$ la fonction caractéristique de $X_2 \setminus \pi_2^{-1}(Y_1)$. C'est bien défini, car $i$ est propre et $\phi$ est $\mu_2$-intégrable. Notant $j$ l'injection $Y_1 \hookrightarrow Y_2$, et $\psi$ la fonction caractéristique de $Y_2 \setminus Y_1$, on a $\pi_2 \circ i = j \circ \pi_1$ et $\phi = \psi \circ \pi_2$. Par [f, i, j], on a donc encore

$$\pi_2(\tilde{\mu}_2) = j(\nu_1) + \psi.\nu_2 = \nu_2.$$

D'autre part, $\phi$ étant nulle sur $X_1$, [h] et [i] donnent $\tilde{\mu}_2|X_1 = \mu_1 + 0 = \mu_1$, comme annoncé.

Cela dit, soit $\dot{\mu}_n$ l'image de $\mu_n$ par $X_n \hookrightarrow X$. Si $f$ est continue à support compact sur $X$, alors $\mathrm{supp}(f)$ est contenu dans un $X_m$. Par notre choix des $\mu_n$, cela entraîne que la suite des $\dot{\mu}_n(f)$ est une constante $\mu(f)$ à partir d'un certain rang. La forme linéaire positive $\mu$ ainsi définie est la *borne supérieure* [a] des $\dot{\mu}_n$,

$$\mu = \sup_n \dot{\mu}_n.$$

D'autre part, par construction des $\mu_n$ et par [i] on a $\pi(\dot{\mu}_n) = \chi_n.\nu$, où $\chi_n$ est la fonction caractéristique de $Y_n$. En particulier l'ensemble des $\pi(\dot{\mu}_n)$ est borné supérieurement, par $\nu$; donc [g] et [e] assurent que $\pi$ est $\mu$-propre et que

$$\pi(\mu) = \sup_n \pi(\dot{\mu}_n) = \sup_n (\chi_n.\nu) = \chi.\nu,$$

où $\chi = \sup_n \chi_n$. Mais ceci est la fonction caractéristique de $\pi(X)$; donc par hypothèse $\chi - 1$ est localement $\nu$-négligeable, et [d] donne $\chi.\nu = \nu$.     **q.e.d.**

### *Références*

# Le théorème du petit groupe

Nous démontrons ci-dessous la proposition 3.9 (p. 26). C'est bien sûr essentiellement un cas particulier du théorème de Mackey [2, p. 1284]. Mais pas tout à fait, car ce théorème est formulé en termes de mesures sur $\hat{A}$, qui pour $A$ non-abélien n'est pas forcément localement compact, alors que nous avons eu soin de ne parler que de mesures de Radon.

On pourrait bien sûr traduire les hypothèses; mais il est au fond plus simple de traduire la preuve de Blattner [1]. Sous les hypothèses de (3.8), et étant donnée une représentation unitaire continue $U$ de $G$, la proposition à établir s'énonce:

**Proposition.** *Si $U \circ \exp_{|\mathfrak{a}}$ est une représentation concentrée sur $X_{|\mathfrak{a}}$, alors il existe une unique représentation unitaire continue $T$ de $H$ telle que*

(a)   $U = \mathrm{Ind}_H^G T$,         (b)   $T \circ \exp_{|\mathfrak{a}} = e^{i\langle p, \cdot \rangle}\mathbf{1}$.

**Preuve.** Soit $E$ la mesure spectrale de $U \circ \exp_{|\mathfrak{a}}$ (dans $\mathfrak{a}^*$). Elle est uniquement déterminée par la relation (1.29):

$$(*) \qquad E(e^{i\langle \cdot, Z \rangle}) = U(\exp(Z)) \qquad \forall Z \in \mathfrak{a}.$$

Comme par hypothèse elle est concentrée sur $X_{|\mathfrak{a}}$, (1.31b) nous autorise à la considérer comme définie sur $X_{|\mathfrak{a}} = G/H$, muni de la topologie quotient. Appliquant (1.30) aux automorphismes par lesquels $G$ agit dans $\mathfrak{a}$, on obtient par ailleurs $E(f \circ g_{\mathfrak{a}^*}^{-1}) = U(g)E(f)U(g^{-1})$ pour tout $g \in G$. Donc le théorème d'imprimitivité [3] s'applique, et affirme qu'il existe une unique représentation unitaire continue $T$ de $H$, dans un espace $\mathcal{H}_T$, telle que le couple $(U, E)$ soit donné (à une équivalence unitaire près) par

(a)   $U = \mathrm{Ind}_H^G T$,         (b')   $(E(f)\varphi)(g) = f(g(p))\varphi(g)$.

Rappelons que (a) signifie que la représentation $U$ se fait par translations à gauche dans un espace $\mathcal{H}_U$ constitué de classes de functions mesurables $\varphi : G \to \mathcal{H}_T$ telles que $\varphi(gh) = T^\natural(h^{-1})\varphi(g)$ pour tout $h \in H$, où $T^\natural$ est le produit de $T$ par le rapport de fonctions modules $(\Delta_G/\Delta_H)^{1/2}$, qui vaut 1 dans $A$ [2, p. 243].

Reste à voir que (b) $\Leftrightarrow$ (b'). Or pour $a \in A$ on a $(U(a)\varphi)(g) = \varphi(a^{-1}g) = \varphi(gg^{-1}a^{-1}g) = T(g^{-1}ag)\varphi(g)$. Substituant ceci dans $(*)$ on conclut que (b') est équivalent à

$$\text{(b'')} \qquad T(\exp(g^{-1}Zg))\varphi(g) = e^{i\langle g(p), Z\rangle}\varphi(g) \qquad \text{p. p.}$$

($\forall\, \varphi \in \mathcal{H}_U$, $\forall\, Z \in \mathfrak{a}$). Cela dit, (b) entraîne clairement (b''); pour voir la réciproque, on applique (b'') à des fonctions continues de la forme

$$\varphi(g) = \int_H \phi(gh)\, T^\natural(h)\xi\, dh,$$

où $\xi \in \mathcal{H}_T$ et où $\phi$ est une fonction continue à support compact sur $G$. (Ces $\varphi$ sont dans $\mathcal{H}_U$.) Faisant $g = e$ dans la formule obtenue, et faisant tendre $\phi(h)dh$ vers la mesure de Dirac, on obtient $T(\exp(Z))\xi = e^{i\langle p, Z\rangle}\xi$. **q.e.d.**

### *Références*

## Systèmes élémentaires galiléens

Nous rappelons ci-dessous les définitions nécessaires à la compréhension de l'exemple 3.13 (p. 27). Le *groupe de Galilée*, $\dot G$, est constitué des matrices de la forme
$$\dot g = \begin{pmatrix} A & \boldsymbol{b} & \boldsymbol{c} \\ 0 & 1 & e \\ 0 & 0 & 1 \end{pmatrix}, \qquad \begin{array}{l} A \in \mathbf{SO}(3), \\ \boldsymbol{b}, \boldsymbol{c} \in \mathbf{R}^3, \quad e \in \mathbf{R}. \end{array}$$

Ce groupe agit sur l'espace-temps par $\dot g(\boldsymbol{r}, t) = (A\boldsymbol{r}+\boldsymbol{b}t+\boldsymbol{c}, t+e)$. D'après les §§ 1.1B-C, ses variétés symplectiques homogènes sont, à d'éventuels revêtements près, les orbites coadjointes de son extension centrale universelle. Celle-ci a été calculée par Bargmann [1]: c'est le groupe $G$ des matrices

$$g = \begin{pmatrix} 1 & \overline{\boldsymbol{b}}A & \frac{1}{2}\|\boldsymbol{b}\|^2 & f \\ 0 & A & \boldsymbol{b} & \boldsymbol{c} \\ 0 & 0 & 1 & e \\ 0 & 0 & 0 & 1 \end{pmatrix}$$

qui se projettent dans $\dot G$ en oubliant la première ligne et la première colonne; on a noté $\overline{\boldsymbol{b}} = \langle \boldsymbol{b}, \cdot \rangle$ ('produit scalaire par $\boldsymbol{b}$'), et $f \in \mathbf{R}$. Son algèbre de Lie $\mathfrak{g}$ est donc constituée de matrices de la forme

$$Z = \begin{pmatrix} 0 & \overline{\boldsymbol{\beta}} & 0 & \varphi \\ 0 & j(\boldsymbol{\alpha}) & \boldsymbol{\beta} & \boldsymbol{\gamma} \\ 0 & 0 & 0 & \varepsilon \\ 0 & 0 & 0 & 0 \end{pmatrix}$$

où $j(\boldsymbol{\alpha}) = \boldsymbol{\alpha} \times \cdot$ ('produit vectoriel par $\boldsymbol{\alpha}$'); et on peut voir $\mathfrak{g}^*$ comme un espace de 5-tuples $x = (\boldsymbol{\ell}, \boldsymbol{g}, \boldsymbol{p}, E, m)$, mis en dualité avec $\mathfrak{g}$ par

$$\langle x, Z \rangle = \langle \boldsymbol{\ell}, \boldsymbol{\alpha} \rangle - \langle \boldsymbol{g}, \boldsymbol{\beta} \rangle + \langle \boldsymbol{p}, \boldsymbol{\gamma} \rangle - E\varepsilon - m\varphi.$$

Dans ces variables, l'action coadjointe s'écrit: $g(x) = x^\star$, avec

$$\begin{cases} \boldsymbol{\ell}^\star = A\boldsymbol{\ell} + A\boldsymbol{g}{\times}\boldsymbol{b} + \boldsymbol{c}{\times}A\boldsymbol{p} + m\boldsymbol{c}{\times}\boldsymbol{b}, \\ \boldsymbol{g}^\star = A(\boldsymbol{g} - \boldsymbol{p}e) + m(\boldsymbol{c} - \boldsymbol{b}e), \\ \boldsymbol{p}^\star = A\boldsymbol{p} + m\boldsymbol{b}, \\ E^\star = E + \langle \boldsymbol{b}, A\boldsymbol{p}\rangle + \frac{1}{2}m\|\boldsymbol{b}\|^2, \\ m^\star = m. \end{cases}$$

Un fait mentionné en (3.13) est ici en évidence: à $m$ fixé, les déplacements euclidiens $(A, \boldsymbol{c})$ agissent sur le 'moment de boost' $\boldsymbol{g}$ par $\boldsymbol{g}^\star = A\boldsymbol{g} + m\boldsymbol{c}$. Cela donne à $\boldsymbol{g}$ (ou plutôt $\boldsymbol{g}/m$) la géométrie affine d'une position lorsque $m$ est non nul, mais pas autrement.

Pour classer les orbites, nous appliquons la proposition 3.14 à l'idéal abélien $\mathfrak{a}$ décrit par $(\boldsymbol{\gamma}, \varepsilon, \varphi)$. Comme le montrent les trois dernières lignes ci-dessus, $G$ agit sur $x_{|\mathfrak{a}} = (\boldsymbol{p}, E, m)$ avec trois types d'orbites selon que $m$, et ensuite $\|\boldsymbol{p}\|$, est nul ou pas. Ceci donne lieu à trois possibilités dans (3.14):

(a) $X_{|\mathfrak{a}}$ est un 3-paraboloïde, et $Y$ une orbite de $\mathbf{SO}(3)$,
(b) $X_{|\mathfrak{a}}$ est un 3-cylindre $\mathbf{R} \times S_2$, et $Y$ une orbite de $\mathbf{E}(2)$,
(c) $X_{|\mathfrak{a}}$ est un point, et $Y (= X)$ une orbite de $\mathbf{E}(3)$.

Les orbites (c) sont simplement celles du groupe d'Euclide (qui intervient ici en tant que groupe de Galilée homogène: $\{\dot{g} \in \dot{G} : \boldsymbol{c} = e = 0\}$). Pour décrire les autres, notons $\boldsymbol{e}_1, \boldsymbol{e}_2, \boldsymbol{e}_3$ la base canonique de $\mathbf{R}^3$.

**Cas a: Particules massives.** Étant donnés 3 nombres $s, c, m$ avec $s \geq 0$ et $m \neq 0$, l'orbite
$$X = G(s\boldsymbol{e}_3, \boldsymbol{0}, \boldsymbol{0}, c, m)$$

décrit une particule de masse $m$, spin $s$ et énergie interne $c$. Lorsque $s > 0$, elle s'identifie à $\mathbf{R}^6 \times S_2$ au moyen de l'application $\Phi$:

$$\Phi(\boldsymbol{r}, \boldsymbol{p}, \boldsymbol{u}) = \left(\boldsymbol{r}{\times}\boldsymbol{p} + s\boldsymbol{u},\ m\boldsymbol{r},\ \boldsymbol{p},\ \frac{\|\boldsymbol{p}\|^2}{2m} + c,\ m\right),$$

où $\boldsymbol{r}, \boldsymbol{p} \in \mathbf{R}^3$ et $\boldsymbol{u} \in S_2$. Sa structure symplectique est la somme de $d\boldsymbol{p} \wedge d\boldsymbol{r}$ sur $\mathbf{R}^6$, et $s$ fois la 2-forme d'aire $\Omega$ sur $S_2$. Même chose lorsque $s = 0$, sauf que la sphère disparaît.

**Cas b$_1$ : Particules de masse nulle.** Étant donnés 2 nombres $s$ et $k > 0$, l'orbite
$$X = G(s\boldsymbol{e}_3, \boldsymbol{0}, k\boldsymbol{e}_3, 0, 0)$$
décrit une particule de masse nulle de couleur $k$, spin $|s|$ et hélicité signe($s$). Elle s'identifie à $TS_2 \times \mathbf{R}^2$ au moyen de l'application $\Phi$ :
$$\Phi(\boldsymbol{r}, \boldsymbol{u}, t, E) = \Big(\boldsymbol{r}\times\boldsymbol{p} + s\boldsymbol{u},\ -\boldsymbol{p}t,\ \boldsymbol{p},\ E,\ 0\Big),$$
où $\boldsymbol{r} \in T_{\boldsymbol{u}}S_2$ et $\boldsymbol{p} = k\boldsymbol{u}$. Sa 2-forme est $d\boldsymbol{p}\wedge d\boldsymbol{r} - dE\wedge dt + s\Omega$, et $G$ agit sur elle via son action spatio-temporelle sur la droite $(\boldsymbol{r} + \mathbf{R}\boldsymbol{u}, t)$ : on peut interpréter le point ci-dessus comme le mouvement d'une particule qui parcourrait instantanément la droite $\boldsymbol{r} + \mathbf{R}\boldsymbol{u}$, à la date $t$, avec une énergie $E$.

**Cas b$_2$ : Autres orbites de masse nulle.** Étant donnés 2 nombres positifs $h$ et $k$, l'orbite
$$X = G(\boldsymbol{0}, h\boldsymbol{e}_1, k\boldsymbol{e}_2, 0, 0)$$
ne décrit aucune particule connue. Elle s'identifie à $T\mathbf{SO}(3) \times \mathbf{R}^2$ au moyen de $\Phi$ :
$$\Phi(L, U, t, E) = \Big(\boldsymbol{\ell},\ \boldsymbol{n} - \boldsymbol{p}t,\ \boldsymbol{p},\ E,\ 0\Big),$$
où $L = j(\boldsymbol{\ell})U \in T_U\mathbf{SO}(3)$ et $(\boldsymbol{n}, \boldsymbol{p}) = (hU\boldsymbol{e}_1, kU\boldsymbol{e}_2)$. Sa 2-forme est la dérivée extérieure de $-\frac{1}{2}\mathrm{Tr}(LdU) - Edt$.

**Remarques.** Les détails de l'exemple 3.13 se vérifient facilement sur ces formules. Notons qu'on aurait aussi pu appliquer la théorie de Mackey symplectique à l'idéal Heisenberg décrit par $(\boldsymbol{\beta}, \boldsymbol{\gamma}, \varphi)$; comme dans [2], les calculs s'en trouvent même simplifiés.

### *Références*

# SUBSETS OF $R^n$ WHICH BECOME DENSE IN ANY COMPACT GROUP

## FRANÇOIS ZIEGLER

**Abstract**

A polynomial curve or variety, contained in no proper affine subspace, becomes dense in the homomorphic image of $R^n$ in any compact group.

## 1. Introduction

Let $M$ be a subset of $R^n$ and let $E$ denote the smallest affine subspace containing $M$. We say that $M$ is *Bohr dense in* $E$ if $\beta(M)$ is dense in $\beta(E)$ whenever $\beta$ is a continuous morphism of $R^n$ into some compact group, $T$. We say that $M$ is a *polynomial variety* if it is the image of a polynomial map $P : R^m \to R^n$. We propose to prove the following:

**Theorem.** *Let $M \subset R^n$ be a polynomial variety. Then $M$ is Bohr dense in its affine hull $E$.*

Thus a cone is Bohr dense in space, and so is a parabola in the plane; in contrast we observe that a hyperbola is not. Functions of the form $f \circ \beta$, with $\beta$ as above and $f$ a continuous function on $T$, are called *almost periodic*; so we may restate the theorem as:

**Corollary.** *Let $M \subset R^n$ be a polynomial variety, contained in no proper affine subspace. Then every almost periodic function on $R^n$ is determined by its restriction to $M$.*

## 2. Proof of the Theorem

Let all notation be as in the introduction. Translating everything so that $M$ contains the origin, we may assume that $E$ is a vector subspace of $R^n$. Replacing $R^n$ by this subspace and $T$ by an abelian subgroup if necessary, we may also assume that $\overline{\beta(E)} = T$. We must then prove that $\overline{\beta(M)} = T$.

The idea is to show that the image of a "renormalized" Lebesgue measure under the given maps $R^m \to E \to T$ coincides with Haar measure $\eta$

386 FRANÇOIS ZIEGLER

on $T$. More precisely let $\lambda$ denote Lebesgue measure on the cube $[0,1]^m$ and let $\mu_a$ be the image of $\lambda$ under the maps

$$[0,1]^m \xrightarrow{a} \mathbf{R}^m \xrightarrow{P} E \xrightarrow{\beta} T,$$

where $a$ stands for dilation by a factor $a \in \mathbf{R}$; then we shall prove that for every continuous function $f$ on $T$ one has

$$(*) \qquad \lim_{a \to \infty} \mu_a(f) = \eta(f).$$

The theorem follows: for if $f \geq 0$ vanishes on $\overline{\beta(M)}$ then so does the left-hand side and hence $\eta(f)$, which forces $f$ to vanish everywhere.

Since linear combinations of characters are uniformly dense in the continuous functions on $T$, it is enough to prove $(*)$ when $f$ is a character, in which case $\eta(f)$ is 1 or 0 according as $f \equiv 1$ or not. Since $(*)$ is clear when $f \equiv 1$, it remains to show that $\mu_a(f) = \lambda(f \circ \beta \circ P \circ a) \to 0$ as $a \to \infty$ whenever $f$ is a nontrivial character of $T$. But then $f \circ \beta$ is a nontrivial character of $E$ (because $\beta(E)$ is dense in $T$), so

$$(f \circ \beta \circ P)(x) = e^{i\langle \varphi, P(x) \rangle}$$

for a nonzero linear form $\varphi$ on $E$. Moreover the polynomial $p(x) = \langle \varphi, P(x) \rangle$ is not constant on $\mathbf{R}^m$, for $M$ is contained in no proper affine subspace of $E$. So matters are reduced to the following:

**Lemma.** *If $p$ is a nonconstant polynomial on $\mathbf{R}^m$ then*

$$\lim_{a \to \infty} \int_{[0,1]^m} e^{ip(ax)} dx = 0.$$

*Proof.* Since $p$ is not a constant, it has degree $k \geq 1$ in at least one of the variables, say $x_i = t$. Writing $y$ for the remaining $x_j$ and $p(x) = p_y(t)$, our integral $I_a$ becomes the integral over $[0,1]^{m-1}$ of

$$I_a(y) = \int_0^1 e^{ip_{ay}(at)} dt = \frac{1}{a} \int_0^a e^{ip_{ay}(t)} dt.$$

Now consider the coefficient $c(y)$ of $t^k$ in $p_y(t)$: being a nonzero polynomial in $y$, it is nonzero for all $y$ in a conull set $Y \subset [0,1]^{m-1}$. Likewise $c(ay)$, for fixed $y \in Y$, is a nonzero polynomial in $a$ and therefore nonzero for all $a$ in a cofinite set $A_y \subset \mathbf{R}$. For fixed $(y,a)$ in $Y \times A_y$ we conclude that the $k$th derivative $k!c(ay)$ of $p_{ay}$ is bounded away from zero. By the generalized van der Corput lemma [1, p. 1258] this implies that

$$\left| \int_u^v e^{ip_{ay}(t)} dt \right| \leq \frac{2^{k+1}}{|k!c(ay)|^{1/k}} \quad \forall u, v \in \mathbf{R}.$$

Taking $[u, v] = [0, a]$ and letting $a \to \infty$, it follows that $I_a(y) \to 0$ for all $y \in Y$, whence $I_a \to 0$ by dominated convergence. This completes the proof.

### Acknowledgements

This work was supported by a grant of the Société Académique Vaudoise. I am also grateful to J. M. Souriau, J. Elhadad and B. Coupet for their help.

UNIVERSITÉ D'AIX-MARSEILLE II ET CNRS-CPT LUMINY, CASE 907 F-13288 MARSEILLE CEDEX 09, FRANCE
*E-mail address*: ziegler@cptvax.in2p3.fr

# ON THE KOSTANT CONVEXITY THEOREM

FRANÇOIS ZIEGLER

(Communicated by Jonathan M. Rosenberg)

ABSTRACT. A quick proof that the coadjoint orbits of a compact connected Lie group project onto convex polytopes in the dual of a Cartan subalgebra.

## 1. Introduction

Let $G$ be a compact connected Lie group, $T$ a maximal torus of $G$, $\mathfrak{g}$ and $\mathfrak{t}$ their Lie algebras and $\pi\colon \mathfrak{g}^* \to \mathfrak{t}^*$ the natural projection. As usual we identify $\mathfrak{t}^*$ with the subspace of all $T$-fixed points in $\mathfrak{g}^*$. Then every coadjoint orbit $X$ of $G$ intersects $\mathfrak{t}^*$ in a Weyl group orbit $\Omega_X$ [4], and in this setting B. Kostant [9] has proved:

**1.1. Theorem.** $\pi(X)$ *is the convex hull of* $\Omega_X$.

Alternative proofs and generalizations have appeared in [2, 5, 7, 8]; see the monograph [3]. Our purpose here is to show that representation theory and the projective embeddings of Borel-Weil-Tits [10, 12] allow for an elementary proof of Theorem 1.1, bypassing the Morse theoretic or asymptotic arguments of *loc.cit.*

## 2. Projective embeddings

If $\Omega_X$ lies in the weight lattice $\Lambda = \{w \in \mathfrak{t}^* : e^{i\langle w, Z\rangle} = 1 \ \forall\ Z \in \ker(\exp|\mathfrak{t})\}$, we say that $X$ is *integral*; then $\Omega_X$ contains the highest weight $w_0$ of a unique irreducible unitary $G$-module $V$ [1]. The corresponding projective space $\mathbf{P}(V)$, regarded as the manifold of all rank one hermitian projectors $\mathbf{p}$ in $V$, carries canonical complex and symplectic structures $J$ and $\sigma$, defined on tangent vectors $\delta\mathbf{p}, \delta'\mathbf{p} \in T_{\mathbf{p}}\mathbf{P}(V)$ by

$$J\delta\mathbf{p} = \frac{1}{i}[\mathbf{p}, \delta\mathbf{p}], \qquad \sigma(\delta\mathbf{p}, \delta'\mathbf{p}) = \operatorname{Tr}(\delta'\mathbf{p} J\delta\mathbf{p}).$$

Writing $\mathbf{E}_0$ for the eigenprojector associated to $w_0$, we know from [10, 12] that the $G$-orbit of $\mathbf{E}_0$ is a *complex* submanifold, $\mathbf{X}$, of $\mathbf{P}(V)$. In particular $\mathbf{X}$ is homogeneous symplectic, with momentum map $\Phi\colon \mathbf{X} \to \mathfrak{g}^*$ readily computed

---

as

$$(*) \qquad \langle \Phi(\mathbf{x}), Z \rangle = \frac{1}{i}\text{Tr}(\mathbf{x}\mathbf{Z}),$$

where $Z \mapsto \mathbf{Z}$ is the differentiated representation of $\mathfrak{g}$ on $V$. By Kirillov-Kostant-Souriau [11] $\Phi$ covers a coadjoint orbit of $G$, namely $X$ since $\Phi(\mathbf{E}_0) = w_0$. But $X$ is simply connected [10], so $\Phi$ is actually a diffeomorphism $\mathbf{X} \leftrightarrow X$.

### 3. Proof of the theorem

If Theorem 1.1 holds when $\Omega_X$ lies in $\Lambda$, it follows also for $\Omega_X$ in $\mathbb{R}\Lambda$ by rescaling, and then for the general $\Omega_X$ in $\mathfrak{t}^* = \overline{\mathbb{R}\Lambda}$ by a straightforward continuity argument. So it is enough to prove Theorem 1.1 when $X$ is integral.

Let, then, $\mathbf{X} \subset \mathbf{P}(V)$ be as above; also let $\Delta \subset \mathfrak{t}^*$ be the weight diagram of $V$, so that we have

$$\frac{1}{i}\mathbf{Z} = \sum_{w \in \Delta} \langle w, Z \rangle \mathbf{E}_w \quad \forall Z \in \mathfrak{t},$$

where $\mathbf{E}_w$ denotes the eigenprojector belonging to $w \in \Delta$. Substituting this in $(*)$ exhibits $\pi(\Phi(\mathbf{x}))$ as a convex combination of elements of $\Delta$; since $\Delta$ lies in the convex hull of $\Omega_X$ [1] so does, therefore, $\pi(X)$.

For the converse inclusion we use a variational method inspired from [6]. Let $\{w_j\}$ be an enumeration of $\Omega_X$ and write $\mathbf{E}_j$ for the projectors $\Phi^{-1}(w_j) = \mathbf{E}_{w_j}$. Given a convex combination $\sum_j \mu_j w_j$ of the $w_j$, we maximize the nonnegative function

$$\rho(\mathbf{x}) = \prod_j \text{Tr}(\mathbf{E}_j \mathbf{x})^{\mu_j}$$

and compute its derivative $D\rho(\mathbf{x})(\delta\mathbf{x})$ in the tangent direction

$$\delta\mathbf{x} = J[\mathbf{Z}, \mathbf{x}], \quad Z \in \mathfrak{t}.$$

Since $\mathbf{X}$ is compact $\rho$ does attain its maximum, which is positive: if $\rho$ vanished identically, so would the product of the real analytic functions $\rho_j(\mathbf{x}) = \text{Tr}(\mathbf{E}_j\mathbf{x})$ and hence also one of the $\rho_j$, whereas $\rho_j(\mathbf{E}_j) = 1$. Now we have

$$D\rho_j(\mathbf{x})(\delta\mathbf{x}) = \text{Tr}(\mathbf{E}_j\delta\mathbf{x}) = \frac{1}{i}\text{Tr}(\mathbf{E}_j[2\mathbf{x}\text{Tr}(\mathbf{x}\mathbf{Z}) - \mathbf{Z}\mathbf{x} - \mathbf{x}\mathbf{Z}])$$

$$= 2\rho_j(\mathbf{x})\langle \Phi(\mathbf{x}) - w_j, Z \rangle,$$

whence

$$D\rho(\mathbf{x})(\delta\mathbf{x}) = 2\rho(\mathbf{x})\left\langle \Phi(\mathbf{x}) - \sum_j \mu_j w_j, Z \right\rangle = 0 \quad \forall Z \in \mathfrak{t}$$

at the maximum. Thus $\Phi(\mathbf{x})$ projects to the given convex combination, and our proof is complete.

### Acknowledgments

This research was supported by the Fonds National Suisse de la Recherche Scientifique and the Sunburst-Fonds der ETH. I would also like to thank Ch. Duval, J. Elhadad, and J. M. Souriau for encouragement and help.[1]

---

[1] *Note added in proof.* Michèle Vergne has kindly pointed out that V. G. Kac & D. H. Peterson [13] also used projective embeddings (but not the short variational argument above) to prove Theorem 1.1.

UNIVERSITÉ D'AIX-MARSEILLE II ET CNRS-CPT LUMINY, CASE 907, F-13288 MARSEILLE CEDEX 09, FRANCE

*E-mail address*: ZIEGLER@CPTVAX.IN2P3.FR

# Abstract

This Thesis consists mainly of two papers.

**1) Quantum representations and the orbit method.** A major success of group representers has been to fit much of their subject into the following scheme: *unitary representations of a Lie group correspond to the symplectic manifolds on which the group acts.* Yet, a statement like "this representation corresponds to (or "quantizes") that symplectic manifold" is usually true by definition, rather than a theorem. Can we get these definitions from some principle—which would generalize the hypothesis of the Stone-von Neumann theorem?

J. M. Souriau has proposed to this effect a principle whose *spectral* nature we first demonstrate (it essentially requires that commuting observables have their quantum spectrum concentrated in their classical range, suitably compactified), and whose relation with the traditional orbit method we then describe.

When the group is compact, we show that our principle effectively selects the expected representation within sections of the line bundle over the orbit. When the group is noncompact, on the other hand, we find many unexpected representations; but we show how the principle can be refined so as to, either, recover the traditional theory for exponential groups, or in another direction, characterize some new discontinuous representations in which states can be "localized" on lagrangian submanifolds of the orbit.

**2) Théorie de Mackey symplectique.** Kazhdan, Kostant and Sternberg introduced in 1978 a construction intended to play, in symplectic geometry, le role played by unitary induction in representation theory. In this second paper, we show that a detailed analogue of all of Mackey's theory (imprimitivity theorem, normal subgroup analysis) can be developed in this purely geometrical framework. In particular, every induced symplectic manifold comes equipped with a distinguished abelian group of automorphisms (system of imprimitivity).

Our main observation then is that the spectral principle of the first paper can be imposed to the generators of this group, and that as a result, the representations attached to an induced symplectic manifold are necessarily induced from the same subgroup.

Five technical appendices complete the Thesis.

# Résumé

Cette thèse se compose principalement de deux articles.

**1) Quantum representations and the orbit method.** Un des grands succès de la théorie des représentations est de s'être presqu'entièrement coulée dans le schéma suivant: *les représentations unitaires d'un groupe de Lie correspondent aux variétés symplectiques sur lesquelles ce groupe agit.* Pourtant, l'assertion que telle représentation correspond à (ou "quantifie") telle variété est en général une définition, plutôt qu'un théorème. Peut-on faire résulter ces définitions de quelque principe—qui généraliserait l'hypothèse du théorème de Stone-von Neumann?

J.-M. Souriau a proposé à cet effet un principe dont nous montrons d'abord la nature *spectrale* (il exige en substance que le spectre quantique d'observables qui commutent soit concentré dans leur ensemble de valeurs classique, convenablement compactifié), et dont nous déterminons ensuite la relation avec la méthode des orbites traditionnelle.

Lorsque le groupe est compact, la condition sélectionne bien la représentation attendue dans les sections du fibré en droites au-dessus de l'orbite. Lorsque le groupe est non compact, nous trouvons par contre de nombreuses représentations inattendues ; mais nous montrons comment on peut renforcer la condition de Souriau, soit de manière à retrouver la théorie traditionnelle pour les groupes exponentiels, soit de manière à caractériser des représentations nouvelles, discontinues, qui admettent des états "localisés" sur des sous-variétés lagrangiennes de l'orbite.

**2) Théorie de Mackey symplectique.** Kazhdan, Kostant and Sternberg ont introduit en 1978 une construction destinée à jouer, en géométrie symplectique, le rôle que tient l'induction unitaire en théorie des représentations. Dans ce second article, nous détaillons l'analogue de la théorie de Mackey (théorème d'imprimitivité, "normal subgroup analysis") dans ce cadre purement géométrique. En particulier, toute variété symplectique induite se trouve munie d'un groupe d'automorphismes abélien privilégié.

Appliquant à ce groupe le principe spectral étudié dans l'article précédent, nous montrons alors qu'à une variété induite doivent nécessairement correspondre des représentations induites.

Cinq annexes techniques complètent l'ensemble.